\documentclass[10pt]{article}
\usepackage{amsmath, amsthm, amsfonts, amssymb}
\usepackage{epsfig}
\usepackage[english]{babel}

\usepackage{enumerate}
\usepackage{bm}
\usepackage{bbm}
\usepackage{color}
\usepackage{graphicx}
\usepackage{caption}
\usepackage{tikz}

\usepackage{verbatim} 

\newtheorem{thm1}{Theorem}

\newtheorem{thm4}{Lemma}
\newtheorem{thm5}{Proposition}
\newtheorem{thm6}{Remark}
\newtheorem{thm7}{Definition}

\begin{document}

\title{Convergence of Diffusion Generated Motion to\\
Motion by Mean Curvature}
\author{\large
Drew Swartz\\
\normalsize
IRI, Chicago, IL, USA\\
\large Nung Kwan Yip\\
\normalsize
Department of Mathematics, Purdue University, USA
}

\maketitle

\newcommand{\nc}{\newcommand}

\nc{\Cal}[1]{{\mathcal {#1}}}
\nc{\Bf}[1]{{\bf {#1}}}
\nc{\Em}[1]{{\em {#1}}}
\nc{\Rm}[1]{{\rm {#1}}}
\nc{\nonu}{\nonumber}
\nc{\Headline}[1]{\noindent{\bf {#1}: }}

\nc{\EqnRef}[1]{(\ref{#1})}
\nc{\DefRef}[1]{\Bf{Definition \ref{#1}}}
\nc{\LemRef}[1]{\Bf{Lemma \ref{#1}}}
\nc{\ProRef}[1]{\Bf{Proposition \ref{#1}}}
\nc{\ThmRef}[1]{\Bf{Theorem \ref{#1}}}
\nc{\CorRef}[1]{\Bf{Corollary \ref{#1}}}
\nc{\SecRef}[1]{\Bf{Section \ref{#1}}}
\nc{\RemRef}[1]{\Bf{Remark \ref{#1}}}
\nc{\ChpRef}[1]{\Bf{Chapter \ref{#1}}}

\nc{\defref}[1]{Definition \ref{#1}}
\nc{\lemref}[1]{Lemma \ref{#1}}
\nc{\proref}[1]{Proposition \ref{#1}}
\nc{\thmref}[1]{Theorem \ref{#1}}
\nc{\corref}[1]{Corollary \ref{#1}}
\nc{\secref}[1]{Section \ref{#1}}
\nc{\remref}[1]{Remark \ref{#1}}
\nc{\chpRef}[1]{Chapter \ref{#1}}

\nc{\Proof}{\begin{proof}}
\nc{\EndProof}{\end{proof}}

\nc{\apriori}{{\it a priori }}
\nc{\Apriori}{{\it A priori }}
\nc{\Holder}{H\"{o}lder }

\nc{\Sty}{\displaystyle}
\nc{\MathSty}[1]{$\Sty{#1}$}
\def\Beqn#1\Eeqn{\begin{equation}#1\end{equation}}

\def\upchi{\raise2pt\hbox{$\chi$}}
\def\upnu{\raise0pt\hbox{$\nu$}}
\def\upstroke{\raise3pt\hbox{$|$}}
\nc{\Bdry}{\partial}
\nc{\Abs}[1]{\left|{#1}\right|}
\nc{\abs}[1]{|{#1}|}
\nc{\Ave}[1]{\bar{#1}}
\nc{\CurBrac}[1]{\left\{{#1}\right\}}
\nc{\curBrac}[1]{\{{#1}\}}
\nc{\Brac}[1]{\left({#1}\right)}
\nc{\brac}[1]{({#1})}
\nc{\SqrBrac}[1]{\left[{#1}\right]}
\nc{\sqrbrac}[1]{[{#1}]}
\nc{\BigBrac}[1]{\Big({#1}\Big)}
\nc{\bigBrac}[1]{\big({#1}\big)}
\nc{\BigCurBrac}[1]{\Big\{{#1}\Big\}}
\nc{\bigCurBrac}[1]{\big\{{#1}\big\}}
\nc{\BigSqrBrac}[1]{\Big[{#1}\Big]}
\nc{\bigSqrBrac}[1]{\big[{#1}\big]}
\nc{\BigAbs}[1]{\Big|{#1}\Big|}
\nc{\bigAbs}[1]{\big|{#1}\big|}
\nc{\converge}{\rightarrow}
\nc{\Converge}{\longrightarrow}
\nc{\WeakConverge}{\rightharpoonup}
\nc{\Del}{\triangle}
\nc{\Div}{\mbox{\rm{div}}}
\nc{\Equivalent}{\Longleftrightarrow}
\nc{\IntOverR}{\int_{-\infty}^{\infty}}
\nc{\IntOverRp}{\int_{0}^{\infty}}
\nc{\IntOverRm}{\int_{-\infty}^{0}}
\nc{\Imply}{\Longrightarrow}
\nc{\InnProd}[2]{\left\langle{#1},\,{#2}\right\rangle}
\nc{\innprod}[2]{\langle{#1},\,{#2}\rangle}
\nc{\BracInnProd}[2]{\left({#1},\,{#2}\right)}
\nc{\bracInnprod}[2]{({#1},\,{#2})}
\nc{\SqrinnProd}[2]{\left[{#1},\,{#2}\right]}
\nc{\sqrinnprod}[2]{[{#1},\,{#2}]}
\nc{\Lap}{\triangle}
\nc{\Lip}{\mbox{\rm{Lip}}}
\nc{\Lover}[1]{\frac{1}{#1}}
\nc{\Grad}{\nabla}
\nc{\MapTo}{\longrightarrow}
\nc{\Min}{\wedge}
\nc{\Max}{\vee}
\nc{\Norm}[1]{\left\|#1\right\|}
\nc{\norm}[1]{\|#1\|}
\nc{\SingleNorm}[1]{\left|#1\right|}
\nc{\singlenorm}[1]{|#1|}
\nc{\Pair}[2]{\left\langle{#1},\,{#2}\right\rangle}
\nc{\pair}[2]{\langle{#1},\,{#2}\rangle}
\nc{\DD}[1]{\frac{d}{d{#1}}}
\nc{\PDD}[1]{\frac{\partial}{\partial{#1}}}
\nc{\Spt}{\mbox{\rm{spt}}}
\nc{\spt}{\mbox{\rm{spt}}}
\nc{\SuchThat}{\ni}
\nc{\Sum}[2]{\sum_{#1}^{#2}}
\nc{\wtilde}[1]{\widetilde{#1}}
\nc{\Text}[1]{{\mbox{{\rm #1}}}}
\nc{\TextMath}[1]{\mbox{\,\,\,{#1}\,\,\,}}
\nc{\intersect}{\cap}
\nc{\Intersect}{\bigcap}
\nc{\union}{\cup}
\nc{\Union}{\bigcup}
\nc{\ColTwo}[2]
{\left(\begin{array}{c}{#1}\\{#2}\end{array}\right)}
\nc{\ColThree}[3]
{\left(\begin{array}{c}{#1}\\{#2}\\{#3}\end{array}\right)}
\nc{\ColFour}[4]
{\left(\begin{array}{c}{#1}\\{#2}\\{#3}\\{#4}\end{array}\right)}
\nc{\ColFive}[5]
{\left(\begin{array}{c}{#1}\\{#2}\\{#3}\\{#4}\\{#5}\end{array}
\right)}
\nc{\coltwo}[2]
{\begin{array}{c}{#1}\\{#2}\end{array}}
\nc{\colthree}[3]
{\begin{array}{c}{#1}\\{#2}\\{#3}\end{array}}
\nc{\colfour}[4]
{\begin{array}{c}{#1}\\{#2}\\{#3}\\{#4}\end{array}}
\nc{\colfive}[5]
{\begin{array}{c}{#1}\\{#2}\\{#3}\\{#4}\\{#5}\end{array}}

\nc{\Vol}[2]{\Bf{#1}({#2})}

\nc{\dt}{{h}}
\nc{\bx}{{\bar{x}}}
\nc{\bz}{{\bar{z}}}
\nc{\E}{{\mathbf{E}}}
\nc{\R}{{\mathbb{R}}}
\nc{\N}{{\mathbf{N}}}
\nc{\T}{{\mathbb{T}}}
\nc{\bbL}{{\mathbb{L}}}
\nc{\CB}{\Cal{B}}
\nc{\dH}{d_{\Cal{H}}}
\nc{\NORMAL}[1]{\frac{\nabla{#1}}{\Abs{\nabla{#1}}}}
\nc{\CO}{\Cal{O}}

\nc{\purple}[1]{\textcolor{purple}{#1}}
\nc{\red}[1]{\textcolor{red}{#1}}
\nc{\RED}[1]{\textcolor{red}{\bf #1}} 
\nc{\blue}[1]{\textcolor{blue}{#1}}
\nc{\BLUE}[1]{\textcolor{blue}{\bf #1}} 
\nc{\green}[1]{\textcolor{green}{#1}}
\nc{\GREEN}[1]{\textcolor{green}{\bf #1}} 


\begin{abstract}
We provide a new proof of convergence to motion by mean curvature (MMC)
for the Merriman-Bence-Osher (MBO) thresholding algorithm. The proof
is elementary and does not rely on maximum principle for the scheme. 
The strategy is to construct a natural ansatz of the solution and then
estimate the error. The proof thus also provides a convergence rate.
Only some weak integrability assumptions of the heat 
kernel, but not its positivity, is used. Currently the result is proved
in the case when smooth and classical solution of MMC exists.
\end{abstract}

\section{Introduction}

Motion by mean curvature (MMC) is a fundamental geometric motion 
arising in a broad range of scientific disciplines. 
Besides its intrinsic geometric interests, in applications, 
it arises naturally in modeling the 
evolution of interfaces such as grain boundaries which are subject to the
effect of surface tension.
It also appears in various aspects of image de-noising
algorithms.
Mathematically, MMC describes the evolution of a manifold with its normal 
velocity $V_N$ equal to its mean curvature, i.e.
\Beqn\label{MMC}
V_N = H = \Lover{n}\sum_{i = 1}^n \kappa_i,
\Eeqn
where the $\kappa_i$'s are the principal curvatures of the manifold 
at a point.
This evolution can also be thought of as the geometric analogue
of the heat equation.  Specifically, given an initial $n$-dimensional
embedded manifold $M_0$ in $\mathbb{R}^{n+1}$,
the time dependent manifolds $M_t = F( M_0 , t)$ with 
$F:M_0\times\R_+\MapTo \R^{n+1}$ solves the MMC equation \eqref{MMC} 
with initial data $M_0$ if $F$ satisfies the following equation,
\Beqn\label{MMC.eqn}
\partial_t F = \Delta_{M_t} F,\,\,\,\,\,\,
F( \cdot, 0) = F_0,
\Eeqn
where $\Delta_{M_t}$ denotes the Laplace-Beltrami operator on the manifold 
$M_t$. 
Alternatively, MMC can be interpreted as the $L^2$-gradient flow, 
or steepest descent for the area functional of a manifold.

It is known that this geometric flow can lead to 
singularity formation and topological changes for the underlying
evolving manifold.
Thus it is desirable to use mathematical formulations which can
handle such events. One method is the \emph{phase-field} approach.
The underlying equation is typically given by the 
\emph{Allen-Cahn equation}
\[
\frac{\partial u}{\partial t} = 
\triangle u - \frac{1}{\epsilon^2}W'(u),
\]
where $W(u) = (1-u^2)^2$ is the double-well potential (with $1$ and 
$-1$ being the two minima of $W$). 
As $u$ evolves under this equation, the domain will separate into two 
regions/phases where $u$ is approximately 
equal to $1$ in one region and $-1$ in the other.
Between these two regions, $u$ will have a diffuse transition layer of thickness $\mathcal{O}(\epsilon)$.
For $\epsilon \ll 1$, this layer will in turn evolve approximately by 
its mean curvature. Convergence of this motion to MMC as 
$\epsilon\Converge 0$
has been shown rigorously by several authors using various approaches, 
for example, 
\cite{ilmanen1993convergence, soner1, soner2} (varifold formulation), 
\cite{bronsard1991motion} (energy approach), 
\cite{de1995geometrical} (asymptotic expansion), 
\cite{chen1992generation} (sub- and super-solutions), 
\cite{evans1992phase} (viscosity solution). 

Another formulation for MMC is to make use of a level set function 
$u: \R^{n+1}\times\R_+\MapTo\R$ which solves the following equation,
$$
\frac{\partial u}{\partial t} = |\nabla u|\text{div}\left(
\frac{\nabla u}{|\nabla u|}
\right).
$$
Each $c$-level set of $u$, $\CurBrac{u=c}$, 
or $\Bdry\CurBrac{u \geq c}$, then evolves 
by its mean curvature (in the viscosity sense). 
The theory surrounding this equation has been developed independently by 
Evans-Spruck (\cite{evans1991motion}) and Chen-Giga-Goto 
(\cite{chen1991uniqueness}). 
See \cite{SouganidisBarles} for more recent development, in particular
the relation between the phase-field and level set formulations.

The \emph{diffusion generated motion} for approximating mean 
curvature flow was first proposed by Merriman-Bence-Osher in 
\cite{BMO} which hereby from now on will be called the MBO algorithm 
or scheme. Essentially it is a time-splitting algorithm for the 
Allen-Cahn equation. Its description is as follows.

Let $\dt > 0$ be a small positive number, called the time step.
Given an initial set $\Omega_0\subseteq\R^{n+1}$, 
a new set $\Omega_{\dt}$ and its boundary $\Gamma_{\dt}=\Bdry\Omega_{\dt}$ 
are constructed by the following two simple procedures:
\begin{description}
\item[Step 1 - linear diffusion:]
given a set $\Omega_0\subseteq\R^{n+1}$ and let
\Beqn
U_0(x) = 2 \cdot \mathbbm{1}_{\Omega_0}(x) - 1
= \left\{
\begin{array}{ll}
1 & \text{for}\,\,x\in\Omega_0,\\
-1 & \text{for}\,\,x\notin\Omega_0.
\end{array}
\right.
\Eeqn
Then solve the linear heat equation:
\Beqn
U_t - \Delta U = 0,\,\,\,\text{in $\mathbb{R}^{n+1} \times (0, \dt)$,}
\,\,\,\,\,\,U(x,0)=U_0(x).
\Eeqn

\item[Step 2 - thresholding:]
project $U(\cdot,\dt)$, the solution from \Bf{Step 1}, 
onto $\CurBrac{-1,1}$ by setting
\Beqn
U(x,\dt^+) = \left\{
\begin{array}{ll}
1 & \text{if}\,\,\,U(x,\dt) \geq 0,\\
-1 & \text{if}\,\,\,U(x,\dt) < 0.
\end{array}
\right.
\Eeqn
Then define
\Beqn
\Omega_{\dt} = \CurBrac{x: U(x, \dt^+) \geq 0}
\,\,\,\text{and}\,\,\,
\Gamma_{\dt} = \Bdry\Omega_\dt.
\Eeqn
Note that
$U(x, \dt^+)  = 2\mathbbm{1}_{\Omega_{\dt}}(x) -1.$
\end{description}
The algorithm then repeats \Bf{Steps 1} and \Bf{2}
but using the result from \Bf{Step 2} as the initial data for 
\Bf{Step 1}. The following collection of sets are thus generated:
\Beqn\label{DGM-sequence}
\Big\{
\Omega_{k\dt} = \CurBrac{x: U(x, k\dt^+) \geq 0}, \,\,\,
\Gamma_{k\dt} = \Bdry\Omega_{k\dt}
\Big\}_{k=0,1,2\ldots}.
\Eeqn

It is proved in \cite{Evans} and \cite{BG} that as 
$\dt\Converge 0$, the sequence \eqref{DGM-sequence}
converges to a solution of MMC in the viscosity 
sense. See also \cite{ishii1995generalization, leoni2001convergence, 
ishii1999threshold} for generalizations to incorporate general kernels and 
anisotropy effects. The recent works \cite{EseOtto, elsey2016threshold,
laux2016convergence} have recast thresholding schemes into a variational 
setting.

In this paper, we provide a new convergence proof of the algorithm. 
Our approach is elementary and does not rely on the theory of
viscosity solution which depends very much on 
comparison principle. Furthermore, it provides a convergence rate.
Even though so far it has only been applied to the case 
when smooth and classical solution of MMC exists, 
it has the potential to be used
in the study of thresholding schemes for other geometric motions.
These include fourth order flows such as 
Willmore and surface diffusion flows \cite{EseRuuthTsai}, 
and higher co-dimension mean curvature flows such 
as filament motions in $\mathbb R^3$ \cite{OsherRuuthXin}.  
Furthermore, the works \cite{Chambolle, EseOtto}, 
for anisotropic MMC shows that the convolution kernel used in the
thresholding scheme are necessarily non-positive.
Consequently comparison principle does not hold.
Our proof can offer reasonable generalizations to these settings.

\section{Main Results and Outline of Proof}
Recall that the MBO scheme produces, for each $\dt > 0$, a sequence of
sets and their boundaries \eqref{DGM-sequence}.
Our main result is that the $\Gamma_{k\dt}$ converges to the solution of
MMC as long as the classical solution exists. 
The precise statement of result is given in the following. 
All the definitions used will be elaborated afterwards.
\begin{thm1}\label{MainTheorem}
Let $\Gamma_0$ be a compact, smooth embedded $n$-dimensional manifold in 
$\R^{n+1}$. Then there is a time $T > 0$ for which the classical solution of 
MMC starting from $\Gamma_0$ exists on $[0,T]$ such that the sequence
$\CurBrac{\Gamma_{k\dt}}_{0 \leq k \leq \lfloor \frac{T}{\dt} 
\rfloor}$ remains an embedded manifold and
as $\dt \Converge 0$, converges to the solution of MMC
in the Hausdorff distance and also in the Bounded-Variation (BV)-sense.
The time $T$ depends on the initial manifold $\Gamma_0$ but not on the
time step $h$.
\end{thm1}

Next we give the definitions of Hausdorff distance, BV-convergence and
some remarks about the theorem.

\begin{thm7}[Hausdorff Distance $d_{\Cal{H}}(\cdot)$]
Let $A$ and $B$ be two subsets of $\R^{n+1}$. Then the Hausdorff distance
between $A$ and $B$ is defined as:
\Beqn\label{Hausdorff}
d_\Cal{H}(A,B) = \max\CurBrac{
\sup_{y\in B\backslash A}\Text{dist}(y,A),
\sup_{x\in A\backslash B}\Text{dist}(x,B)
}
\Eeqn
where
$\Text{dist}(y,A) = \inf_{x\in A}|y-x|$ and likewise,
$\Text{dist}(x,B) = \inf_{y\in B}|x-y|$.
\end{thm7}
It is known that $d_{\Cal{H}}$ gives rise to a complete 
metric space. In addition, for 
$\Omega_1$, $\Gamma_1=\Bdry\Omega_1$ and
$\Omega_2$, $\Gamma_2=\Bdry\Omega_2$, it holds that
$\dH(\Omega_1, \Omega_2) = \dH(\Gamma_1, \Gamma_2).$
Hence it does not matter if we are using $\Omega$ or $\Gamma$ to measure
the Hausdorff distance.

To formulate the notion of convergence, we define
\Beqn\label{approx1}
{\Gamma}_{\dt}(t) = M(\Gamma_{k\dt}, t-k\dt),\,\,\,
\,\,\,\,\,\,
\text{for each}\,\,\,k \geq 1,\,\,\,\text{and}\,\,\,k\dt \leq t < (k+1)\dt
\Eeqn
where $M(\Gamma_{k\dt},t-k\dt)$ is the solution at time $t-k\dt$
of the classical MMC \eqref{MMC.eqn} with initial data $\Gamma_{k\dt}$.
Note that due to the thresholding step, $\Gamma_\dt(t)$ is discontinuous in 
time, in particular, at each $t=k\dt$: 
$\Gamma_{\dt}(k\dt^-) \neq \Gamma_{\dt}(k\dt^+)$.
Let $\Omega_\dt(t)$ be the interior of $\Gamma_{\dt}(t)$, 
i.e. the bounded subset of $\R^{n+1}$ such that
$\Gamma_{\dt}(t) = \Bdry\Omega_{\dt}(t)$.
Then we define
\Beqn\label{approx2}
\chi_{\dt}: \R^{n+1}\times\R_+\MapTo\CurBrac{0,1}:
\,\,\,\,\,
\chi_{\dt}(x,t) = \mathbbm{1}_{{\Omega}_{\dt}(t)}(x).
\Eeqn
(The above definition is to facilitate a more efficient usage of
integration by parts formula for MMC,
in particular \eqref{MMC.nonlinear.1}-\eqref{MMC.nonlinear.2}. 
But from the perspective of 
understanding the convergence statement, it is essentially the same 
as if we had used the ``piece-wise-constant'' definition:
$\Gamma_{\dt}(t) = \Gamma_{k\dt}$,
$\Omega_{\dt}(t) = \Omega_{k\dt}$ for $k\dt < t < (k+1)\dt$.)

\begin{thm7}[BV-Convergence to MMC]\label{luckhaus} \cite{luckhaus}
There is a $0 < T < \infty$ and a
\Beqn
\chi^*\in L^\infty\Big(0,T; BV\big(\R^{n+1}, \CurBrac{0,1}\big)\Big)
\Eeqn
such that as $\dt\Converge 0$, $\chi_{\dt}$ converges to $\chi^*$
in $L^1(\mathbb{R}^{n+1}\times[0,T])$.
Furthermore there exists a function
$v \in L^1\big(0,T; L^1\big( | \nabla \chi^* | \big) \big)$
such that for all 
$\zeta \in C^{\infty}( \bar{ \Lambda \times [0,T] }, \mathbb{R}^n )$ 
with $\zeta = 0$ on $\partial \Lambda \times [0,T]$ and
$\xi \in C^{\infty}( \bar{ \Lambda \times [0,T] }, \mathbb{R} )$ 
with $\xi = 0$ on $\partial \Lambda \times [0,T] \cup \Lambda \times \{0\}$,
the following two properties hold,
\begin{eqnarray}
\int_0^T \int_{\Lambda} \big( \Text{div} \zeta - \frac{\nabla \chi^*}{| \nabla \chi^* |} \nabla \zeta \frac{\nabla \chi^*}{| \nabla \chi^* |} \big) |\nabla \chi^* | 
& = & - \int_0^T \int_{\Lambda} v \zeta \cdot \nabla \chi^*,
\label{WeakFormulation-a}\\
\int_0^T \int_{\Lambda} \chi^* \partial_t \xi + \int_{\Lambda} \chi_0 \xi(0) 
& = & - \int_0^T \int_{\Lambda} v \xi |\nabla \chi^* |.
\label{WeakFormulation-b}
\end{eqnarray}
In the above, $\Lambda$ is a fixed and sufficiently large ball in $\R^{n+1}$.
\end{thm7}

Some remarks are in place.

(i) The property of $\chi^*$ above implies that for a.e. $t$,
$\chi^*(t) = \mathbbm{1}_{\Omega^*(t)}$ for some $\Omega^*(t)\subseteq\R^{n+1}$
which is a \emph{set of finite perimeter}. In fact, since the results are proved
within the realm of classical solution, we actually have
$\chi^*\in C\Big(0,T; BV\big(\R^{n+1}, \CurBrac{0,1}\big)\Big)$.
Furthermore, the set $\Omega^*(t)$, or more exactly its boundary 
$\Gamma^*(t)=\Bdry\Omega^*(t)$, evolves smoothly in time.
Identity \eqref{WeakFormulation-b} shows that $v(t)$ is the
velocity function of $\Gamma^*(t)$ and
\eqref{WeakFormulation-a} shows that
$v(t)$ is given by the mean curvature of $\Gamma^*(t)$.
We refer to \cite{luckhaus} and also Section \ref{ProofConvSec}
for more detail explanation about these concepts and notations.

(ii) For some quantitative estimate and reasoning, 
more precise condition on the initial manifold $\Gamma_0$ 
will be given in Section \ref{InitialData}.

(iii) The time $T$ appearing in the Theorem depends only on geometric 
properties of the initial manifold $\Gamma_0$ but not on $\dt$. 
We have not yet shown that $T$ coincides with the 
maximal time $T^*$ for which the smooth MMC flow starting from 
$\Gamma_0$ exists.
This is because with the current estimate, we have only 
$C^{1,\alpha}$-convergence (in space) of the $\Gamma_{\dt}$. 
We expect that this can be improved to $C^{2,\alpha}$-convergence so that
the curvature of $\Gamma_\dt$ will also converge. 
Then $T$ would be the same as $T^*$.
See the remark (ii) after \Bf{Theorem \ref{Finite Time Stability}} 
for more detail discussion.

(iv) The convergence statement is that $\Gamma_\dt$ converges to a 
solution of MMC in the BV-sense. It will be ideal if we can show that this 
solution coincides with the classical (strong) solution. Again, this can be 
done if we can demonstrate $C^{2,\alpha}$-convergence.

In order to present the strategy of our proof, we introduce the 
following notations. 
Let $\Gamma$ be a smooth compact $n$-dimensional embedded manifold
in $\R^{n+1}$. We use $M(M_0, t)$, or simply $M_t$ if the choice of $M_0$ is clear, 
to denote the 
\Em{solution of MMC at 
time $t$ starting from $M_0$}. By regularity theory of MMC, 
$M_t$ is a smooth embedded manifold for small $t > 0$. 
Each manifold $\Gamma_{k\dt}$, obtained from the MBO scheme
will be called a \emph{numerical manifold}.

The proof of Theorem \ref{MainTheorem} relies on two main results: 
\emph{consistency} and \emph{stability} statements. 
The former states that for $\dt \ll 1$, 
the Hausdorff distance between
$M(\Gamma_{k\dt}, \dt)$ and 
$\Gamma_{(k+1)\dt}$ is of order $o(\dt)$. The exact order 
of $o(\dt)$ will be made precise during the proof.
The latter states that the curvatures of the numerical manifolds
$\Gamma_{k\dt}$ are uniformly bounded so that the error 
does not accumulate over repeated iterations. 
We next give an outline of the proof.

\subsection{Step I. Construction of ansatz.}
Recall that at each \Bf{Step 1} of the MBO scheme, 
we solve the linear heat equation
\begin{eqnarray}
U_t - \Delta U & = & 0 \space \text{,} \indent \text{for } (x, t) \in \mathbb{R}^{n+1} \times (0, \dt),
\label{Heat Eq 1}
\\
U( \cdot, 0 ) & = & \left\{ 
  \begin{array}{l l}
    1 & \quad \text{if $x \in \Omega_0$}\\
    -1 & \quad \text{if $x \notin \Omega_0$}
  \end{array} \right.,
\label{Heat Eq 2}
\end{eqnarray}
where $\Omega_0$ is an open set in $\mathbb{R}^{n+1}$ with smooth boundary
$\Gamma_0 = \Bdry\Omega_0$.  
The main idea is to formulate an appropriate ansatz for the
solutions to \eqref{Heat Eq 1}-\eqref{Heat Eq 2}
which can be easily compared to the solution of MMC.

For this, we will make use of $M(\Gamma_0,t)$, the solution of MMC
starting from $\Gamma_0$. By regularity of MMC, for a short time,
the solution $M(\Gamma_0,t)$ will remain smooth and continue to bound 
a set $\Omega(t)$.
Now let $r(\cdot, t)$ be the signed distance function to 
$M(\Gamma_0, t )$, with the convention that $r > 0$ in the interior 
$\Omega(t)$ of $M(\Gamma_0, t )$: 
\Beqn\label{sign.dist}
r(x,t) = \text{sdist}(x, M(\Gamma_0,t)) 
=\left\{
\begin{array}{ll}
\text{dist}(x,M(\Gamma_0,t)), & \text{for}\,\,\,x\in \Omega(t),\\
-\text{dist}(x,M(\Gamma_0,t)), & \text{for}\,\,\,x\notin \Omega(t).
\end{array}
\right.
\Eeqn
Then we decompose $U$ as
\Beqn\label{U.decompose}
U = U^0 + U^1
\Eeqn
where the leading order term $U^0$ is given by
\begin{align}\label{U^0}
U^0(x , t ) = \frac{2}{\sqrt{\pi}}  \int^{ \frac{r(x,t)}{2\sqrt{t}} }_0 
\exp(-y^2) dy.
\end{align}
Note that $U^0$ is the modified error function which
solves the following one dimensional linear heat equation 
in $(r, t)$:
\Beqn\label{U0HeatEqn}
U^0_t = \partial_{rr} U^0,\,\,\,\,\,\,
U^0(r, 0 ) = 2 \cdot \mathbbm{1}_{(0, \infty)}(r) - 1.
\Eeqn
For future reference, we write down the following formula:
\Beqn\label{U0_r}
\partial_r U^0 = \Lover{\sqrt{\pi t}}\exp\Brac{-\frac{r^2}{4t}}.
\Eeqn

The error term of the above ansatz is given by the function 
$U^1 = U-U^0$. Since $U$ solves the linear heat equation, we have that 
\begin{align*}
U^1_t - \Delta U^1 = - (U^0_t - \Delta U^0).
\end{align*}
Furthermore, as the initial data for $U$ is the same as that for $U^0$,
we have $U^1(\cdot, 0) \equiv 0$. 
Hence by means of variation of parameters, $U^1$ can be represented as
\begin{align}\label{U^1prelim}
& U^1 (x, t) = - \int_0^t \int_{\mathbb{R}^{n+1}} G(x -y , t - \tau) 
\big(\partial_\tau - \Delta\big)U^0(y, \tau)dy d\tau,
\end{align}
where \MathSty{G(x,t) = \frac{e^{-\frac{\Abs{x}^2}{4t}}}{\brac{4\pi t}^{\frac{n+1}{2}}}}
is the heat kernel defined on $\mathbb{R}^{n+1}\times\mathbb{R}_+$.
The explicit representations (\ref{U^0}) and (\ref{U^1prelim}) allow us 
to carry out a detailed analysis of the solution to 
(\ref{Heat Eq 1})-(\ref{Heat Eq 2}).

\subsection{Step II. Consistency estimate.}
The consistency result in this paper is phrased in terms of 
the Hausdorff distance $d_{\mathcal{H}}$ 
between $\Gamma_{\dt}$ and $M( \Gamma_0, \dt)$. 
This is stated by the following estimate 
(\Bf{Theorem \ref{consistency}}):
\Beqn\label{consistency0-d}
d_{\mathcal{H}} \big( \Gamma_{\dt} , M(\Gamma_0, \dt) \big) 
\leq \| \CB_0 \|^2 \CO\big( \dt ^{\frac{3}{2}}\big).
\Eeqn
Here $\| \CB_0 \|^2$ is a norm placed on the Weingarten map $\CB_0$ for 
$\Gamma_0$ which incorporates curvature information of $\Gamma_0$
(see \Bf{Section \ref{diff.geom}} for definition).
Heuristically this means that $\Gamma_{\dt}$ lies within a tubular 
neighborhood of $M(\Gamma_0, \dt)$ of radius 
$ \| \CB_0 \|^2 \CO\big( \dt^{\frac{3}{2}}\big)$.  

The proof of \eqref{consistency0-d} makes use of the following
formula (\Bf{Lemma \ref{HeatOperatorFormLem}}, \eqref{HeatOperatorForm}):
\[
\big(\partial_t - \Delta\big)U^0(x, t) 
= 
\frac{1}{\sqrt{\pi t}}\exp\left(-\frac{r(x,t)^2}{4t}\right) \sum_{i = 1}^n 
\frac{r(x,t)\kappa_i^2(\Pi(x),t)}{1 - r(x,t)\kappa_i(\Pi(x),t)}
\]
where the $\kappa_i(\Pi(x),t)$'s are the principal curvatures of the 
point $\Pi(x)$ on $M(\Gamma_0, t)$ which is closest to $x$.
By \eqref{U^1prelim}, the above
gives $|| U^1 ||_{L^\infty(\R^{n+1})} \leq  \| \CB_0\|^2 \CO(\dt)$. 

Next, near $M(\Gamma_0, \dt)$, 
$U^0 ( \cdot, \dt)$ roughly equals
$\dfrac{r(\cdot, \dt)}{\sqrt{\dt}}$.
Hence, $\Gamma_{\dt}$, which is given by the 
\Em{zero set of $U^0 (\cdot, \dt) + U^1 (\cdot, \dt)$}, 
corresponds to \Beqn\label{consistency0-r}
r(\cdot, \dt) 
= \| \CB_0 \|^2 \CO\big(\dt^{\frac{3}{2}}\big) 
\ll \CO(h)
\Eeqn
giving the consistency of the scheme (see Figure 1).

\begin{figure}[h]
\center
\includegraphics[scale=0.3]{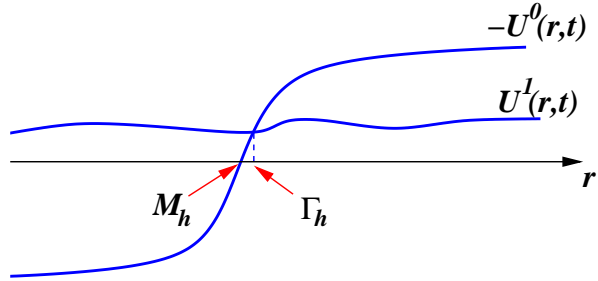}
\caption{Illustration of $U^0$, $U^1$. 
Note that $U^0$ vanishes at $M_h = M(\Gamma_0, h)$, the solution at time $h$
of MMC, starting from $M_0$.
The intersection between $-U^0$ and $U^1$ gives $\Gamma_h$. 
}
\end{figure}

The convergence rate in (\ref{consistency0-r}) can in fact be improved to 
$\CO(\dt^2)$.
Estimate \eqref{consistency0-r} simply makes use of the $L^{\infty}$-norm of 
$\big(\partial_t - \Delta\big)U^0(x, t)$. 
The improved rate comes by performing a more precise point-wise analysis. This will involve more analytical computation but
it is quite similar to what is done in the stability estimates 
-- 
see Remark \ref{improved_convergence_rate} in Section \ref{RefinedU12}.

\subsection{Step III. Stability Estimates}
Note that the consistency statement \eqref{consistency0-d} 
involves the curvature of
$\Gamma_{\dt}$. In order to ensure that such a statement can be extended
to multiple iterations, we need to estimate the curvature of 
$\Gamma_{\dt}$ in terms of the curvature of $\Gamma_0$.

The first step in doing this is to describe 
$\Gamma_{\dt}$ as a graph over 
$M_\dt = M(\Gamma_0, \dt)$. This is achievable due to the fact that
near $M_\dt$, we have
\[
\nabla U\approx \partial_r U^0 
= \Lover{\sqrt{\pi h}}\exp\Brac{-\frac{r^2}{4h}} \gg 1
\,\,\,\,\,\,\text{(since $r\sim\CO(h^\frac{3}{2})$).}
\]
Hence Implicit Function Theorem (IFT) gives us that
locally $\Gamma_\dt$ can be written as a graph of a function $g$ 
over the tangent plane $\mathbb{T}_{M_\dt}(p_0)$ at some point $p_0$ on 
$M_\dt$ (see Figure 2). The next step is to relate the curvature
of $\Gamma_\dt$ to the second derivatives of $g$. 

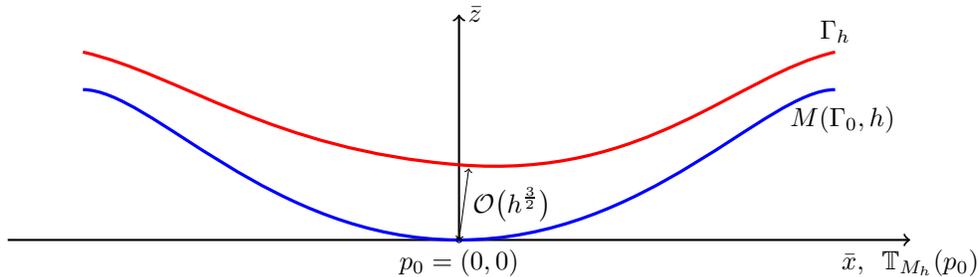
\begin{figure}[h]
\begin{center}
\begin{tikzpicture}
\draw[->, thick] (-6,0) -- (6,0);
\draw[->, thick] (0,0) -- (0,3);
\draw[blue, very thick] (-5, 2) .. controls (-4.1, 2) and (-2.5, 0) .. (0, 0);
\draw[blue, very thick] (0,0) .. controls (2.5,0) and (4.1, 2) .. (5, 2);
\draw[red, very thick] ( -5, 2.5) .. controls (-3.7, 2.2) and (-2.5, 1.2) .. (0, 1); 
\draw[red, very thick] (0, 1) .. controls (2.5, 0.8) and (3.7, 2.2) .. (5, 2.5);
\draw (6,0) node[below]{$\bar{x},\,\,\,\mathbb{T}_{M_\dt}(p_0)$};
\draw (0,3) node[right]{$\bar{z}$};
\draw (0,0) node[below]{$p_0=(0,0)$};
\draw (0,0) circle (1pt);
\draw (5.1, 1.9) node[below]{$M(\Gamma_0, \dt)$};
\draw (5, 2.5) node[above]{$\Gamma_{\dt}$};
\draw[very thin, <->] (0,0) -- node[midway, right]{$\mathcal{O}\big(\dt^{\frac{3}{2}}\big)$} (.13, .95);
\end{tikzpicture}
\caption{The manifolds $M_h = M(\Gamma_0,h)$ and $\Gamma_h$. Note that
$d_{\Cal{H}}(M_h, \Gamma_h) = \CO(h^\frac{3}{2})$ and $\Gamma_h$ locally
can be written as a graph over $\mathbb{T}_{M_\dt}$: 
$\Gamma_h = \CurBrac{(\bx, g(\bx)): \bx\in\mathbb{T}_{M_\dt}}$.}
\label{fig:Graph}
\end{center}
\end{figure}

The necessary computations are also facilitated by IFT.
Letting $\CB_\dt$ be the Weingarten map
of $\Gamma_\dt$, and $\bar{x}$ be the coordinates of the tangent plane
$\mathbb{T}_{M_\dt}(p_0)$, we have 
(\ref{A}, \ref{1st deriv g}, \ref{2nd deriv g}) 
\Beqn
\CB_\dt
=
\nabla^2_{\bar{x}} g 
+ \nabla_{\bar{x}}^2 U^1 \CO(\sqrt{\dt})
+ \| \CB_0 \|^3 \CO(\dt).
\Eeqn
Careful analysis of $\nabla_{\bar{x}}^2 U^1$ via 
(\ref{U^1prelim}) and \eqref{HeatOperatorForm} gives
(\Bf{Lemma \ref{U1-2ndDerEst}}):
\Beqn\label{ddU1.est}
\| \nabla_{\bar{x}}^2 U^1 \|_{L^\infty}
\leq \| \CB_0 \|^3 \CO(\sqrt{ \dt })
\Eeqn
leading to the following \emph{one-step stability estimate}
(\Bf{Theorem \ref{Stability Estimate}}):
\Beqn\label{OneStepStabEst}
\| \CB_{\dt} \| \leq \| \CB_0 \| + \|\CB_0\|^3 \CO(\dt).
\Eeqn

The most daunting computation is the estimate \eqref{ddU1.est}. 
This is because \eqref{U^1prelim} expresses $U^1$ in terms of 
higher order derivatives of $U^0$. The analysis thus needs to 
make crucial use of some regularity theory of MMC, in particular
the decay estimates for the Weingarten map of
$M_t$ (\Bf{Lemmas \ref{curvature bound}} and 
\Bf{\ref{regularity of curvature}}). 

Once we have \eqref{ddU1.est}, it can be iterated by means of a discrete
Gronwall inequality. This leads to that the curvatures
of our numerical manifolds $\Gamma_{k\dt}$ are bounded uniformly
over multiple iterations of the scheme (\Bf{Theorem \ref{Finite Time Stability}}).

\subsection{Step IV. Convergence.}
By the consistency statement and curvature bound from the previous steps, 
together with some geometric argument to exclude the occurrence of 
self-intersection, the sequence of $\Gamma_{n\dt}$ can be shown to 
converge to a limit manifold in the Hausdorff distance and also 
the $C^{1,\alpha}$-norm.
The final step is to identify the equation satisfied by the limit.

We find the weak formulation of MMC using BV-functions 
or sets of finite perimeter as used in \cite{luckhaus} 
the most convenient for our purpose. 
Using the notation of Theorem \ref{MainTheorem}, we show that 
$\chi_{k \dt}$ converges to a limiting function 
$\chi^*_t \in BV( \mathbb{R}^{n+1}; \{0,1\})$ as $k \dt \rightarrow t$.  
Furthermore, $\chi^*$ solves MMC in the sense of 
\eqref{WeakFormulation-a}-\eqref{WeakFormulation-b}.  
The key step in establishing this is to prove that the area converges
in the sense that
\Beqn
\int | \nabla \chi_{k \dt} | \rightarrow \int | \nabla \chi^*_t |.
\Eeqn
(The above is assumed in \cite{luckhaus}.)
The main ingredient in doing this is the \emph{Ball Lemma}
(\Bf{Lemma \ref{BallLemma}}) by which we may place a tubular
neighborhood with uniform thickness such that 
$\Gamma_{k\dt}$ remains embedded.

\subsection{Some notations from geometry of surfaces}\label{diff.geom}
The reference for this section is \cite{ecker2004regularity}, in particular
Appendix A.  
From our perspective and application point of view, we take the definition
of the mean curvature $H$ as the negative first variation of the area 
functional.
But for the sake of performing analytical computation, we will relate $H$
to the Weingarten map of a manifold. Essentially $H$ is defined to be
the trace of the Weingarten map.

Recall that for an $n$-dimensional manifold $M=\Bdry\Omega$ embedded in 
$\mathbb{R}^{n+1}$ given by an embedding map: $M=F(D)$ where
$F: D\subseteq \mathbb{R}^n \rightarrow \mathbb{R}^{n+1}$, the second fundamental form is the symmetric bilinear form on the tangent 
bundle $TM$ of $M$ given by
\Beqn
A_{ij} = \langle \partial_i \N , \partial_j F \rangle 
= - \langle \N, \partial^2_{ij} F \rangle
\Eeqn
where $\N$ is the unit outward normal for $M$.
Inherently related is the \emph{Weingarten map}, which is the mapping $L$ 
from $TM$ to $TM$ determined by,
\Beqn
L(\mathbf{u}) = -\nabla_{\mathbf{u}} \N.
\Eeqn
In the coordinate system determined by $F$, the matrix corresponding to the Weingarten map is given by,
\Beqn
A^i_j = g^{ik} A_{kj}, 
\Eeqn
where $ ( g^{ij} ) = ( g_{ij} )^{-1} = 
\big( \langle \partial_i F , \partial_j F \rangle \big)^{-1}$.
The eigenvalues $\kappa_1, \ldots, \kappa_n$ of $L$ are called the
principal curvatures of $M$.

With the above notations, the mean curvature $H$ of $M$ is given by
\Beqn
H = \Lover{n}\Text{div}_M \N := \Lover{n}\sum_{i,j} g^{ij} \langle \partial_i \N, \partial_j F \rangle
\Eeqn
which can be related to the trace of the Weingarten map of $M$ as follows,
\Beqn
H = \frac{1}{n} \sum_i A^i_i =\Lover{n}\Brac{\kappa_1 + \cdots + \kappa_n}.
\Eeqn

Of particular relevance is the case when $F$ is the graph of a function 
$f$ over $\mathbb{R}^n$, i.e. $F(x) = ( x, f(x) )$ for $x\in\R^n$.  
In this case, we have
\[
A_{ij}  =   \dfrac{ \nabla ^2 f }{\sqrt{ 1 + | \nabla f |^2 }},
\] 
\[
( g_{ij} ) = I + \nabla f \otimes \nabla f,
\]
and
\Beqn\label{weingarten.graph}
( A^i_j ) =  
\left ( I - \dfrac{ \nabla f \otimes \nabla f}{1 + | \nabla f |^2 }\right ) 
\dfrac{ \nabla ^2 f }{\sqrt{ 1 + | \nabla f|^2 } }.
\Eeqn
The mean curvature $H$ is then given by,
\Beqn\label{mean.graph}
H 
= \frac{1}{n} \sum_{i=1}^n A_i^i 
= \text{div}\frac{\nabla f}{\sqrt{1+\Abs{\nabla f}^2}}.
\Eeqn

For the rest of this paper, 
we will use $A_M$ (or simply $A$) to denote the \emph{Weingarten map}
of some general manifold $M$.
But to emphasize the importance of the numerical manifolds $\Gamma_{k\dt}$'s, 
we will use $\CB_{k\dt}$ to denote their Weingarten maps.
For the analysis in the rest of the paper, we will use the following norm 
for the Weingarten map of a manifold $M$,
\Beqn\label{A2}
|A(p)| 
= \sqrt{\text{tr}(AA^T)}
= \sqrt{\sum_{i,j} A^j_i(p)A_j^i(p)},
\,\,\,\,\,\,\text{and}\,\,\,\,\,\,
\| A\| = \sup_{ p \in M} |A(p)|.
\Eeqn

A useful observation is that if the Weingarten map is bounded, then we can 
have a quantitative estimate about the size of the region over which the 
manifold can be written as a graph. To be specific, let $p\in M$, $\T_p(M)$ be 
its tangent plane, and locally near $p$, $f$ be the graph of $M$ over 
$\T_p(M)$. From 
\eqref{weingarten.graph}, we have
\[
\nabla ^2 f
=
\sqrt{ 1 + | \nabla f|^2 }
\left ( I + \nabla f \otimes \nabla f\right)
( A^i_j ).
\]
Hence
\MathSty{\Norm{\nabla ^2 f}
\leq
C\Brac{1 + | \nabla f|^2}^{\frac{3}{2}}
\Norm{A}
}.
Upon integrating in space, there is a universal constant $C_*=C_*(n)$ 
such that 
\Beqn\label{RadiusGraph}
\CurBrac{
x\in\T_p(M),\,\,\abs{x} \leq \frac{C_*}{\Norm{A}}}
\subseteq
\big\{x: \Abs{\nabla f(x)} < \infty\big\}.
\Eeqn
The above implies that the connected component of 
$M\Intersect B_{\frac{C_*}{\Norm{A}}}(p)$
containing $p$ is completely given by the graph of $f$.

\subsection{Time Step Size $\dt$ in Relation to the Geometry of Initial 
Manifold $\Gamma_0$}
\label{InitialData}
Here we discuss the requirement for the initial manifold $\Gamma_0$
which is a compact smooth embedded $n$-dimensional 
manifold in $\mathbb{R}^{n+1}$.  The time step $\dt$ will be sent to 
zero in the convergence statement. But to be quantitative, we will 
specify its smallness with respect to two geometric quantities pertaining 
to $\Gamma_0$.
In the following, we introduce a small constant $\delta$ such that 
$\dt \leq \delta$.

The first requirement is a local property. The value of 
$\delta$ satisfies
\Beqn\label{InitialDataAssump}
 \| \CB_0 \| ( \delta  | \log \delta  | )^{\frac{1}{4}} \leq 1.
\Eeqn
The second is a more global condition. To describe this precisely, 
we first define the following \emph{Ball Property}.

\begin{thm7}{\bf Ball Property.}
\label{BallProperty}
Given an embedded n-dimensional manifold $M$ which is the boundary of a
subset $\Omega\subseteq \R^{n+1}$, 
i.e., $M = \Bdry\Omega$, we say that $M$ satisfies the 
\Em{ball property with radius $r$} if for every $p \in M$, 
there are two balls $B_{r,p}^\Text{int}$ and $B_{r,p}^\Text{ext}$ 
(interior and exterior) 
with radius $r$ such that $B_{r,p}^\Text{int} \subseteq \Omega$ and
$B_{r,p}^\Text{ext}\subseteq\R^{n+1}\backslash\Omega$ and
\Beqn
B_{r,p}^\Text{int}\intersect M = B_{r,p}^\Text{ext}\intersect M = \curBrac{p}.
\Eeqn
\end{thm7}
\noindent
Note that the above condition is stronger than simply requiring that the curvature of $M$ is bounded from above by $\Lover{r}$. It is some kind of ``uniform embeddedness'' condition.

Now let $R(\Gamma_0)$ be the maximal radius for which $\Gamma_0$ satisfies the ball property. Then we require that there is a small constant $\rho$ such that,
\Beqn\label{InitialDataAssump2}
( \delta | \log \delta | )^{\frac{1}{4}} \leq \rho \leq \frac{R(\Gamma_0)}{2}.
\Eeqn

The stability analysis to be carried out in
Section \ref{Stability.Sec} which includes the regularity statement 
and the Ball Lemma will imply that for $\dt \leq \delta$
and $k = \CO(\Lover{\dt})$, we have that
$$ 
\|\CB_{k \dt} \| ( \dt | \log \dt | )^{\frac{1}{4}} \leq 1
$$
and the $\Gamma_{k\dt}$'s will all be embedded manifolds.

The two conditions \eqref{InitialDataAssump} and 
\eqref{InitialDataAssump2} will be assumed for the rest of the paper.

\subsection{Notations and Conventions}\label{Notation}
Throughout the estimates in the paper, constants and bounded functions 
will be grouped together as $C$ and $\mathcal{O}(1)$ respectively.  
These are terms bounded by constants that do not depend on
$\|\CB_0\|$ or $\dt$.
From one line to the next, the $C$ and $\mathcal{O}(1)$ terms may change, 
but we will still use $C$ and $\mathcal{O}(1)$. We will also use the following conventions:

\begin{enumerate}
\item[(i)]
$A=\Cal{O}(1)$ or $A(\cdot)=\Cal{O}(1)f(\cdot)$: there is some constant $C$ such that
\[
\Abs{A} \leq C
\,\,\,\,\,\,\text{or}\,\,\,\,\,\,
\Abs{A(\cdot)} \leq C\Abs{f(\cdot)}.
\]

\item[(ii)]
$A \lesssim B$ or $A\gtrsim B$: there is some constant $C$ such that
\[
A\leq C B
\,\,\,\,\,\,\text{or}\,\,\,\,\,\,
A \geq CB
\]

\item[(iii)]
$A \ll 1$ or $A \gg 1$: $A$ is a sufficiently small or large constant.

\item[(iv)]
$A\approx B$: the following limiting behavior holds:
\[
\lim \frac{A}{B} = 1.
\]
The meaning of the limit will be specified or clear from the context.
\end{enumerate}

\section{Ansatz and Its Consistency}

In this section, we prove the consistency of the scheme by analyzing
the ansatz $U^0$ and the error term $U^1 = U-U^0$ 
defined in \eqref{U^0} and \eqref{U.decompose}.

Recall that the initial manifold $\Gamma_0$ is a smooth, embedded
$n$-dimensional manifold in $\R^{n+1}$. By the regularity of MMC,
there exists a positive number $\rho = \rho(\Gamma_0) > 0$ such that 
for any $0 < t < \dt \ll 1$, the signed distance function $r$ to 
$M(\Gamma_0, t)$ is a smooth function in the $\rho$-tubular 
neighborhood $T_\rho(\Gamma_0)$ of $\Gamma_0$,
\Beqn
T_\rho(\Gamma_0) = \CurBrac{
q\in\R^{n+1}: \text{dist}(q, \Gamma_0) \leq \rho
}.
\Eeqn
Hence inside $T_\rho(\Gamma_0)$,
the $U^0$ and $U^1$ are smooth functions of $(x,t)$.

Inside $T_\rho(\Gamma_0)$, we have the following representation formula 
for $U^1$ in terms of the heat kernel $G$ on $\mathbb{R}^{n+1}\times\R_+$ 
(c.f. \cite{jost2002partial}):
\begin{align}
U^1(x,t) = &\int_0^t \int_{T_\rho(\Gamma_0)} G(x -y , t - \tau) 
\big( \partial_t  - \Delta\big) U^0 (y, \tau) dy d\tau \nonumber\\
& + \int_0^t \int_{\partial T_\rho(\Gamma_0)} G( x - y, t - \tau) \frac{\partial U^1}{\partial \nu}(y , \tau) - U^1(y, \tau) \frac{\partial G}{\partial \nu} (x - y, t - \tau) dS d\tau.
\label{U^1}
\end{align}
where $\nu$ is the unit outward normal to $\Bdry T_\rho(\Gamma_0)$.

We give some remarks about the above ansatz.

(i) Recall that $\rho$ is some small constant satisfying the condition 
\eqref{InitialDataAssump2}. Hence the signed distance function to 
$M(\Gamma_0,t)$, $r(\cdot, t)$, is smooth in the set 
$\{ r(\cdot, t) < 2 \rho\}$ (for $0 \leq t \leq \dt \leq \delta$).  
It will be a consequence of Lemma \ref{BallLemma} that 
the same $\rho$ will work for multiple iterations, 
i.e. a tubular neighborhood of radius $\rho$ may be placed about 
$\Gamma_{k \dt}$ with the ansatz constructed in the same manner 
(see Theorem \ref{Finite Time Stability}).

(ii)
The second term in (\ref{U^1}), integration on the boundary 
$\Bdry\Gamma_\rho(\Gamma_0)$, produces exponentially small terms. 
Specifically, the following estimates hold for $U$ and $U^0$
on $\partial T_{\rho}(\Gamma_0)$,
\begin{center}
$\left | U - U^0 \right | \lesssim e^{\frac{-\rho^2}{16 \dt}}$,\,\,\,\,\,\,
$\left| \partial_rU^0 \right | \lesssim \frac{1}{\sqrt{\dt}} e^{\frac{-\rho^2}{16 \dt}}$,
\,\,\,and\,\,\,
$\left| \nabla U \right | \lesssim \frac{1}{\sqrt{\dt}} e^{\frac{-\rho^2}{8 \dt}}$.
\end{center}
Furthermore, the area of $\partial T_{\rho} (\Gamma_0)$ is roughly
two of that of $\Gamma_0$.
Hence we will simply omit it in our analysis.  

For the following, given any constant $c$, we define 
$\CurBrac{r=c}$ to the manifold which is at a distant $c$ to $M_t$:
\Beqn\label{manifold.at.distance}
\CurBrac{r=c} = \CurBrac{y\in\R^{n+1}: \Text{sdist}(y, M(\Gamma_0, t))=c}.
\Eeqn
We note that for $c\ll 1$, $\CurBrac{r=c}$ is a smooth manifold.
Let also $H_{\{r=c\}}: \CurBrac{r=c}\MapTo\R$ be the mean curvature
function of $\CurBrac{r=c}$. 

\begin{thm4}\label{HeatOperatorFormLem}
For any $(x,t) \in T_\rho(\Gamma_0) \times [0, \dt]$,
the following formula holds,
\Beqn\label{HeatOperatorForm}
(\partial_t - \Delta)U^0(x, t) 
= \Big( 
H_{\{ r = r( x, t) \}}(x) - H_{\{ r = 0 \}}\big(\Pi(x)\big)\Big)\partial_rU^0(x,t),
\Eeqn
where $\Pi(x)$ is the point in $M(\Gamma_0,t)$ closest to $x$.
In addition, we have
\Beqn \label{H-Diff-r}
H_{\{ r = r(x, t) \}}(x) -
H_{\{ r = 0 \}}(\Pi(x))
= r(x,t) \psi \Big(r(x,t),\kappa_1(\Pi(x),t), \kappa_2(\Pi(x),t), \ldots, 
\kappa_n(\Pi(x),t) \Big)
\Eeqn
where the $\kappa_i$ are the principal curvatures of $M( \Gamma_0, t)$ 
evaluated at $\Pi(x)$, and 
\Beqn\label{psi-formula}
\psi(r,\kappa_1, \kappa_2, \ldots, \kappa_n)
= \sum_{i=1}^n\frac{\kappa_i^2}{1-r\kappa_i}.
\Eeqn 
\end{thm4}
See Figure 3 for an illustration of some of the notations appearing above.

\begin{figure}[h]
\center
\includegraphics[scale=0.3]{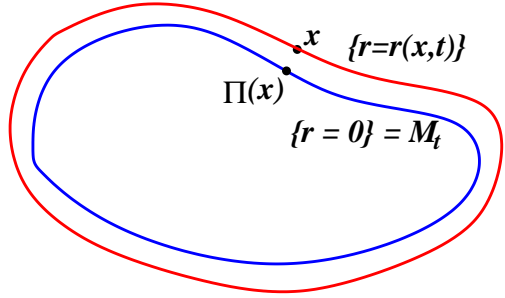}
\caption{$\CurBrac{r=r(x,t)} =
\CurBrac{y\in\R^{n+1}: \Text{sdist}(y, M(\Gamma_0, t))=r(x,t)}$.
$\Pi(x)$ is the projection of $x$ onto $M_t$}
\end{figure}
The above will be proved in \Bf{Section \ref{sdist-sec}}.
We will frequently abuse the notation by simply writing $\psi(x,t)$ 
in place of $\psi \big( r(x,t),\kappa_1(\Pi(x),t), \kappa_2(\Pi(x),t), \ldots, 
\kappa_n(\Pi(x),t) \big)$. By the formula \eqref{U0_r} for 
$\partial_r U^0$, the expression $(\partial_t - \Delta)U^0(x, t)$
takes the following form
\Beqn\label{HeatOperatorFormRHS}
\psi(r, \kappa) \dfrac{r( x, t)}{\sqrt{\pi t}} 
\exp\Brac{ -\frac{r^2(x, t)}{4t}}.
\Eeqn
We also note the following estimates:
\Beqn\label{psi-exp-r-est}
\Norm{\psi} \lesssim \Norm{A_{M(\Gamma_0,t)}}^2
\,\,\,\Text{for}\,\,\,\Abs{r}\ll \Lover{\Norm{A_{M(\Gamma_0,t)}}},
\,\,\,\,\,\,\text{and}\,\,\,\,\,\,
\Abs{\frac{r}{\sqrt{t}}\exp\Brac{-\frac{r^2}{4t}}} = \CO(1)
\,\,\,\Text{for all $r,\,t$}.
\Eeqn

With the above, we now proceed to prove the following consistency 
statement for the MBO scheme.
\begin{thm1}\label{consistency}
For any $\dt\leq \delta$, there is a constant $C$ depending only on the
spatial dimension such that the following estimate holds,
\Beqn
d_{\mathcal{H}}( \Gamma_{\dt} , M(\Gamma_0, \dt) ) 
\lesssim  
\| \CB_0 \|^2 \dt ^{\frac{3}{2}}
\Eeqn
where we have used $\CB_0$ to denote $A_{\Gamma_0}$.
\end{thm1} 
\begin{proof}
By \eqref{HeatOperatorFormRHS} and \eqref{psi-exp-r-est},
we have for $0 < t < \dt$ that
\[
\Abs{(\partial_t-\Delta)U^0 (x, t)}
= 
\Abs{\psi(r, \kappa) \dfrac{r( x, t)}{\sqrt{\pi t}} 
\exp\Brac{ -\frac{r^2(x, t)}{4t}}}
\lesssim  \| A_{M(\Gamma_0, t)}\|^2.
\]
Substituting this estimate into the integral in (\ref{U^1}) and by
the $L^\infty$-$L^\infty$ estimate of the heat operator, we obtain that,
\begin{eqnarray*}
\Abs{U^1(x, \dt)}
& \leq & 
\Abs{\int_0^\dt \int_{T_\rho(\Gamma_0)} G(x -y , \dt - \tau)
\big( \partial_t  - \Delta\big) U^0 (y, t) dy d\tau}\\
& \leq & 
\int_0^\dt
\Norm{G(\cdot, \dt-\tau)}_{L^1(\R^{n+1})}
\Norm{(\partial_t-\Delta)U^0 (\cdot, \tau)}_{L^\infty(T_\rho(\Gamma_0))}
\,d\tau\\
& \lesssim &
\sup_{\tau \in [0, \dt]} \| A_{M( \Gamma_0, \tau)} \|^2 \dt.
\end{eqnarray*}
By Lemma \ref{curvature bound} (which appears in the next section), we have
\begin{align*}
\sup_{t \in [0, \dt]} \| A_{M( \Gamma_0, t)} \big 
\|& \leq \big( \| \CB_0 \| + C\| \CB_0 \|^3 \dt \big)
=  \| \CB_0 \| \Big (1 + C\| \CB_0 \|^2 \dt \Big )
\lesssim \| \CB_0 \|
\end{align*}
leading to
\Beqn \label{est1}
\big \| U^1(x, \dt) \big \|_{L^\infty(T_\rho(\Gamma_0))} \lesssim \| \CB_0 \|^2 \dt.
\Eeqn

Finally, by the representation formula \eqref{U^0} of $U^0$, 
for $r\ll 1$, we have
\Beqn
\label{est2}
U^0(x,t) 
= \frac{2}{\sqrt{\pi}} \int_0^{\frac{r}{2\sqrt{t}}} \exp(-y^2)dy
\approx \frac{r}{\sqrt{\pi t}}. 
\Eeqn Hence $U^0 + U^1$ can be zero only if $r$ satisfies
$$
\left | r(\cdot,\dt) \right | \lesssim \| \CB_0 \|^2 \dt^\frac{3}{2}.
$$
\end{proof}

\section{Stability}\label{Stability.Sec}
The stability estimates of this section will allow us to extend the 
consistency estimate from the previous section to 
multiple iterations in the MBO algorithm.
The first, key step is show that the curvature of 
$\Gamma_{\dt}$ can be controlled by the curvature of $\Gamma_0$.  
Then we prove a geometric result preventing $\Gamma_{\dt}$ to 
have self-intersection and hence $\Gamma_\dt$ remains embedded.
Lastly, tying this together with a discrete non-linear Gronwall inequality,
we are able to show that the curvatures of the $\Gamma_{k\dt}$'s
are uniformly bounded over multiple iterations. The crucial computation
and analysis rely on the regularity property of MMC.

We now state the two main results in this section. (We recall the notations about
differential geometry from \Bf{Section \ref{diff.geom}}.)

\begin{thm1}{\bf Stability over one time step.}\label{Stability Estimate}
There is a constant $C$ depending on the spatial dimension such that
for $\dt \leq \delta$, we have
\Beqn\label{onetimestability}
\| \CB_{\dt} \| \leq \| \CB_0 \| ( 1 + C \| \CB_0 \|^2 \dt ).
\Eeqn
\end{thm1}
\begin{thm1}{\bf Stability over multiple time steps.}\label{Finite Time Stability}
Let $\rho$ be some small number. For any constant
\Beqn
C_0 \in \left (\| \CB_0 \|, \min \left \{\frac{1}{( \delta | \log \delta | )^{\frac{1}{4}}}, \frac{1}{2 \rho} -1 \right\} \right],
\Eeqn
there exists a time $T = T(\Gamma_0, C_0)$ independent of $\dt$,
such that for $0 \leq k \leq \lfloor \frac{T}{\dt} \rfloor$, 
the following two statements hold.
\begin{itemize}
\item $\Gamma_{k \dt}$ is an embedded manifold.  More specifically, there is a uniform radius
$r_0 = r_0(\Gamma_0, C_0)$ so that $\Gamma_{k \dt}$ satisfies the Ball Property
(see Definition \ref{BallProperty}) with radius $r_0$.
\item $\| \CB_{k\dt} \| \leq C_0$.
\end{itemize}
(Note that for $\rho$ small enough, the interval for the choice of
$C_0$ is non-empty which by \eqref{InitialDataAssump2} can in fact 
be further simplified to 
\MathSty{C_0 \in \Big(\Norm{\CB_0}, \frac{1}{2\rho}-1\Big]}.)
\end{thm1}

Some remarks are in place.

(i) The specific manner in which $r_0$ and $T$ depend on $\Gamma_0$, and $C_0$ 
will be apparent in the proof of Theorem~\ref{Finite Time Stability}.
We note also that the established curvature bound can increase with each 
iteration, so that the larger the choice of $C_0$, the larger we may choose 
the convergence time $T$.

(ii) Theorem~\ref{Finite Time Stability} essentially proves that
$\sup_{k}\Norm{\Gamma_{k\dt}}_{C^2} < \infty$ for appropriate range of $k$. 
This only implies that $\Gamma_{k\dt}$ converges in $C^{1,\alpha}$ in space.
This is not sufficient
to show that the convergence time $T$ for our numerical scheme coincides 
with the maximum time $T^*$ for the existence of classical solution. 
This is because the constant $C$ in the
discrete one-time-step-stability estimate \eqref{onetimestability} 
might not be optimal and can be different from the continuous time case. 
Note however 
that the scaling in the growth rate of the curvature is the same in both
the discrete and continuous cases as illustrated by the MMC of a circle. 
Finite time blow-up
in the estimate over multiple time steps will definitely occur, just as in
the continuous case. But the two blow-up times might not be the same. 
We believe that with more refined analysis, we can in fact have 
$\sup_{k}\Norm{\Gamma_{k\dt}}_{C^{2,\alpha}} < \infty$
so that the $\Gamma_{k\dt}$ converges in $C^{2,\alpha'}$ for 
$\alpha' < \alpha$, i.e. the curvature $\CB_{k\dt}$ also converges 
(in $C^{\alpha'}$).
Then such a discrepancy between $T$ and $T^*$ can be removed.

As mentioned before, the proof of Theorem \ref{Stability Estimate} requires 
some regularity results of MMC. Specifically, we need the following 
two lemmas for surfaces following MMC.
For the next two results, we use $M_t$ to denote
$M(M_0,t)$, the solution of MMC starting from $M_0$. 

\begin{thm4}{\bf Bound on Curvature Growth of MMC.}
\label{curvature bound}
There is a constant $C$ which depends only on the spatial dimension such 
that for $0 \leq t \leq \dt$, we have
\begin{align*}\label{MMC assumption}
\|A_{M_t} \|  \leq \| A_{M_0} \| ( 1 + C \| A_{M_0} \|^2 t).
\end{align*}
\end{thm4}

\begin{thm4}{\bf Regularity of Higher Derivatives of MMC.}
\label{regularity of curvature}
Let $\N$ be the normal vector of $M_t$.
Suppose
\Beqn\label{ABdOneStep}
\|A_{M_t}\| \leq c_0 \text{ for $t \in [0, \dt]$.}
\Eeqn
Then for $t\leq \dt$, there is a constant $C$ depending on the spatial 
dimension such that,
\begin{eqnarray}
\| \nabla_{M_t} A_{M_t} \| &\leq& \frac{C \cdot c_0}{\sqrt{t}},\\
\text{and}\,\,\,\,\,\,
\| \partial_t \N \|_{L^{\infty}(M_t)} &\leq& \frac{C \cdot c_0}{\sqrt{t}}.
\end{eqnarray}
In the above,
$\|\nabla_{M_t} A_{M_t}\| = \sup_{M_t} \sqrt{ |  \nabla_{M_t} A_{M_t} |^2 }$
where $ |  \nabla_{M_t} A_{M_t} |^2 $ is the squared norm of the tensor 
$ ( \nabla_k A_{M_t})$. (See \cite[Appendix A]{ecker2004regularity} for 
detail explanation of the notations.)
\end{thm4}
\noindent
Note that assumption \eqref{ABdOneStep} holds by 
Lemma \ref{curvature bound} together with assumption 
\eqref{InitialDataAssump}.
The above Lemmas will be proved in \Bf{Sections \ref{ProofCurvatureBound}}
and \Bf{\ref{ProofCurvatureRegularity}}.

The stability analysis makes use of the Implicit Function Theorem (IFT).
From the formula \eqref{U^0} and \eqref{U^1} for $U^0$ and $U^1$, 
it can be seen that,
\[
\left | \partial_r U^0(x,t) \right | = \left |\frac{1}{\sqrt{\pi t}}\exp \left( -\frac{r^2(x,t)}{4t} \right) \right |
\gtrsim \Lover{\sqrt{t}}\,\,\,\,\,\,
\text{  for $|r(x,t)| \lesssim \sqrt{t}$,}
\]
and
\[
\| \nabla U^1(\cdot,t) \|_{L^{\infty}(T_\rho(\Gamma_0))} 
\lesssim \| \CB_0 \|^2 \sqrt{t}.
\]
We can then conclude that,
\[
\left | \partial_r U(x,t) \right | > 0 \text{  when  }\left | r(x,t) \right | 
\lesssim \sqrt{t}.
\]
Hence, we can apply the IFT to locally write 
$\Gamma_{\dt}$ as the graph of a function $g$ over the tangent plane to 
$M(\Gamma_0, \dt)$ at some reference point (see Fig. \ref{fig:Graph}).

Next we give some preparations for the proof of 
Theorem \ref{Stability Estimate}.

\subsection{Further Notations and Preparations}
\label{further.notations}

For simplicity, we use $M_t$ to denote $M( \Gamma_0, t )$ for $0 < t \leq \dt$. 
The Weingarten maps of $M_t$ and $\Gamma_\dt$ are denoted by 
$A_t$ and $\CB_\dt$, respectively. We further set $\CB_0 = A_0$.

The main goal of the next few sections is to estimate the derivatives of $U$ at
any point $p\in\Gamma_\dt$. By the consistency statement, we can assume that 
$p\in T_\rho(\Gamma_0)$.
To simplify the computation, we introduce the following coordinate system.
Let the origin be the point $p_0=\Pi_{M_\dt}(p)\in M_\dt$ 
which is closest to $p$. Then we write
\Beqn\label{coord_s0}
\R^{n+1} \cong \R^n\times\R \cong\T_{M_\dt}(p_0)\times \bbL(\N(p_0))
\Eeqn
where $\T_{M_\dt}(p_0)$ and $\bbL(\N(p_0))$ are the tangent plane 
and the line spanned by the outward normal vector 
$\N(p_0)$, both to $M_\dt$ at $p_0$.
We use $\bx\in\R^n$ and $\bz\in\R$ to denote the
coordinate variables of $\T_{M_\dt}(p_0)$ and $\bbL(\N(p_0))$. 
Then we have the following representations,
\Beqn\label{coord_s0_xz}
x=(\bx, \bz)\in\R^{n+1},\,\,\,
p_0 = (\bx=0, \bz=0),
\Eeqn
and near $p_0$,
\Beqn
\Gamma_\dt = \CurBrac{(\bx, g(\bx)): \bx\in\T_{M_\dt}(p_0)},
\,\,\,\text{in particular,}\,\,\,
p=(0,g(0)) \in \bbL(\N(p_0)).
\Eeqn
(See Figure 2.) 
With the above notations, we write
\[
U = U(\bx, \bz, t),\,\,\,
U^0 = U^0(\bx, \bz, t),\,\,\,
U^1 = U^1(\bx, \bz, t),\,\,\,\text{and}\,\,\,
r = r(\bx, \bz, t).
\]
We further note the following statements,
\Beqn
\nabla_\bx\Big|_{p_0} = \nabla_{M_\dt}\Big|_{p_0},
\Eeqn
and on $\bbL(\N(p_0))$ and in particular at $p$,
\Beqn\label{r-properties}
\nabla_\bx r = 0,\,\,\,
\partial_\bz = \partial_r,\,\,\,
\partial_\bz\nabla_\bx r = 0,\,\,\,
\partial_\bz r = -1,\,\,\,
\partial_\bz^2 r = 0.
\Eeqn

We now proceed to prove Theorem \ref{Stability Estimate} by
estimating the derivatives of $U$ at $p$.
First, by \eqref{weingarten.graph}, $\CB_{\dt}$ can be expressed 
in terms of $g$ via the following formula,
\Beqn\label{A}
\CB_{\dt} =  \left ( I - \dfrac{ \nabla_{\bar{x}} g \otimes \nabla_{\bar{x}} g}{\sqrt{1 + | \nabla_{\bar{x}} g | ^2}} \right ) \dfrac{ \nabla_{\bar{x}}^2 g}{\sqrt{1 + | \nabla_{\bar{x}} g | ^2}}.
\Eeqn
To analyze the above, we need the following
expressions which are obtained by differentiating $g$ 
using the implicit function
(recall that $U(\bx, g(\bx), \dt) = 0$):
\begin{eqnarray}
\nabla_{\bar{x}} g (\bar{x} )
& = & - \Lover{\partial_{\bar{z}}U (\bar{x}, g(\bar{x}), \dt )}
\nabla_{\bar{x}} U (\bar{x}, g (\bar{x}), \dt ),
\label{1st deriv g}
\\
\nabla^2_{\bar{x}} g(\bar{x})
&= &
- \Lover{(\partial_{\bar{z}} U)} \big( \nabla^2_{\bar{x}} U + \partial_{\bar{z}} \nabla_{\bar{x}} U \otimes \nabla_{\bar{x}} g \big) \nonumber\\
& & +\Lover{(\partial_{\bar{z}}U)^2} \big(\partial_{\bar{z}}\nabla_{\bar{x}}U \otimes \nabla_{\bar{x}} U + \partial_{\bar{z} \bar{z}} U \nabla_{\bar{x}} g \otimes \nabla_{\bar{x}} U \big)\Big|_{(\bar{x},g(\bar{x}),\dt)}.
\label{2nd deriv g}
\end{eqnarray}

Theorem \ref{Stability Estimate} is proved by obtaining precise estimates for 
$\nabla_{\bar{x}}g, \nabla^2_{\bx}g$ and hence $\CB_{\dt}$.
The central technical part is the analysis of the term
$\nabla_{\bar{x}}^2 U^1$.

We will frequently use the following estimates concerning the Green's 
function: for some constant $C$ that depend only on the spatial dimension,
it holds that
\Beqn\label{GreensL1Estimates}
\Norm{G(\cdot, t)}_{L^1(\R^{n+1})} = 1\,\,\,(\leq C)
\,\,\,\text{and}\,\,\,
\Norm{\nabla G(\cdot, t)}_{L^1(\R^{n+1})} \leq \frac{C}{\sqrt{t}}.
\Eeqn

\subsection{Estimates for terms without $\nabla^2_{\bar{x}}U^1$}\label{easy estimates}
This section gives several preliminary estimates for various
derivatives associated with $U$.

\begin{thm4}[Estimates for $U^0$]
\label{gradU0} 
\begin{eqnarray}
\nabla_{\bar{x}} U^0(p) 
& = & 0,\label{U0x}\\
\nabla_{\bar{x}}^2 U^0(p) 
& = & - \frac{1}{\sqrt{\pi \dt}} \exp\left(- \frac{g^2(0)}{4\dt}\right) 
\nabla^2_\bx r(p)
\label{U0xx}\\
\partial_\bz\nabla_\bx U^0(p) & = & 0,
\label{U0zx}\\
\partial_{\bar{z}} U^0(p)
& = & -\frac{1}{\sqrt{\pi \dt}} \exp\left(- \frac{g^2(0)}{4\dt}\right),
\label{U0z}\\ 
|\partial^2_{\bar{z}} U^0(p)| 
& \lesssim & \|\CB_0\|^2.
\label{U0zz}
\end{eqnarray}
\end{thm4}
\begin{proof}
Recall that
\MathSty{U^0 = \frac{2}{\sqrt{\pi}}  
\int^{ \frac{r}{2\sqrt{t}} }_0 \exp(-y^2) dy}
and $r(x,t)=\text{sdist}(x,M_t)$.
Then we have,
\begin{eqnarray*}
\nabla_{\bx} U^0 
& = & 
\Lover{\sqrt{\pi t}}\exp\Brac{-\frac{r^2}{4t}}\nabla_{\bx}r,\\
\nabla_{\bar{x}}^2 U^0 
&=&
\frac{1}{\sqrt{\pi t}} \left[ 
\exp\left(- \frac{r^2}{4 t}\right) 
\nabla_{\bar{x}}^2 r + \frac{r}{2 t} \exp\left(- \frac{r^2}{4 t}\right) 
\nabla_{\bar{x}}r \otimes \nabla_{\bar{x}}r 
\right],\\
\partial_\bz\nabla_{\bar{x}} U^0
& = & 
\Lover{\sqrt{\pi t}}\exp\Brac{-\frac{r^2}{4t}}\partial_\bz\nabla_{\bx}r
-\Lover{\sqrt{\pi t}}\exp\Brac{-\frac{r^2}{4t}}\frac{r}{2t}
\partial_\bz r\nabla_{\bx}r,
\\
\partial_{\bar{z}} U^0 
&=& \frac{1}{\sqrt{\pi t}} \exp \left( - \frac{r^2}{4 t} \right) 
\partial_\bz r,\\
\partial^2_{\bar{z}} U^0 
&=& \frac{1}{\sqrt{\pi t}} \SqrBrac{
\exp\left(- \frac{r^2}{4 t}\right) \partial_{\bz}^2 r
- \exp\left(-\frac{r^2}{4t}\right)\Brac{\frac{r}{2t}}(\partial_\bz r)^2
}
\end{eqnarray*}
Noting the property \eqref{r-properties}, 
and the facts $r(p) = g(0)$, 
and $\|g\|_{L^\infty}\lesssim \|\CB_0\|^2 \dt^{\frac{3}{2}}$
(Theorem \ref{consistency}), all the claims in the Lemma follow.
\end{proof}

\begin{thm4}[Estimates for $U^1$]\label{L infty grad U^1} 
\begin{eqnarray}
\| \nabla U^1 \|_{L^{\infty}(T_{\rho}(\Gamma_0))} 
&\lesssim& \| \CB_0 \|^2 \sqrt{\dt},
\label{gradU1.simple}\\
\| \nabla^2 U^1 \|_{L^{\infty}(T_{\rho}(\Gamma_0))} 
&\lesssim& \|\CB_0\|^2.
\label{grad2U1.simple}
\end{eqnarray}
\end{thm4}
\noindent
(In the above, the gradient operator $\nabla$ is defined with respect 
to the spatial variable $x\in\R^{n+1}$.)
\begin{proof}
The statements are consequences of the 
$L^1$-estimates of the Green's function \eqref{GreensL1Estimates}.

Taking the spatial gradient in (\ref{U^1}), we get,
\begin{equation*}
\nabla U^1 (x, \dt) 
= \int_0^{\dt} \int_{T_{\rho}(\Gamma_0)} 
\Big(\nabla_y G(x -y , \dt - \tau)\Big)
\left(\frac{1}{\sqrt{\pi \tau}}\exp\left(-\frac{r^2(y,\tau)}{4 \tau} \right) 
r(y,\tau)\psi(y,\tau) \right)
dy d\tau.
\end{equation*}
Therefore for all $x \in T_{\rho}(\Gamma_0)$,
\begin{align*}
\big | \nabla U^1 (x, \dt) \big | 
&\lesssim
\int_0^{\dt} 
\Norm{\nabla G(\cdot, \dt-\tau)}_{L^1(\R^{n+1})}
\Norm{
\frac{r(\cdot,\tau)}{\sqrt{\tau}}\exp\left(-\frac{r^2(\cdot,\tau)}{4 \tau} 
\right) 
\psi(\cdot,\tau)
}_{L^\infty}\,d\tau
\\
&\lesssim  \int_0^{\dt} \frac{\Norm{A_\tau}^2}{\sqrt{\dt-\tau}} d\tau
= \| \CB_0 \|^2 \Cal{O}(\sqrt{\dt})
\end{align*}
which is \eqref{gradU1.simple}.
(Note that we have used \eqref{psi-exp-r-est}.)

Differentiate (\ref{U^1}) again to get,
\begin{eqnarray*}
& & \left |\nabla^2 U^1 (x, \dt)\right |\\
&\leq& \int_0^{\dt} \int_{T_{\rho}(\Gamma_0)}  \left| \nabla_y G(x -y , \dt - \tau) \right | \left| \nabla_y \left(\frac{r(y,\tau)}{\sqrt{\pi \tau}}\exp\left(-\frac{r^2(y,\tau)}{4 \tau} \right) \psi(y,\tau) \right) \right| dy d\tau.
\end{eqnarray*}
Note that for $f(y) = ye^{-\frac{y^2}{4}}$, we have
\begin{eqnarray*}
\nabla\SqrBrac{f\Brac{\frac{r}{\sqrt{\tau}}}\psi(y,\tau)}
= f'\Brac{\frac{r}{\tau}}\frac{\nabla\tau}{\sqrt{\tau}}\psi
+ f\Brac{\frac{r}{\tau}}\nabla\psi
\end{eqnarray*}
so that
\begin{eqnarray*}
\Abs{\nabla\SqrBrac{f\Brac{\frac{r}{\sqrt{\tau}}}\psi(y,\tau)}}
\lesssim 
\frac{\Norm{\psi}_{L^\infty}}{\sqrt{\tau}}
+ \Norm{\nabla\psi}_{L^\infty}
\lesssim 
\frac{\Norm{A_\tau}^2}{\sqrt{\tau}} + \Norm{A_\tau\nabla A_\tau}
\lesssim
\frac{\Norm{\CB_0}^2}{\sqrt{\tau}}.
\end{eqnarray*}
Hence
\begin{eqnarray*}
\Norm{\nabla^2 U^1(\cdot,\dt)}_{L^\infty}
& \lesssim & \|\CB_0\|^2 \int_0^{\dt} \frac{1}{\sqrt{\dt - \tau}} \frac{1}{\sqrt{\tau}} d \tau\,\,\,\,\, 
\text{(by Lemma \ref{regularity of curvature} and \eqref{GreensL1Estimates})}\\
& \lesssim & \|\CB_0\|^2
\end{eqnarray*}
which is \eqref{grad2U1.simple}.
\end{proof}

\begin{thm4}[First derivative estimates for $U$ and $g$]\label{assorted_terms} 
\begin{eqnarray}
\big( \partial_{\bar{z}}U \big)^{-1}(p)
&=& \left( \frac{1}{\sqrt{\pi\dt}} \exp(-\frac{g^2}{4\dt}) 
\right)^{-1} \Big ( 1 + \| \CB_0 \|^2 \Cal{O}(\dt) \Big ),\label{partial_z U}
\\
|\nabla_{\bar{x}}g(0)|
&\lesssim&
\| \CB_0 \|^2 \dt.\label{grad g estimate}
\end{eqnarray}
\end{thm4}
\begin{proof}
First note that
$\partial_{\bar{z}}U(p) = \partial_\bz U^0 + \partial_\bz U^1$.
Then by \eqref{U0z} and \eqref{gradU1.simple}, we have
\begin{eqnarray*}
\partial_{\bar{z}}U(p) 
& = & \frac{1}{\sqrt{\pi\dt}} \exp\Brac{-\frac{g(0)^2}{4\dt}} 
+ \| \CB_0 \|^2 \Cal{O}(\sqrt{\dt})\\
& = & \frac{1}{\sqrt{\pi\dt}} \exp\Brac{-\frac{g(0)^2}{4\dt}} 
\Brac{1 + \sqrt{\dt}\exp\Big(\frac{g(0)^2}{4\dt}\Big)
\| \CB_0 \|^2 \Cal{O}(\sqrt{\dt})
}
\\
& = & 
\frac{1}{\sqrt{\pi\dt}} \exp\Brac{-\frac{g(0)^2}{4\dt}} 
\Big ( 1 + \| \CB_0 \|^2 \Cal{O}(\dt) \Big )
\end{eqnarray*}
leading to \eqref{partial_z U}.
Note that we have used the estimate 
$\Abs{g(0)}\lesssim\CO(\dt^\frac{3}{2})$
from Theorem \ref{consistency}.

For \eqref{grad g estimate}, the statement follows from 
\eqref{1st deriv g}, \eqref{U0x}, \eqref{partial_z U}, 
and \eqref{gradU1.simple}:
\begin{eqnarray*}
\Abs{\nabla_x g}
= \Abs{\frac{\nabla_\bx U^0 + \nabla_\bx U^1}{\partial_z U}}
&=& \Abs{\Norm{\CB_0}^2\CO(\sqrt{\dt})
\Brac{\Lover{\sqrt{\dt}}\exp\Brac{-\frac{g(0)^2}{4\dt}}}^{-1}
(1+\Norm{\CB_0}\CO(\dt))
}\\
&=&\Abs{\Norm{\CB_0}^2\CO(\dt)
\exp\Brac{\frac{g(0)^2}{4\dt}}(1+\Norm{\CB_0}^2\dt)}\\
& = & \Norm{\CB_0}^2\CO(\dt).
\end{eqnarray*}
\end{proof}

\subsection{Refined estimates for $\nabla^2_{\bar{x}} U^1$}
\label{RefinedU12}
We now begin the most substantial computation in estimating 
$A_{\dt}$,
which is directly related to $\nabla^2_{\bar{x}} U^1$.  
The following is the key result of this section. 
\begin{thm4}\label{U1-2ndDerEst}
The following estimate holds,
\Beqn
\Norm{\nabla_{\bar{x}}^2 U^1} \lesssim \|\CB_0\|^3 \sqrt{\dt}.
\Eeqn
\end{thm4}

Before proceeding with the proof of Lemma \ref{U1-2ndDerEst}, 
we show how Lemma \ref{U1-2ndDerEst} along with the estimates of Section \ref{easy estimates} can lead to Theorem \ref{Stability Estimate}.
First, by substituting the estimates from Lemma \ref{L infty grad U^1} into
\eqref{2nd deriv g}, we have
\begin{eqnarray*}
& & \nabla^2_{\bar{x}} g(\bar{x})\Big|_{(\bar{x},g(\bar{x}),\dt)}\\
&= &
- \Lover{(\partial_{\bar{z}} U)} \big( \nabla^2_{\bar{x}} U + \partial_{\bar{z}} \nabla_{\bar{x}} U \otimes \nabla_{\bar{x}} g \big)
+\Lover{(\partial_{\bar{z}}U)^2} \Big(\partial_{\bar{z}}\nabla_{\bar{x}}U \otimes \nabla_{\bar{x}} U + 
(\partial^2_{\bz} U) \nabla_{\bar{x}} g \otimes \nabla_{\bar{x}} U \Big)
\\
& = & 
- \Lover{(\partial_{\bar{z}} U)} 
\big( \nabla^2_{\bar{x}} U^0 + \nabla^2_{\bar{x}} U^1 
+ \partial_{\bar{z}} \nabla_{\bar{x}} U^1 \otimes \nabla_{\bar{x}} g \big)
\\
& & 
+\Lover{(\partial_{\bar{z}}U)^2} \Big(\partial_{\bar{z}}\nabla_{\bar{x}}U^1 
\otimes \nabla_{\bar{x}} U^1 
+ (\partial^2_\bz U^0 + \partial^2_\bz U^1) 
\nabla_{\bar{x}} g \otimes \nabla_{\bar{x}} U \Big)\\
& = & 
-\Lover{(\partial_\bz U)}\SqrBrac{
-\Lover{\sqrt{\pi\dt}}\exp\Brac{-\frac{g(0)^2}{4\dt}}
\nabla^2_\bx r(p)
+\nabla^2_{\bar{x}} U^1 + \Norm{\CB_0}^2\Norm{\CB_0}^2\dt
}\\
& &
+\Lover{(\partial_{\bar{z}}U)^2} \Big(
\Norm{\CB_0}^2\Norm{\CB_0}^2\sqrt{\dt} + 
\Norm{\CB_0}^2\Norm{\CB_0}^4\CO(\dt^\frac{3}{2})
\Big)\\
& = & 
-\SqrBrac{\Lover{\sqrt{\pi\dt}}\exp\Brac{-\frac{g(0)^2}{4\dt}}}^{-1}\SqrBrac{
-\Lover{\sqrt{\pi\dt}}\exp\Brac{-\frac{g(0)^2}{4\dt}}
\nabla^2_\bx r(p)
+\nabla^2_{\bar{x}} U^1 + \Norm{\CB_0}^4\CO(\dt)
}\\
& & 
+ \SqrBrac{\Lover{\sqrt{\pi\dt}}\exp\Brac{-\frac{g(0)^2}{4\dt}}}^{-2}
\SqrBrac{
\Norm{\CB_0}^4\CO(\sqrt{\dt}) + \Norm{\CB_0}^6\CO(\dt^\frac{3}{2})
}\\
& = & 
\nabla^2_\bx r(p)
+ \nabla^2_{\bar{x}} U^1\CO(\sqrt{\dt}) 
+ \Norm{\CB_0}^4\CO(\dt^\frac{3}{2})\\
& = & 
\nabla^2_\bx r(p)
+ \Norm{\CB_0}^3\CO(\dt) 
+ \Norm{\CB_0}^4\CO(\dt^\frac{3}{2}).
\end{eqnarray*}

Now apply the above to \eqref{A} which relates $\CB_\dt$ to $\nabla_\bx^2 g$, 
we get
\begin{eqnarray*}
| \CB_{\dt} |
& = & 
\frac{\Abs{\text{tr} \big ( \nabla_{\bar{x}}^2 g\big )}}
{\sqrt{1+\Abs{\nabla_\bx g}^2}}\Abs{
I - \frac{\nabla_\bx g\otimes\nabla_\bx g}{\sqrt{1+\Abs{\nabla_\bx g}^2}}
}\\
& \lesssim &
\Big(\Abs{\Lap_\bx r(p)} 
+ \Norm{\CB_0}^3\CO(\dt)
+ \Norm{\CB_0}^4\CO(\dt^\frac{3}{2})\Big)
\Brac{1+C\Norm{\CB_0}^2\CO(\dt^2)}\\
& \lesssim &
\Big(\Abs{H_{\CurBrac{r=r(p)}}(p)} 
+ \Norm{\CB_0}^3\CO(\dt)
+ \Norm{\CB_0}^4\CO(\dt^\frac{3}{2})\Big)
\Brac{1+C\Norm{\CB_0}^2\CO(\dt^2)}\\
&\lesssim&
\Big(\Norm{A_\dt} +\Norm{\CB_0}^3\CO(\dt)
+ \Norm{\CB_0}^4\CO(\dt^\frac{3}{2})\Big)
\Brac{1+C\Norm{\CB_0}^2\CO(\dt^2)}\\
&\lesssim&
\Big(\Norm{\CB_0}(1 + C\Norm{\CB_0}^2\CO(\dt))
+ \Norm{\CB_0}^3\CO(\dt)
+ \Norm{\CB_0}^4\CO(\dt^\frac{3}{2})\Big)
\Brac{1+C\Norm{\CB_0}^2\CO(\dt^2)}\\
& \lesssim &
\Norm{\CB_0}\Brac{1+C\Norm{\CB_0}^2\CO(\dt)}
\end{eqnarray*}
which is the statement of Theorem \ref{Stability Estimate}.
In the above, $\Lap_\bx r(p)=H_{\CurBrac{r=r(p)}}(p)$ 
is the mean curvature of the manifold $\CurBrac{r=r(p,\dt)}$ 
at $p$. It is estimated in terms of $\Norm{A_\dt}$ by using 
\eqref{H.dist.form}:
\[
\Abs{H_{\CurBrac{r=r(p,\dt)}}(p)}
= \Abs{\sum_i\frac{\kappa_i(\Pi(p))}{1-r(p,\dt)\kappa_i(\Pi(p))}}
\leq \Norm{A_\dt}\Brac{1 + \Norm{A_\dt}\CO(\dt^\frac{3}{2})}
\]
and then in terms of $\Norm{\CB_0}$ by 
using Lemma \ref{curvature bound}.

The proof of \Bf{Lemma \ref{U1-2ndDerEst}} is divided 
into several sub-sections. 

\subsubsection{Preparation for Estimating $\nabla_{\bar{x}}^2 U^1(p,\dt)$,
$p\in\Gamma_\dt$}  

Evaluating $U^1(x,t)$ at $x=p = (0,g(0))$ and $t=\dt$, we get,
\begin{align*}
U^1(p,\dt) 
&= - \int_0^{\dt} \int_{T_{\rho}(\Gamma_0)} 
G(p - y, \dt - \tau) \psi( y,\tau) \dfrac{r(y, \tau)}{\sqrt{\pi \tau}}
\exp\left(-\frac{r^2(y, \tau)}{4\tau}\right) dy d\tau. 
\end{align*}
For the integration in the $y$ variable, we use the coordination system 
\eqref{coord_s0} (see also Figure 2):
\Beqn\label{coord_s0_2}
\R^{n+1} 
\cong \R^n\times\R 
\cong \CurBrac{(y_1, y_2, \ldots y_n)\in\R^n}\times\CurBrac{y_{n+1}\in\R}.
\Eeqn
Taking the second partial derivatives in tangential directions 
$\bx_i$ and $\bx_j$ for $1 \leq i, j\leq n$ leads to
\begin{multline}\label{2nd deriv U^1}
\partial^2_{\bx_i \bx_j} U^1(p,\dt)\\
= - \int_0^{\dt} \int_{T_{\rho}(\Gamma_0)} 
\Big(\partial_{y_j} G(p - y, \dt - \tau)\Big)
\partial_{y_i}\SqrBrac{\psi( y,\tau) \dfrac{r(y, \tau)}{\sqrt{\tau}}
\exp\left(-\frac{r^2(y, \tau)}{4\tau}\right)} dy d\tau.
\end{multline}
In addition, for 
$\left | x - y \right | \geq \sqrt{2(n+5) (\dt - \tau)
\left| \log (\dt -  \tau) \right |}$, we have
\[
\left |\nabla G(x - y, \dt - \tau) 
\right | \leq \frac{1}{2 \left(4 \pi \right)^{\frac{n+1}{2}}}
\left( \dt - \tau \right )^\frac{3}{2} \sqrt{(\dt - \tau)\left| \log (\dt - \tau ) \right |} = o(\dt).
\]
Hence we can restrict our analysis of (\ref{2nd deriv U^1}) to the 
following region 
\Beqn\label{restricted-region}
\Big\{y\in\R^{n+1}:
| p - y | \leq \sqrt{2( n+5) \dt | \log \dt | }
\Big\}.
\Eeqn

Now we compute:
\begin{eqnarray}
& & \partial_{y_i}\SqrBrac{ \psi(y,\tau) \dfrac{r(y, \tau)}{\sqrt{\pi\tau}}
\exp\Brac{-\frac{r^2(y, \tau)}{4\tau}}}\nonumber\\
&=&
\psi( y , \tau )\Brac{1 - \frac{r^2(y, \tau)}{2\tau}}
\exp\Brac{-\frac{r^2(y, \tau)}{4\tau}}
\frac{\partial_{y_j}r}{\sqrt{\pi\tau}}\label{I.def}\\
& &  +
\dfrac{r(y, \tau)}{\sqrt{\pi\tau}} \exp\Brac{-\frac{r^2(y, \tau)}{4\tau}}
\partial_{y_i}\psi(y,\tau)\label{J.def}
\end{eqnarray}
Then we decompose (\ref{2nd deriv U^1}) into the following
two terms:
\[
\partial^2_{x_i x_j} U^1(0, g(0)) = -I -J
\]
where
\begin{eqnarray}
I & = & 
- \int_0^{\dt} \int_{T_{\rho}(\Gamma_0)} \partial_{y_j} \Big ( G(p - y, \dt - \tau) \Big)\BigBrac{\text{expression \eqref{I.def}}}\,dyd\tau,
\label{I.int}\\
J & = & 
- \int_0^{\dt} \int_{T_{\rho}(\Gamma_0)} \partial_{y_j} \Big ( G(p - y, \dt - \tau) \Big)\BigBrac{\text{expression \eqref{J.def}}}\,dyd\tau.
\label{J.int}
\end{eqnarray}
Observe that $I$ does not involve any derivatives of the curvature term, 
while $J$ does. 

The main estimates we will arrive at are:
\[
\Abs{I} \lesssim \|\CB_0\|^3\sqrt{\dt}
\,\,\,\,\,\,\text{and}\,\,\,\,\,\,
\Abs{J} \lesssim \|\CB_0\|^3\sqrt{\dt}
\]
which together will give the result of \Bf{Lemma \ref{U1-2ndDerEst}}.

\subsubsection{Analysis of $I$ \eqref{I.int}, 
term without derivatives of curvature}
Note that 
\Beqn
\partial_{y_i} r(y,\tau) = \InnProd{\N(M_\tau, \Pi_{M_\tau}y)}{\E_i}
\Eeqn
where 
\begin{eqnarray*}
\Pi_M(y) & = & \text{closest point on $M$ to $y$,}\\
\N(M,q) & = & \text{(outward) normal vector to $M$ at $q\in M$.}\\
\E_i & = & \text{the $i$-th coordinate vector of $\T_{M_\dt}(p_0)$.}
\end{eqnarray*}
We will also use the following abbreviation,
\[
\Pi_\tau(y) = \Pi_{M_\tau}(y),
\,\,\,\,\,\,
\N_\tau(y) = \N(M_\tau, \Pi_{M_\tau}(y)).
\]

The main observation is that 
$\N_\tau(y) \approx \N_\dt(p)$ which is
orthogonal to $\E_i$ for $i=1,2,\ldots n$.  
This is made precise as follows.
First extend the normal vector to $M_\tau$ as a vector field over 
$T_{\rho}(\Gamma_0)$ by taking it to be constant in directions normal to 
$M_\tau$. Let $\nabla$ be the gradient operator on $\mathbb{R}^{n+1}$.  
Then 
\begin{eqnarray*}
\N_\tau(y)
& = &
\N_\dt(p) 
+ \int_0^1 \nabla\N_\tau(y + \lambda(p-y))(p-y)\,d\lambda
- \int_\tau^\dt \frac{d}{ds}\N_s(y)\,ds\\
& = & \N_\dt(p) + 
\mathcal{O}(1) \| \CB_0\| 
\left( (\dt - \tau)^{\frac{1}{2}} + | p - y | \right)
\,\,\,\,\,\,\text{(by Lemma \ref{regularity of curvature})}
\end{eqnarray*}
which leads to
\Beqn
\Big|\InnProd{\N_\tau(y)}{\E_i}\Big|
\lesssim
\| \CB_0\| \left( (\dt - \tau)^{\frac{1}{2}} + | p - y | \right).
\label{NormalTaylor}
\Eeqn

Substituting the above into $I$ \eqref{I.int}, we get
\begin{eqnarray*}
\Abs{I}
& = & 
\left|\int_0^{\dt}  \int_{T_{\rho} (\Gamma_0)}
\Big(\partial_{y_j} G(p - y, \dt - \tau) \Big)
\psi( y , t )\times\right.\\
& & 
\left.\times
\Big( 1 - \frac{r^2(y, \tau)}{2\tau} \Big) 
\exp\left(-\frac{r^2(y, \tau)}{4\tau}\right)
\frac{\InnProd{\N_\tau}{\E_i}}{\sqrt{\tau}}
\,dy\,d\tau
\right|
\\
& \lesssim & 
\Abs{
\int_0^\dt\int_{T_\rho(\Gamma_0)}
\Abs{\nabla_y G(p-y, \dt-\tau)}
\Norm{\psi}_{L^\infty}
\frac{
\| \CB_0\| \left( (\dt - \tau)^{\frac{1}{2}} + | p - y | \right)}
{\sqrt{\tau}}
\,dy\,d\tau
}\\
& \lesssim & \| \CB_0 \|^3 
\Abs{
\int_0^\dt\int_{T_\rho(\Gamma_0)}
\Abs{\nabla_y G(p-y, \dt-\tau)}
\frac{
\left( (\dt - \tau)^{\frac{1}{2}} + | p - y | \right)}
{\sqrt{\tau}}
\,dy\,d\tau
}.
\end{eqnarray*}

Now for
\begin{eqnarray*}
& & \Abs{
\int_0^\dt\int_{T_\rho(\Gamma_0)}
\Abs{\nabla_y G(p-y, \dt-\tau)}
\frac{(\dt - \tau)^{\frac{1}{2}}}{\sqrt{\tau}}
\,dy\,d\tau
}\\
& \lesssim &
\Abs{
\int_0^\dt
\Norm{\nabla G(\cdot, \dt - \tau)}_{L^1(\R^{n+1})}
\frac{(\dt - \tau)^\Lover{2}}{\sqrt{\tau}}
\,dy\,d\tau}
\lesssim 
\int_0^\dt \frac{1}{\sqrt{\tau}}\,d\tau
= \CO(\sqrt{\dt}),
\end{eqnarray*}
while for
\begin{eqnarray*}
& & 
\Abs{
\int_0^\dt\int_{T_\rho(\Gamma_0)}
\Abs{\nabla_y G(p-y, \dt-\tau)}
\frac{\abs{p-y}}{\sqrt{\tau}}
\,dy\,d\tau
}\\
& \lesssim &
\Abs{
\int_0^\dt\int_{\R^{n+1}}
\Lover{(\dt-\tau)^{\frac{n+1}{2}}}
\exp\left(-\frac{\abs{y}^2}{4(\dt-\tau)}\right)
\frac{\abs{y}}{\dt -\tau}
\frac{\abs{y}}{\sqrt{\tau}}
\,dy\,d\tau
}\\
& \lesssim &
\Abs{
\int_0^\dt\int_{\R^{n+1}}
\exp\left(-\frac{\abs{y}^2}{4}\right)
\abs{y}^2
\frac{1}{\sqrt{\tau}}
\,dy\,d\tau
}
\lesssim 
\Abs{\int_0^\dt\Lover{\sqrt{\tau}}\,d\tau}
= \CO(\sqrt{\dt})
\end{eqnarray*}
which combined together gives the estimate for $I$ stated in 
\eqref{I.int}.

\subsubsection{Analysis of $J$ \eqref{J.int}, 
term involving derivatives of curvature}
Recall the formula \eqref{psi-formula} for the function $\psi$
where the curvatures are evaluated at $\Pi_{\tau}(y)$.
Then we have 
\begin{eqnarray*}
\partial_{y_j}\psi
& = &
\sum_l^n \partial_{y_j}\Brac{\frac{\kappa_l^2}{1-r\kappa_l}}\\
& = & 
2\sum_{l=1}^n \frac{\kappa_l \partial_{y_i} \kappa_l}{ 1 - r \kappa_l } + 
\sum_{l=1}^n \frac{ \kappa_l^2 }{ (1 - r \kappa_l)^2} 
\Big( r \partial_{y_i} \kappa_l + \kappa_l\partial_{y_j} r\Big)\\
& = & 2\sum_{l=1}^n \kappa_l \partial_{y_i} \kappa_l 
+ \sum_{l=1}^n 
\left[ 
r \left( \frac{ 2\kappa_l^2 \partial_{y_i}\kappa_l}{1- r \kappa_l} 
+ \kappa_l^2 \frac{\partial_{y_i} \kappa_l}{(1 - r \kappa_l)^2} \right)+ 
\frac{ \kappa_l^3 \langle \N_\tau(y), \E_i\rangle}{ (1 - r \kappa_l)^2} 
\right] \\
& = & 2\sum_{l=1}^n \kappa_l \partial_{y_j} \kappa_l + \mathcal{O}(1)\|\CB_0\|^3 \left[ \frac{r}{\sqrt{t}} + \|\CB_0\|( | x_0 - y | + | \dt - t |^{\frac{1}{2}}) 
\right],
\end{eqnarray*}
where (\ref{NormalTaylor}) and Lemma \ref{regularity of curvature} 
have been used. Substituting this back into \eqref{J.int}, we obtain,
\begin{eqnarray}
- J 
&=& 2\int_0^{\dt} \int_{T_{\rho}(\Gamma_0)} 
\partial_{y_j} G(p - y, \dt - \tau)\left[ \dfrac{r(y, \tau)}{\sqrt{\tau}}\exp\left(-\frac{r^2(y, \tau)}{4\tau}\right) 
\sum_{l=1}^n \kappa_l \partial_{y_i} \kappa_l\right] dy d\tau
\label{J1.int}
\\
& & + \int_0^{\dt} \int_{T_{\rho}(\Gamma_0)} \partial_{y_j} 
G(p - y, \dt - \tau) \left[\dfrac{r(y, \tau)}{\sqrt{\tau}}\exp\left(-\frac{r^2(y, \tau)}{4\tau}\right)\times\right.\nonumber\\
& & 
\hspace{80pt}\left.
\times\mathcal{O}(1)\|\CB_0\|^3 \Brac{\frac{r(y,\tau)}{\sqrt{\tau}} 
+ \|\CB_0\|\Brac{ | p- y | + | \dt - \tau |^{\frac{1}{2}}}}
\right] dy d\tau \label{J2.int}\\
& :=& J_1 + J_2.\nonumber
\end{eqnarray}
Note that both $J_1$ and $J_2$ contains derivatives of $\kappa_l$'s. 
However, $J_2$ is easier to deal with as it contains the pre-factor $r$
which makes the integrand small. Hence it is analyzed first.

\Bf{Analysis of $J_2$ \eqref{J2.int}.}
Notice that,
\begin{eqnarray*}
&& \dfrac{r(y, \tau)}{\sqrt{\tau}}\exp(-\frac{r^2(y, \tau)}{4\tau}) 
\|\CB_0\|^3 \left[ \frac{r(y,\tau)}{\sqrt{\tau}} 
+ \|\CB_0\|( | p - y | + | \dt - \tau |^{\frac{1}{2}}) \right] \\
& \lesssim & 
\frac{r^2}{\tau} \exp\Brac{-\frac{r^2}{4\tau}}
\|\CB_0 \|^3 
+ \| \CB_0 \|^4 \frac{r}{\sqrt{\tau}} \exp\Brac{-\frac{r^2}{\sqrt{4\tau}}} 
( | p - y | + | \dt - \tau |^{\frac{1}{2}}) \\
& \lesssim & 
\|\CB_0\|^3 + \| \CB_0 \|^4 ( | p - y | + | \dt - \tau |^{\frac{1}{2}}).
\end{eqnarray*}
Hence
\begin{eqnarray*}
| J_2 | 
& \lesssim & 
\int_0^\dt\int_{\R^{n+1}}\Abs{\nabla G(\cdot, \dt-\tau)}
\SqrBrac{
\|\CB_0\|^3 + \| \CB_0 \|^4 ( | p - y | + | \dt - \tau |^{\frac{1}{2}})
}\\
& \lesssim & 
\|\CB_0 \|^3 \sqrt{ \dt} + 
\|\CB_0 \|^4 \dt 
\lesssim 
\|\CB_0 \|^3 \sqrt{ \dt}.
\end{eqnarray*}

\Bf{Analysis of $J_1$ \eqref{J1.int}.}
First note that the following simple bound for $J_1$
\begin{eqnarray*}
\Abs{J_1}
&\lesssim&
\Abs{\int_0^{\dt} \int_{T_{\rho}(\Gamma_0)}
\partial_{y_j} G(p - y, \dt - \tau)\left[ \dfrac{r(y, \tau)}{\sqrt{\tau}}\exp\left(-\frac{r^2(y, \tau)}{4\tau}\right)
\sum_{l=1}^n \kappa_l \partial_{y_i} \kappa_l\right] dy d\tau}\\
& \lesssim &
\Abs{
\int_0^\dt
\Norm{\nabla G(\cdot, \dt-\tau)}_{L^1(\R^{n+1})}
\frac{\Norm{\CB_0}^2}{\sqrt{\tau}}
\,dy\,d\tau
}\\
& \lesssim & 
\Norm{\CB_0}^2\int_0^\dt\Lover{\sqrt{\dt-\tau}}\Lover{\sqrt{\tau}}\,d\tau
\lesssim
\Norm{\CB_0}^2
\end{eqnarray*}
is too crude for our purpose. The following more refined computation
takes into account the different scalings of the integrand along the 
tangential and normal directions to $M_\tau$.

For this purpose, we will use the co-area formula to perform the 
integration in $J_1$ over $\CurBrac{r=r'}$ for $-\rho \leq r' \leq \rho$:
\begin{multline*}
J_1 = 2\int_0^{\dt} \int_{-\rho}^\rho
\int_{\{r=r'\}} 
\Big(\partial_{y_j} G(p - y, \dt - \tau) \Big) 
\Big(\dfrac{r'(y,\tau)}{\sqrt{\tau}}\exp(-\frac{r'^2(y,\tau)}{4\tau}) \Big) 
\Big(\sum_{l=1}^n \kappa_l \partial_{x_i} \kappa_l\Big)_{\Pi_{\tau}(y)}
\times\\
\times dA_{r'}(y) dr' d\tau.
\end{multline*}
In the above, we recall the sign distance function $r$ 
\eqref{sign.dist} to $M_\tau$ and the notation 
\eqref{manifold.at.distance}
$\CurBrac{r=r'} = \CurBrac{y\in\R^{n+1}: r(y,\tau)=r'}$.
Furthermore, $dA_{r'}$ is the volume form of $\CurBrac{r=r'}$.
Then we parametrize the collection of manifolds 
$\Big\{\CurBrac{r=r'}: -\rho \leq r' \leq \rho\Big\}$ 
using $M_\tau$ as follows:
\Beqn\label{tilde.s.def}
\tilde{s}:\,\,\,
M_\tau\times(-\rho, \rho) \MapTo \tilde{s}(s,r') := s + r'\N(s,\tau)
\Eeqn
where $\N(s,\tau)$ is the unit normal to $M_\tau$ at $s$. 
Note that $s=\Pi_{\tau}(\tilde{s})$.
Then
\begin{multline*}
J_1  =
2\int_0^{\dt} \int_{-\rho}^\rho \int_{M_\tau} 
\partial_{y_j} \Big ( G(p-\tilde{s}, \dt - \tau) \Big)
\Big(\sum_{l=1}^n \kappa_l(s, \tau) \partial_{y_i} \kappa_l(s,\tau)\Big)
\dfrac{r'}{\sqrt{\tau}}\exp(-\frac{r'^2}{4\tau})\times\\
\times\Abs{\frac{dA_{r'}(\tilde{s})}{dA_0(s)}}dA_0(s)\, dr'\, d\tau
\end{multline*}
where $dA_0$ is the volume form of $M_\tau$ and 
$\Abs{\frac{dA_{r'}(\tilde{s})}{dA_0(s)}}$ is Jacobian for
the change of variable from $\tilde{s}\in\CurBrac{r=r'}$ to
$s\in M_\tau$.

Next we decompose 
\[
J_1 = J_{11} + J_{12}
\]
where
\begin{multline}
J_{11} = 
2\int_0^{\dt} \int_{-\rho}^\rho \int_{M_\tau} 
\partial_{y_j} \Big ( G(p-\tilde{s}, \dt - \tau) \Big)
\Big(\sum_{l=1}^n \kappa_l(s, \tau) \partial_{y_i} \kappa_l(s,\tau)\Big)
\dfrac{r'}{\sqrt{\tau}}\exp(-\frac{r'^2}{4\tau})\times\\
\times dA_0(s)\, dr'\, d\tau,
\end{multline}
and
\begin{multline}
J_{12} = 
2\int_0^{\dt} \int_{-\rho}^\rho \int_{M_\tau} 
\partial_{y_j} \Big ( G(p-\tilde{s}, \dt - \tau) \Big)
\Big(\sum_{l=1}^n \kappa_l(s, \tau) \partial_{y_i} \kappa_l(s,\tau)\Big)
\dfrac{r'}{\sqrt{\tau}}\exp(-\frac{r'^2}{4\tau})\times\\
\times\SqrBrac{\Abs{\frac{dA_{r'}(\tilde{s})}{dA_0(s)}}-1}
dA_0(s)\, dr'\, d\tau.
\end{multline}
In a sense, $J_{12}$ is the error term due to the deviation of the
Jacobian from $1$. Hence $J_{11}$ is the dominating term.
The remainder of this subsection is to show that both 
$| J_{11} |$ and $| J_{12} |$ can be bounded from above by 
$C \Norm{\CB_0}\sqrt{\dt}$. 

To proceed, we decompose $\nabla G$ into two exponential 
kernels, corresponding to integrating in the tangential and normal directions 
to $M_\tau$. For this purpose, we introduce the operator
$\Sigma_{s,\tau}(p)$ which gives the orthogonal projection of $p$ onto
$\text{Span}\CurBrac{\N(s,\tau)}$, the normal line at $s\in M_\tau$.
(See Figure 4.)
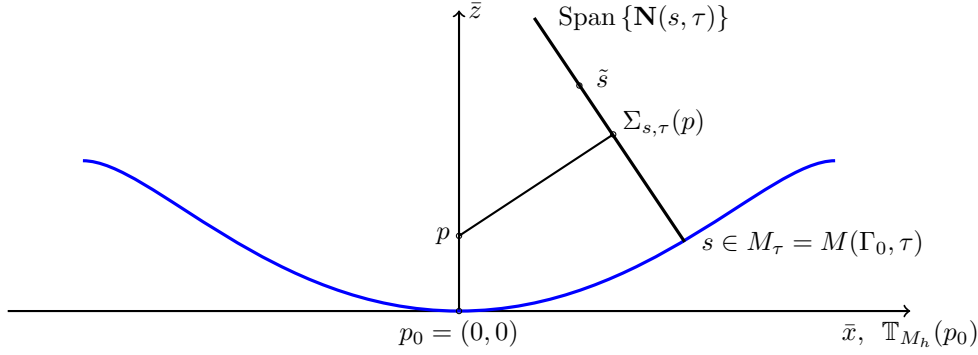
\begin{figure}[h]
\begin{center}
\begin{tikzpicture}
\draw[->, thick] (-6,0) -- (6,0);
\draw[->, thick] (0,0) -- (0,4);
\draw[blue, very thick] (-5, 2) .. controls (-4.1, 2) and (-2.5, 0) .. (0, 0);
\draw[blue, very thick] (0,0) .. controls (2.5,0) and (4.1, 2) .. (5, 2);
\draw (6,0) node[below]{$\bar{x},\,\,\,\mathbb{T}_{M_\dt}(p_0)$};
\draw (0,4) node[right]{$\bar{z}$};
\draw (0,0) node[below]{$p_0=(0,0)$};
\draw (0,0) circle (1pt);
\draw (0,1) circle (1pt);
\draw (0,1) node[left]{$p$};
\draw (2.05, 2.35) circle (1pt);
\draw[very thick, -] (3,0.92) -- node[midway, right]{}(1, 3.9);
\draw[thick, -] (0, 1) -- node[midway, right]{}(2.05, 2.35);
\draw (1.6,3) circle (1pt);
\draw (1.7,3.1) node[right]{$\tilde{s}$};
\draw (2.05, 2.5) node[right]{$\Sigma_{s,\tau}(p)$};
\draw (3.1, 0.90) node[right]{$s\in M_\tau = M(\Gamma_0,\tau)$};
\draw (1.2, 3.9) node[right]{$\text{Span}\CurBrac{\N(s,\tau)}$};
\end{tikzpicture}
\caption{Illustration of $\Sigma_{s,\tau}$. Note that for any $\tilde{s}
\in T_\rho(\Gamma_0)$, 
it can be decomposed as $\tilde{s} = s + (\tilde{s} - s)$ where
$s=\Pi_\tau(\tilde{s})\in M_\tau$ and 
$\tilde{s}-s\in\text{Span}(\N(s,\tau))$.}
\end{center}
\end{figure}

Then for any $\tilde{s}=s+r'\N(s,\tau)\in\CurBrac{r=r'}$, we write,
\Beqn
p-\tilde{s} 
= \big[p-\Sigma_{s,\tau}(p)\big] + \big[\Sigma_{s,\tau}(p) - \tilde{s}\big].
\Eeqn
Note that
\begin{eqnarray*}
\partial_{y_j}G(p-y,\dt-\tau) 
= \Lover{(\dt-\tau)^{\frac{n+1}{2}}}
\exp\Brac{-\frac{\abs{p-y}^2}{4(\dt-\tau)}}
\frac{\InnProd{p-y}{\E_j}}{\dt-\tau}.
\end{eqnarray*}
Hence
\begin{eqnarray*}
& & \partial_{y_j}G(p-\tilde{s},\dt-\tau) \\
&=& 
\SqrBrac{ \frac{\langle p - \Sigma_{s,\tau}(p), \E_j \rangle}
{(\dt - \tau)^{\frac{n+1}{2} + 1}} 
+ \frac{\langle \Sigma_{s,\tau}(p) - \tilde{s}, \E_j \rangle}
{(\dt - \tau)^{\frac{n+1}{2} + 1}}}
\exp\Brac{-\frac{|p - \Sigma_{s,\tau}(p)|^2}{4(\dt - \tau)}}
\exp\Brac{-\frac{|\Sigma_{s,\tau}(p) - \tilde{s}|^2}{4(\dt - \tau)}}.
\end{eqnarray*}

Substituting the above into the expression for $J_{11}$, we have
\[
J_{11} = J_{111} + J_{112}
\]
where 
\begin{eqnarray*}
J_{111} &=& \int_0^{\dt} \int_{M_\tau} 
\frac{\langle p - \Sigma_{s,\tau}(p), \E_j \rangle}
{(\dt - \tau)^{\frac{n}{2} + 1} } 
\exp\Brac{-\frac{|p - \Sigma_{s,\tau}(p)|^2}{4(\dt - \tau)}} 
\Brac{\sum_{l=1}^n \kappa_l \partial_{y_i} \kappa_l}_{(s,\tau)}\times \\
&&\times \SqrBrac{\Lover{(\dt - \tau)^{\frac{1}{2}}}  \int_{r'} 
\exp\Brac{-\frac{|\Sigma_{s,\tau}(p) - \tilde{s}|^2}{4(\dt - \tau)}}
\frac{r'}{\sqrt{\tau}}\exp(-\frac{r'^2}{4\tau}) dr'}
dA_0(s)\,d\tau, \\
J_{112} &=&
 \int_0^{\dt} \int_{M_\tau} 
\Lover{(\dt - \tau)^{\frac{n}{2}}} 
\exp\Brac{-\frac{|p - \Sigma_{s,t}(p)|^2}{4(\dt - \tau)}} 
\Brac{\sum_{l=1}^n \kappa_l \partial_{x_i} \kappa_l}_{(s,\tau)}\times \\
&&\times \SqrBrac{\Lover{(\dt - \tau)^{\frac{1}{2}}}
\int_{r'} \frac{\langle \Sigma_{s,t}(p) - \tilde{s}, \E_j \rangle}
{(\dt - \tau) }
\exp\Brac{-\frac{|\Sigma_{s,t}(p) - \tilde{s}|^2}{4(\dt - \tau)}}
\frac{r'}{\sqrt{\tau}}\exp(-\frac{r'^2}{4\tau}) dr'}dA_0(s)\, d\tau.
\end{eqnarray*}
which essentially correspond to integrations along the tangential and 
normal directions.

\Bf{Analysis of $J_{111}$.}
We first evaluate the inner integration with respect to $r'$.
From \Bf{Section \ref{Integrations}}, we have that,
\begin{multline}
\Lover{(\dt - \tau)^{\frac{1}{2}}}
\int_{-\infty}^\infty 
\exp\Brac{-\frac{|\Sigma_{s,\tau}(p) - \tilde{s}|^2}{4(\dt - \tau)}}
\frac{r'}{\sqrt{\tau}}\exp\Brac{-\frac{r'^2}{4\tau}} dr'\\
= 
2\sqrt{\pi} \frac{\tau}{\dt^{\frac{3}{2}}}r(\Sigma_{s,t}(p), \tau) 
\exp\Brac{-\frac{r^2(\Sigma_{s,\tau}(p), \tau)}{4\dt}}.
\label{Tricky Integral 1}
\end{multline}
Substituting this into $J_{111}$, we get
\begin{multline}\label{J_{11}}
J_{111} =   \frac{\sqrt{\pi}}{(\dt)^{\frac{3}{2}}} \int_0^{\dt} 
\tau \left[\int_{M_\tau} 
\frac{\langle p - \Sigma_{s,\tau}(p), \E_j \rangle}
{(\dt - \tau)^{\frac{n}{2} + 1} }
\exp\Brac{-\frac{|p - \Sigma_{s,\tau}(p)|^2}{4(\dt - \tau)}} 
\Brac{\sum_{l=1}^n \kappa_l \partial_{y_i} \kappa_l}_{(s,\tau)}
\right.\times\\
\left.\times r(\Sigma_{s,\tau}(p), \tau)
\exp\Brac{-\frac{r^2(\Sigma_{s,\tau}(p), \tau)}{4\dt}}dA_{0}(s)
\right]d\tau.
\end{multline}

Heuristically, the inner integration on $M_\tau$ is essentially
a convolution with the $n$-dim Green's function integrated 
on the tangent plane $\T_{M_\tau}(p_0)$ to $M_\tau$ at 
$p_0=\Pi_{M_\tau}(p)$. Recall the coordinate system \eqref{coord_s0_2}.
Note also that by (\ref{2nd deriv U^1}), we may restrict our attention to
$ |y| \leq \sqrt{2(n+5)\dt | \log \dt | } \ll 1$.
Now we change the integration variable from $s$ to $y$.
The following estimates can be established,
\begin{eqnarray}
\exp\Brac{ -\dfrac{|p - \Sigma_{s, \tau}(p)|^2}{4(\dt - \tau)}}
& \lesssim & 
\exp \left( -\dfrac{C |y|^2}{ \dt - \tau } \right)
\label{Est1}\\
| p - \Sigma_{s, \tau}(p) | 
& \lesssim & | y|,
\label{Est2}\\
dA_0
& \lesssim & d^n y.
\label{Est3}
\end{eqnarray}
We further note the following two statements:
\begin{enumerate}
\item Let $f$ be the graph of $M_\tau$ over $\T_{M_\tau}(p_0)$. Then
\begin{eqnarray*}
\big| r(\Sigma_{s,\tau}(p), \tau) \big| 
& \lesssim & \Abs{r( p, \tau)} + \Abs{\langle y, \nabla f \rangle}
\lesssim \Abs{r( p, \tau)} + \Abs{y}\Abs{\nabla f}\\
& \lesssim & \| \CB_0 \|^2 \dt^{\frac{3}{2}} + \|\CB_0\|^2 \dt \Abs{y}\\
& \lesssim & \| \CB_0 \|^2 \dt^{\frac{3}{2}} + \|\CB_0\|^2 (\dt^2 +\Abs{y}^2).
\end{eqnarray*}

\item By the regularity estimate Lemma \ref{regularity of curvature}
for MMC, we have
\begin{align}\label{Est5}
| \nabla_{y_i} \kappa_l | 
\lesssim \frac{\| B_\tau \|}{\sqrt{\tau}}
\lesssim \frac{\| B_0 \|}{\sqrt{\tau}}.
\end{align}
\end{enumerate}

Now substituting (\ref{Est1}) - (\ref{Est5}) into (\ref{J_{11}}), 
we obtain,
\begin{multline*}
|J_{111}| 
\lesssim \frac{1}{\dt^{\frac{3}{2}}} \int_0^{\dt} \tau 
\int_{y} \frac{|y|}{(\dt - \tau)^{\frac{n}{2} + 1} } 
\exp\Brac{-\frac{3 |y|^2}{4(\dt - \tau)}} 
\frac{\|\CB_0\|^2}{\sqrt{\tau}}\times\\
\times\Big \{ \| \CB_0 \|^2 (\dt)^{\frac{3}{2}} +  \|\CB_0\|^2 \big( \dt^2 + | y |^2\big)  \Big \} d^n y d\tau.
\end{multline*}
Using the fact that
$
\int_{\R^n}
\frac{| y |^K}{(\dt - \tau)^{n/2 + 1} } 
\exp\Brac{-\frac{3 | y |^2}{4(\dt - \tau)}}
\,d^n y
\lesssim
(\dt - \tau)^{\frac{K}{2} - 1}$,
we have
\begin{eqnarray*}
&&| J_{111} | \\
&\lesssim &
\frac{\Norm{\CB_0}^4}{\dt^{\frac{3}{2}}} \int_0^{\dt} 
\frac{\sqrt{\tau}\dt^{\frac{3}{2}}}{\sqrt{\dt - \tau}} d\tau
+ \frac{\Norm{\CB_0}^4}{\dt^{\frac{3}{2}}} \int_0^{\dt} 
\dfrac{\sqrt{\tau}\dt^2}{ \sqrt{\dt - \tau} }d\tau 
+ \frac{\Norm{\CB_0}^4}{\dt^{\frac{3}{2}}} \int_0^{\dt} \sqrt{\tau} 
\sqrt{\dt - \tau} d\tau \\
&\lesssim& 
\| \CB_0 \|^4 \dt + + \|\CB_0\|^4\dt^\frac{3}{2} + \| \CB_0 \|^4 \sqrt{\dt}\\
&\lesssim& \| \CB_0 \|^3 \sqrt{\dt}.
\end{eqnarray*}

\Bf{Analysis of $J_{112}$.}
This is very similar to that for $J_{111}$, but with slightly different 
terms. The inner integral in $r'$ is given by 
(see \Bf{Section \ref{Integrations}}):
\begin{multline}
\Lover{\sqrt{\dt - t}} 
\int_{-\infty}^\infty 
\frac{\InnProd{\Sigma_{s,\tau}(p) - \tilde{s}}{\E_j}}
{(\dt - \tau) }
\exp\Brac{-\frac{|\Sigma_{s,\tau}(p) - \tilde{s}|^2}{4(\dt - \tau)}}
\frac{r'}{\sqrt{\tau}}\exp(-\frac{r'^2}{4\tau}) dr'\\ 
= 2\sqrt{\pi} \InnProd{\N(s,\tau)}{\E_j}
\exp\Brac{ -\frac{ r^2(\Sigma_{s,\tau}(p), \tau) }{ 4\dt } }
\Brac{ \frac{ \tau }{\dt^{\frac{5}{2}}} r^2(\Sigma_{s,\tau}(p), \tau) 
- \frac{2\tau}{\dt^{\frac{3}{2}}} 
}.
\label{Tricky Integral 2}
\end{multline}
As $\frac{a^2}{b} \exp(-\frac{a^2}{b} )$ 
is a bounded function of $a$ and $b$, the above can be estimated as
\[
\mathcal{O}(1) \frac{\tau}{(\dt)^{\frac{3}{2}}} 
\InnProd{\N(s,\tau)}{\E_j}.
\]
Hence 
\begin{align*}
| J_{112} | \lesssim \frac{\Norm{\CB_0}^2}{(\dt)^{\frac{3}{2}}} 
\int_0^{\dt} 
\sqrt{\tau} \int \Lover{(\dt - \tau)^{\frac{n}{2}}}
\exp \Brac{ - \frac{ | y|^2 }{4(\dt - \tau)} }
\Big|\InnProd{\N(s,\tau)}{\E_j}\Big| d^ny d\tau,
\end{align*}

Similar to \eqref{NormalTaylor}, we can infer that
\begin{align*}
\InnProd{\N(s, \tau)}{\E_j}
= \mathcal{O}(1)\|\CB _0\| \Big(\sqrt{\dt} +  |s| \Big)
= \mathcal{O}(1)\|\CB _0\| \Big(\sqrt{\dt} +  |y| \Big)
\,\,\,\,\,\,\text{(note $s=(y,h(y))$).}
\end{align*}
Substituting this into the integral and using the $L^1-L^{\infty}$ 
estimate for the Green's function, we obtain,
\begin{align*}
| J_{112} | 
\lesssim & 
\frac{\Norm{\CB_0}^3}{(\dt)^{\frac{3}{2}}} 
\int_0^{\dt} \sqrt{\tau}  
\int \Lover{(\dt - \tau)^{\frac{n}{2}}}
\exp \Brac{- \frac{ C |y|^2 }{\dt - \tau} } 
\big( \sqrt{\dt} + |y|\big)  d^ny d\tau\\
\lesssim &
\frac{\Norm{\CB_0}^3}{(\dt)^{\frac{3}{2}}} 
\int_0^{\dt} \sqrt{\tau}  
\Brac{\sqrt{\dt} + \sqrt{\dt - \tau}}\,d\tau
\lesssim \| \CB_0 \|^3 \sqrt{\dt}
\end{align*}
completing the analysis for $J_{11}$.

\Bf{Analysis of $J_{12}$.}
Note that \MathSty{\Abs{\Abs{\frac{dA_{r'}(\tilde{s})}{dA_0(s)}}-1}
\lesssim \Norm{\CB_0}\Abs{r'}}.
Hence $J_{12}$ is of smaller order. Using similar analysis as for
$J_{11}$, we can then conclude that
\[
\Abs{J_{12}} \lesssim \Norm{\CB_0}^3\sqrt{\dt}.
\]

The above concludes the proof of Theorem \ref{Stability Estimate}.

\begin{thm6}\label{improved_convergence_rate}
Note that in the proof of Theorem \ref{Stability Estimate}, we establish
point-wise estimates for $\nabla^2 U^1$ of the form 
$ \|\CB_0\|^3 \CO(\sqrt{\dt})$.  
Applying the same point-wise analysis to $U^1$ and $\nabla U^1$, 
we could have established that 
$U^1 = \|\CB_0 \|^3 \CO((\dt)^{\frac{3}{2}})$ and 
$\|\nabla U^1\| = \|\CB_0\|^3 \CO(\dt)$.  
This will be improvements over the estimates
\eqref{est1} and \eqref{gradU1.simple} which follow from
standard $L^{\infty}$-estimates given by convolving 
$(\partial_t - \Delta)U^0$ with the Green's function $G$.
\end{thm6}

\subsection{Stability over successive iterations}
In this section, we prove that the algorithm can be iterated over
$\lfloor \frac{T}{\dt} \rfloor$ steps, over which the numerical 
manifold stays embedded and has a uniform curvature bound.  
There are two tools we use toward this. 
The first is a ``Ball Lemma'' which can ensure the
embeddedness of the numerically manifolds $\Gamma_{n\dt}$'s. 
The second is a discrete Gronwall-type inequality as given in 
\cite{cheung2004some}. We state the latter first in the following.
\begin{thm1}\label{Gronwall}
(\cite[Theorem 2.1]{cheung2004some})
Let $D$ be a constant. Suppose $\{ x_k \}_{k\geq 0}$ 
is a sequence of numbers satisfying for all $k\geq 1$ that
\[
| x_k | \leq |x_{k-1} | \big ( 1 + D |x_{k-1} |^2 \dt).
\]
Let further $\Phi$ be the following monotone (and hence invertible)
function,
\[
\Phi(x) = \int_1^x \dfrac{ ds }{s + D s^3 }.
\]
Then for all $k \geq 0$, we have,
\[
| x_k | \leq \Phi^{-1} \Big( \Phi ( | x_0 | ) + k \dt \Big).
\]
\end{thm1}

Next is the statement of the Ball Lemma.
\begin{thm4}[Ball Lemma]\label{BallLemma} 
(For one step.)
Suppose we have a constant $C_0$ such that 
$\| \CB_0 \|, \|\CB_\dt\| \leq C_0$.
Then $\Gamma_{\dt}$ satisfies the Ball Property 
(Definition \ref{BallProperty}) with radius $r$ given by,
\Beqn\label{r1}
r = \min \left( \frac{1}{C_0+1},\,\,m_0 - C(n, C_0)\dt \right),
\Eeqn
where 
\Beqn
m_0 = \min\left\{
\Abs{s_1 - s_2}:
s_1, s_2\in\Gamma_0\subset \R^{n+1},\,\,
d_0(s_1, s_2) \geq \frac{C_*}{C_0} 
\right\}, 
\Eeqn
$d_0$ is the distance between points on 
$\Gamma_0$ defined in \eqref{DistOnManifold},
$C_*$ is the universal constant in \eqref{RadiusGraph}
and $C(n,C_0)$ is some constant depending
only on the spatial dimension $n$ and the curvature bound $C_0$.

(For multiple steps.)
Suppose we have iterated the algorithm $N$ times and that 
$\| \CB_{k\dt} \| \leq C_0$ for $0 \leq k \leq N$ for some $N$.  Then $\Gamma_{N\dt}$ will satisfy the ball property with radius
\Beqn\label{rm}
r = \min \left( \frac{1}{C_0 + 1}, m_0 - C(n, C_0) N \dt \right).
\Eeqn
\end{thm4}

The proof is given in \Bf{Section \ref{ProofBallLemma}}.
\begin{thm6}
It is evident in the proof of the above lemma that $M(\Gamma_0, t)$, for $0 \leq t \leq \dt$, will also satisfy the Ball Property with the same radius 
$r$ as in \eqref{r1}.
\end{thm6}

We now show how Lemma \ref{BallLemma} and Theorem \ref{Gronwall} 
allow for the algorithm to be iterated repeatedly
with uniform curvature estimate.

\subsubsection{\bf{Proof of Theorem \ref{Finite Time Stability}}}
Let $\tilde{T} = \Phi(C_0) - \Phi(\Norm{\CB_0})$ which is also the maximal
$t$ such that
$\Phi^{-1} \left ( \Phi ( \| \CB_0 \| ) + t \right ) \leq C_0$.
Now define
\Beqn
T := 
\min \left \{ \tilde{T},\,\,\,\dfrac{m_0 - 2 \rho}{ C(n, C_0) } \right \}
=
\min \left \{\Phi(C_0) - \Phi(\Norm{\CB_0}),\,\,\,
\dfrac{m_0 - 2 \rho}{ C(n, C_0) } 
\right \}
\Eeqn
where $C(n, C_0)$ is the constant coming from Lemma \ref{BallLemma}.

From Theorem \ref{Stability Estimate}, we have
$$
\| \CB_{\dt} \| \leq \| \CB_0 \| \left( 1 + C\|\CB_0\|^2 \dt \right).
$$
Applying Theorem \ref{Gronwall} we have that,
$$
\| \CB_{\dt} \| \leq \Phi^{-1} \left(  \Phi( \| \CB_0 \| ) + \dt  \right) \leq C_0.
$$
Consequently, the Ball Lemma implies that $\Gamma_{\dt}$ is an embedded manifold, which
satisfies the ball property with radius,
$$
r = \min \left(  \frac{1}{C_0 + 1}, m_0 - C(n,C_0) \dt \right).
$$
Notice that $\rho \leq \frac{1}{2} \min \left(  \frac{1}{C_0 + 1}, m_0 - C(n,C_0) \dt \right)$.
As a consequence, we may repeat all the preceding analysis with 
$\Gamma_{0}$ and $T_{\rho}(\Gamma_{0})$ replaced by
$\Gamma_{\dt}$ and $T_{\rho}(\Gamma_{\dt})$.
Iterating for $k$ ($0\leq k\leq\lfloor \frac{T}{\dt} \rfloor$) steps, 
we arrive at the result.

\begin{thm6}
The above proof reveals the simultaneous preservation of uniform curvature 
estimates and the embeddedness of the numerical manifolds $\Gamma_{k\dt}$.
The bigger the $C_0$ and smaller the $\rho$ are, 
the larger the maximum convergence time $T$ can be chosen.
The smaller $\rho$ ensures that the Ball Property will hold for longer
time interval. The only constraint is \eqref{InitialDataAssump} which will 
impose a smaller maximum time-step size $\delta$.
\end{thm6}


\section{Convergence to MMC}\label{ProofConvSec}

In this section we prove that the algorithm converges to MMC.
Recall that the convergence is phrased in the weak form using the 
framework bounded variation (BV-)functions or sets of finite perimeter
\eqref{WeakFormulation-a}-\eqref{WeakFormulation-b}
\cite{luckhaus}.
We first explain some key notations.
We refer to \cite{evans-gariepy, Giusti} for more detail exposition
about the function space. 

Let $\chi:\R^{n+1}\MapTo\R$ with $\Norm{\chi}_{L^1(\R^{n+1})} < \infty$. 
It is called a function with bounded variation,
written as $\chi\in BV(\R^{n+1}, \R)$, or simply $\chi\in BV(\R^{n+1})$,
if its variational derivative is given by a (Radon) measure.
Precisely, there is a finite constant $C$ such that
\Beqn
\sup\CurBrac{\int \chi\text{div}\varphi: 
\varphi\in C^1_0(\R^{n+1}, \R^{n+1}),\,\Abs{\varphi} \leq 1}
\leq C.
\Eeqn
In this case, there is a Borel measure $\mu$ and a vector valued function
$\nu: \text{supp}(\mu)\MapTo \mathbb{S}^{n}$ such that
for any $\varphi\in C^1_0(\R^{n+1}, \R^{n+1})$,
\Beqn
\int\chi\text{div}\varphi = -\int\varphi\cdot\nu\,d\mu.
\Eeqn
It is also customary to denote $\mu$ and $\nu$ by $\Abs{\nabla\chi}$ and 
\MathSty{\NORMAL{\chi}}.

In the present paper, $\chi$ takes its values from $\CurBrac{0,1}$ so that
we write $\chi \in BV(\R^{n+1}, \CurBrac{0,1})$. 
Then the set
$\Omega=\CurBrac{x: \chi(x)=1}$ (so that $\chi = \mathbbm{1}_\Omega$)
is called a \Em{set of finite perimeter}.
We state here two fundamental facts about such sets.

(i) There exists
a notion of \Em{reduced boundary} $\Gamma=\Bdry^*\Omega$ which is
$\Cal{H}^{n}$-rectifiable such that
$\mu = \Cal{H}^{n}\lfloor_\Gamma$ and $\nu$ is a continuous function 
on $\Gamma$ $\mu$ almost everywhere.
Then $\Abs{\nabla \chi}$ and $\NORMAL{\chi}$
are essentially the area measure and 
outward unit normal vector function of $\Gamma$.
The area (also commonly known as the \Em{perimeter}) 
and the mean curvature of $\Gamma$ are given by
\Beqn
\int_\Gamma\,d\Cal{H}^n
= 
\int \Abs{\nabla \chi}
\,\,\,\,\,\,\text{and}\,\,\,\,\,\,
H = \text{div}\NORMAL{\chi}.
\Eeqn
The mean curvature function $H$ can also be defined via the following 
definition: for any $\zeta\in C^1(\R^{n+1},\R^{n+1})$,
\[
\int H\zeta\cdot\nu\Abs{\nabla\chi}
= \int\Brac{\Text{div}\zeta - \nu\cdot\nabla\zeta\nu}\Abs{\nabla\chi}.
\]

(ii) The space of sets of finite perimeter (and more generally BV-functions)
with uniform bounded perimeters is compact in $L^1$: 
if 
\MathSty{\sup_\alpha
\int_{\R^{n+1}}\Abs{\nabla\chi_\alpha} 
+\int_{\R^{n+1}}\Abs{\chi_\alpha} 
< \infty}, then
there is a subsequence $\alpha_i$ and a $\chi^*\in BV(\R^{n+1})$ such that
\[
\chi_{\alpha_i} \Converge \chi^*
\,\,\,\,\,\,\text{in $L^1_{loc}(\R^{n+1})$}.
\]
In particular, if $\chi_\alpha = \mathbbm{1}_{\Omega_\alpha}$, then 
$\chi^* = \mathbbm{1}_{\Omega^*}$ for some $\Omega^*\subseteq\R^{n+1}$.
Furthermore, the area measure is lower-semicontinuous,
i.e. for all Borel set $B\subset\R^{n+1}$, it holds that
\[
\int_B \Abs{\nabla\chi^*} 
\leq \liminf_i 
\int_B \Abs{\nabla\chi_{\alpha_i}}.
\]

In the following, for simplicity, we will simply use the terminology
BV-function with the understanding that we are dealing exclusively with
sets of finite perimeter. 

Now we will make use of the above formulation to prove the convergence of 
$\chi_\dt$ and identify the equation satisfied by its limit. 
We first recall the definition \eqref{approx1} and 
\eqref{approx2} of $\chi_\dt$, $\Omega_\dt$ and $\Gamma_\dt$. 
Note that now $\chi_\dt: \R^{n+1}\times[0,T]
\MapTo\CurBrac{0,1}$ is a function of both spatial and temporal variables. 
For convenience, we denote $\chi_\dt(t) = \chi_\dt(\cdot,t)$.
Since $\chi_\dt(t)$ is a classical solution of MMC for 
$t\in(k\dt, (k+1)\dt)$, denoting $H_\dt(t)$ to be the mean curvature of
$\Gamma_\dt(t)$, we have for any
$\zeta \in C^{\infty}( \bar{ \Lambda \times [0,T] }, \mathbb{R}^n )$, 
$\zeta = 0$ and
$\xi \in C^{\infty}( \bar{ \Lambda \times [0,T] }, \mathbb{R} )$, $\xi = 0$ on $\partial \Lambda \times [0,T] \cup \Lambda \times \{0\}$ that
\begin{multline}\label{MMC.nonlinear.1}
\int_{k\dt}^{(k+1)\dt} \int_{\Lambda} \big( \text{div} \zeta - \frac{\nabla \chi_{\dt}(t)}{| \nabla \chi_{\dt}(t) |} \nabla \zeta \frac{\nabla \chi_{\dt}(t)}{| \nabla \chi_{\dt}(t) |} \big) | \nabla \chi_{\dt}(t) |\\
=- \int_{k\dt}^{(k+1)\dt} \int_{\Lambda} H_{\dt}(t) \zeta \cdot   \nabla \chi_{\dt}(t),
\end{multline}
and
\begin{multline}\label{MMC.nonlinear.2}
\int_{k\dt}^{(k+1)\dt} \int_{\Lambda} \chi_{\dt}(t) \partial_t \xi + \int_{\Lambda} \chi_{\dt}(k\dt^+) \xi(k\dt) -  \int_{\Lambda}  \chi_{\dt} ((k+1)\dt^-) \xi((k+1)\dt^-)\\
= - \int_{k\dt}^{(k+1)\dt} \int_{\Lambda} H_{\dt}(t)  \xi |\nabla \chi_{\dt}(t) |.
\end{multline}
Summing the above over 
$k = 0, 1, ..., \lfloor \frac{T}{\dt} \rfloor$, we obtain:
\Beqn
\int_0^T \int_{\Lambda} \Brac{\text{div} \zeta - \frac{\nabla \chi_{\dt}(t)}{| \nabla \chi_{\dt}(t) |} \nabla \zeta \frac{\nabla \chi_{\dt}(t)}{| \nabla \chi_{\dt}(t) |} }
| \nabla \chi_{\dt}(t) | 
= - \int_0^T \int_{\Lambda} H_{\dt}(t) \zeta \cdot \nabla \chi_{\dt}(t),
\label{NumericalWeakSolution2}
\Eeqn
and
\Beqn
\int_0^T \int_{\Lambda} \chi_{\dt}(t) \partial_t \xi + \int_{\Lambda} \chi_{\dt}(0) \xi(0) = - \int_0^T \int_{\Lambda} H_{\dt}(t) \xi |\nabla \chi_{\dt}(t) | + E(\dt),
\label{NumericalWeakSolution22}
\Eeqn
where 
\Beqn
E(\dt) = \sum_{k=1}^{k = \lfloor \frac{T}{\dt} \rfloor} \int_{\Lambda} \zeta(k\dt) \Big ( \chi_{\dt}(k \dt^+) - \chi_{\dt}(k\dt^-) \Big )
\Eeqn
is the sum of the ``jump'' errors made between iterations, precisely
at the thresholding steps. The above are the discrete analog of 
\eqref{WeakFormulation-a}-\eqref{WeakFormulation-b}.

By the consistency and stability estimates
Theorem \ref{consistency} and 
\ref{Finite Time Stability}, as $\dt\Converge 0$, we have that
\begin{align*}
\left | E(\dt) \right |  \lesssim 
\sum_{k=1}^{k = \lfloor \frac{T}{\dt} \rfloor} 
\dt^{\frac{3}{2}} \Converge 0.
\end{align*}
Thus (\ref{NumericalWeakSolution2})-(\ref{NumericalWeakSolution22}) is 
``almost'' a solution to \eqref{WeakFormulation-a}-\eqref{WeakFormulation-b}.  
The remaining step is to show that $\chi_{\dt}$ and $H_{\dt}$ exhibit appropriate compactness in $\dt$, so that we may pass to the limit 
$\dt \rightarrow 0$ in 
(\ref{NumericalWeakSolution2})-(\ref{NumericalWeakSolution22}).

With the consistency and stability estimates, we can already conclude that
the sequence of manifolds 
$\CurBrac{\Gamma_{k\dt}: 0 \leq k \lfloor\frac{T}{h}\rfloor}$
converges to some limit in the Hausdorff distance 
$d_\Cal{H}$ \eqref{Hausdorff}. It remains to show that the limit
satisfies the equation of MMC. We find the framework of BV-convergence
as stated in Definition \ref{luckhaus} to be the most convenient.

The outline of proof is as follows. We first show that 
$\{\chi_{\dt}\}$ is compact in 
$L^1\Big(\R^{n+1}\times[0,T], \CurBrac{0,1}\Big)$
and hence has a limit $\chi^*$.  
Next we show that the area measure converges in measure:
$\Abs{\nabla \chi_{\dt}}\rightharpoonup \Abs{\nabla\chi^*}$.
This enables us to prove that both the normal vectors
$\dfrac{\nabla \chi_{\dt}}{| \nabla \chi_{\dt} |}$ and 
the mean curvature $H_{\dt}$ are also convergent.

For convenience, we use $rca(\Lambda\times[0,T])$ to denote the space
of regular Radon measures on $\Lambda\times[0,T]$.

\subsection{$L^1$-Compactness}
The main conclusion in this section is that up to subsequence,
$\chi_\dt$ converges to some $\chi^*$ in 
$L^1(\R^{n+1}\times[0,T], \CurBrac{0,1})$.
Furthermore, for each $t$, $\chi^*(t)$ is in $BV(\R^{n+1}, \CurBrac{0,1})$.
This is a consequence of the Kolmogorov-Riesz-Frechet Theorem 
\cite[Thm. 4.26]{Brezis} together 
with the compactness property of BV-functions.
The following sequence of propositions facilitate the use of this theorem.
\begin{thm5}\label{PerimeterBound}
The perimeters of $\chi_{\dt}(t)$ are uniformly bounded, i.e.
\begin{align*}
\sup_{t \in [0, T] } \int_{\Lambda} | \nabla \chi_{\dt} (t) | & < + \infty.
\end{align*}
\end{thm5}
\begin{proof}
Using the implicit function theorem, $\Gamma_{\dt}(k \dt)$ 
can be parametrized over $X_{\dt}(k \dt^-)$ via a map,
$$
q\in X_\dt(k\dt^-) \MapTo q + g(q) \N(q),
$$
where $\N(q)$ is the unit normal of $X_\dt(k\dt^-)$ at $q$.  
By the consistency and stability estimates, we have that
$g = \mathcal{O}\big(\dt^{\frac{3}{2}}\big)$ and 
$\nabla g = \mathcal{O}(\dt)$.
(Here $\nabla g$ is the gradient computed over 
$\chi_{\dt}(k \dt^-)$.)  
Hence,
\begin{align*}
\int | \nabla \chi_{\dt} (k \dt^+) | 
& = \int | \nabla \chi_{\dt} ( k \dt^- ) | 
\big ( 1 + \Cal{O}(\dt) ).
\end{align*}
Furthermore, as the area decreases through MMC, we have
\begin{align*}
\int | \nabla \chi_{\dt} ( k \dt^- ) | &\leq \int | \nabla \chi_{\dt} \big( (k-1) \dt^+ \big) |.
\end{align*}
By iterating, we obtain,
\begin{align*}
\int | \nabla \chi_{\dt} (k \dt^+) | &\leq (1 + C\dt) \int | \nabla \chi_{\dt} ( 0 ) |.
\end{align*}
\end{proof}

\begin{thm5}\label{SpaceCompactness}
For all $t \in [0,T]$, $\dt > 0$ and $w \in \mathbb{S}^n$,
the following spatial continuity statement holds,
\begin{align*}
\lim_{s \rightarrow 0} \int_{\Lambda} | \chi_{\dt} (t, y + s w) - \chi_{\dt} (t, y) | dy = 0.
\end{align*}
\end{thm5}
\begin{proof} 
This follows from the estimate,
\begin{align*}
\int_{\Lambda} | \chi_{\dt} (t, y + s w) - \chi_{\dt} (t, y) | dy 
\lesssim s\int | \nabla \chi_{\dt} ( t ) |
\end{align*}
and the uniform perimeter bound just proved.
\end{proof}

\begin{thm5}\label{TimeLipschitz}
The collection $\{ \chi_{\dt}(t): t\in[0,T] \}$ 
satisfies the following Lipschitz in time estimate,
\begin{align*}
\int_{\Lambda} | \chi_{\dt} ( t + \epsilon) -  \chi_{\dt} ( t ) | 
&\lesssim \epsilon,
\end{align*}
whenever $\dt \leq \epsilon \leq T - t$.
\end{thm5}
\begin{proof}
We have the following estimates,
\[
\int_{\Lambda} | \chi_{\dt} ( k \dt^+ )  - \chi_{\dt} ( k \dt^- ) | 
\lesssim \dt^{\frac{3}{2}},
\,\,\,\,\,\,\text{and}\,\,\,\,\,\,
\int_{\Lambda} | \chi_{\dt} ( k \dt^- )  - \chi_{\dt} \big( (k-1) \dt^+ \big) | \lesssim \dt.
\]
The first follows by the consistency estimate 
while the second follows from the regularity of $\Gamma_\dt(t)$ which
solves the MMC for $t \in [ (k-1) \dt, k \dt)$.
The result follows by iterating these estimates.
\end{proof}

From the above, the Kolmogorov-Riesz-Frechet Theorem 
\cite[Thm 4.26]{Brezis} implies that there is a 
subsequence such that $\chi_\dt$ convergences to $\chi^*$ in 
$L^1(\R^{n+1}\times[0,T], \CurBrac{0,1})$. By the uniform boundedness of the
perimeters of $\chi_\dt(t)$, compactness of sets of finite perimeters
implies that $\chi^*(t)\in BV(\R^{n+1}, \CurBrac{0,1})$ for almost 
every $t$. The Lipschitz continuity in time implies that this holds 
\Em{for every $t$}. 
In particular, we have a fixed subsequence 
$\dt_i\Converge 0$ such that
\Em{for all} $t$,
\Beqn\label{fixed.sequence.dt}
\chi_{\dt_i}(t)\Converge\chi^*(t)
\,\,\,\,\,\,
\text{in $L^1(\Lambda)$.}
\Eeqn
In the following, the notation $\lim_\dt$ and $\dt\Converge 0$ 
refer to $\dt_i\Converge 0$. In several occasions, this subsequence
will be further refined. 
Hence for simplicity, the subscript $i$ will be omitted.

\subsection{Convergence of Area}

In this section, we will prove that 
$\Abs{\nabla \chi_{\dt}}$ converges weakly to $\Abs{\nabla\chi^*}$
in measure. This is a stronger statement than just 
$\nabla\chi_\dt\WeakConverge\nabla\chi^*$. 
It implies that the area converges: 
\[
\int_{\R^{n+1}\times[0,T]}\Abs{\nabla\chi_\dt} \Converge
\int_{\R^{n+1}\times[0,T]}\Abs{\nabla\chi^*} .
\]
By \cite{luckhaus}, this gives a sufficient condition 
for (\ref{NumericalWeakSolution2})-(\ref{NumericalWeakSolution22}) 
to converge to \eqref{WeakFormulation-a}-\eqref{WeakFormulation-b}.
By the uniform boundedness of the perimeter 
(Proposition \ref{PerimeterBound}) and the
Lebesgue Dominated Convergence Theorem, it suffices
to prove that for each $t\in[0,T]$, 
$\Abs{\nabla\chi_\dt(t)} \WeakConverge \Abs{\nabla\chi^*(t)}$.

The first step toward this goal is the observation that
the normal vectors $\nu_\dt(t)$ to $\Gamma_\dt(t)$
converges \Em{strongly}.
By the Lemma \ref{BallLemma} (Ball Lemma), we may extend 
$\nu_{\dt}(t)$ to be a smooth function
defined on $\mathbb{R}^{n+1}$.  
By the uniform $C^2$-bound of the $\Gamma_\dt(t)$, we can invoke the 
Arzela-Ascoli theorem 
to pick a further subsequence of $\dt$ such that
$\nu_\dt(t)$ converges to a 
$\tilde{\nu}(t)$ in $C^{1,\alpha}(\mathbb{R}^{n+1}; \mathbb{R}^{n+1})$.

We are now ready to prove the main theorem of this section.
\begin{thm1}\label{AreaConvergence}
For every $t \in [0,T]$, $| \nabla \chi_{\dt} (t)| \rightharpoonup | \nabla \chi^*(t)|$ in $rca(\Lambda)$ as $\dt \rightarrow 0$, i.e.
for all open set $B\subseteq \R^{n+1}$ such that
$\Abs{\nabla\chi_\dt(t)}(\Bdry B) = 0$, it holds that
\Beqn
\label{area.conv.B}
\lim_\dt\int_B \Abs{\nabla\chi_\dt(t)} = \int_B \Abs{\nabla\chi^*(t)}.
\Eeqn
\end{thm1}
\noindent
We emphasize that the same sequence $\dt\Converge0$ works for every $t$.
\begin{proof}
We fix a $t\in[0,T]$.
First, by the lower semi-continuity of area under $L^1$-convergence, 
we have,
\begin{align*}
\int_B | \nabla \chi^*(t) | &\leq \liminf_{\dt \rightarrow 0} \int_B | \nabla \chi_{\dt} (t) |.
\end{align*}
Next, let $L(B)$ and a subsequence $\dt_j$ of $\dt$ be such that
\[
L(B) 
= \limsup_\dt\int_B\Abs{\nabla\chi_\dt(t)}
= \lim_{\dt_j}\int_B\Abs{\nabla\chi_{\dt_j}(t)}.
\]
Note that the subsequence can depend on $t$. 
But this does not matter as $t$ is fixed.

Let $\nu^*(t)$ be the normal vector function of $\chi^*(t)$.
By the weak convergence of $\nabla\chi_{\dt_j}(t)$ to 
$\nabla\chi^*$,
it holds that
\Beqn
\nu_{(\dt)_j}(t) | \nabla \chi_{(\dt)_j} (t) | 
= \nabla \chi_{(\dt)_j} (t) \rightharpoonup \nabla \chi^*(t) = \nu^*(t) | \nabla \chi^*(t)|.
\Eeqn
We now compute,
\begin{eqnarray*}
L(B) 
& = & \lim_{j \rightarrow \infty} \int_B | \nabla \chi_{(\dt)_j} (t) |\\
&=& \lim_{j \rightarrow \infty} \int_B \langle \nu_{(\dt)_j}(t), \nu_{(\dt)_j}(t) \rangle | \nabla \chi_{(\dt)_j} (t) | 
= \lim_{j \rightarrow \infty} \int_B 
\langle \nu_{(\dt)_j}(t), \nabla \chi_{(\dt)_j} (t)\rangle\\
&=&\lim_{j \rightarrow \infty} 
\int_B \langle \nu_{(\dt)_j}(t) - \tilde{\nu}(t), \nabla \chi_{\dt_j}\rangle
+ \int_B \langle \tilde{\nu}(t), \nabla \chi_{\dt_j}\rangle
\\
& = & 
\int_B \langle \tilde{\nu}(t), \nabla \chi^*(t) \rangle 
=
\int_B \langle \tilde{\nu}(t), \nu^*(t) \rangle | \nabla \chi^*(t) | \\
& & 
\text{(since 
$\nu_\dt(t)\Converge\tilde{\nu}(t)$ strongly in $C^{1,\alpha}$
and $\nabla\chi_{\dt_j}\rightharpoonup\nabla\chi^*$)}\\
& \leq & \int_B | \nabla \chi^*(t) |
\,\,\,\,\,\,\text{(since $\abs{\tilde{\nu}(t)}, \abs{\nu^*(t)} \leq 1$).}
\end{eqnarray*}
Hence 
\[
\limsup_\dt\int_B\Abs{\nabla\chi_\dt(t)}
= L(B) 
\leq \int_B|\nabla \chi^*(t)|
\leq
\liminf_\dt\int_B\Abs{\nabla\chi_\dt(t)}
\]
from which \EqnRef{area.conv.B} follows.
\end{proof}

\subsection{Convergence of $\chi_{\dt}$ to 
(\ref{WeakFormulation-a})-(\ref{WeakFormulation-b}) }
The procedure now largely follows the ideas of the proof of 
\cite[Theorem 2.3]{luckhaus}. Two additional technical but crucial results will 
be used.

The first is that we can control the normal vectors appropriately
so that they can be passed to the limit.  
For this purpose, consider smooth functions $\nu_\epsilon$
which approximate $\nu^*$ in the sense that
\begin{align}\label{normal.approx}
\lim_{\epsilon\converge0}\int_0^T \int_{\Lambda} | \nu^* - \nu_{\epsilon} |^2 |\nabla \chi^*| dt = 0.
\end{align}
Then we claim that
\Beqn\label{normal.error.conv}
\lim_{\epsilon\converge 0}\lim_{\dt}
\int_0^T \int_{\Lambda} | \nu_{\dt} - \nu_{\epsilon} |^2 | \nabla \chi_{\dt} | dt = 0
\Eeqn
which follows from
\begin{eqnarray*}
& & 
\lim_{\epsilon\converge 0}\lim_{\dt}
\int_0^T \int_{\Lambda} | \nu_{\dt} - \nu_{\epsilon} |^2 | \nabla \chi_{\dt} | dt\\
& = &
\lim_{\epsilon\converge 0}\lim_{\dt}
\int_0^T \int_{\Lambda} | \nabla  \chi_{\dt} | dt+ \int_0^T \int_{\Lambda} | \nu_{\epsilon} |^2 | \nabla \chi_{\dt} |dt - 2 							\int_0^T \int_{\Lambda} \langle \nu_{\dt}, \nu_{\epsilon} \rangle | \nabla \chi_{\dt} | dt\\
& = &
\lim_{\epsilon\converge 0}
\int_0^T \int_{\Lambda} | \nabla \chi^* | dt + \int_0^T \int_{\Lambda} | \nu_{\epsilon} |^2 | \nabla \chi^* | dt - 2 \int_0^T \int_{\Lambda} \langle \nu^* , \nu_{\epsilon} \rangle 						| \nabla \chi^* | dt\\
& = & \lim_{\epsilon\converge 0}
\int_0^T \int_{\Lambda} | \nu^* - \nu_{\epsilon} |^2 | \nabla \chi^*| dt \\
& = & 0.
\end{eqnarray*}

The second is the existence of a limiting mean curvature function.  
For this, notice that the measures $H_{\dt} | \nabla \chi_{\dt} |dt$ 
are uniformly bounded in $rca(\Lambda \times [0,T])$.  
Therefore, we may pass to a further subsequence such that 
$ H_{\dt} | \nabla \chi_{\dt} |dt \rightharpoonup \sigma$, 
for some $\sigma\in rca(\Lambda\times[0,T))$. 
We note the following two facts about $\sigma$:

(i) Since the sequence of mean curvatures $H_{\dt}$ are uniformly 
bounded, $\sigma$ is absolutely continuous with respect to 
$| \nabla \chi^*(t)|dt$. Indeed, let $B$ be such that 
\MathSty{\int_B\Abs{\nabla\chi^*(t)}dt = 0}. Then
\[
\Abs{\int_B H_\dt\Abs{\nabla\chi_\dt(t)}dt}
\lesssim \int_B\Abs{\nabla \chi_\dt(t)}dt
\Converge_{\dt\converge 0}
\int_B\Abs{\nabla \chi^*(t)}dt = 0.
\]
Hence \MathSty{\int_B d\sigma 
= \lim_\dt\int H_\dt\Abs{\nabla\chi_\dt}=0}.

(ii) Let $H^*$ to be the Radon-Nikodym derivative of $\sigma$ 
with respect to $| \nabla \chi^* |dt$. Then 
$H^*$ belongs to $L^2( | \nabla \chi^*|dt) $. 
Indeed, let $\eta$ be some smooth test function. Then,
\begin{align*}
\Abs{\int_0^T \int_{\Lambda} H_{\dt} \eta | \nabla \chi_{\dt} | dt}
&\lesssim \Big( \int_0^T \int_{\Lambda} | \eta |^2 | \nabla \chi_{\dt} | dt \Big)^{\frac{1}{2}}.
\end{align*}
Sending $\dt$ to $0$ we get, 
\begin{align*}
\Abs{\int_0^T \int_{\Lambda} H^* \eta | \nabla \chi^* | dt}
&\lesssim \Big( \int_0^T \int_{\Lambda} | \eta |^2 | \nabla \chi^* | dt \Big)^{\frac{1}{2}}.
\end{align*}
The conclusion follows as a consequence of the Riesz-representation theorem.

We now show the convergence of 
(\ref{NumericalWeakSolution2}) to \eqref{WeakFormulation-a}
as follows. 
Consider the right hand side. 
For any $\zeta\in C^\infty(\R^{n+1}\times[0,T], \R^{n+1})$,
\begin{eqnarray*}
& & 
\limsup_{\dt}
\Big | \int_0^T \int_{\Lambda} H_{\dt} \zeta \cdot \nabla \chi_{\dt} dt - \int_0^T \int_{\Lambda} H^* \zeta \cdot \nabla \chi^* dt \Big |\\
& = &  
\limsup_{\dt}
\Big | \int_0^T \int_{\Lambda}  H_{\dt} \zeta \cdot \nu_{\dt} |\nabla \chi_{\dt} | dt -  \int_0^T \int_{\Lambda} H^* \zeta \cdot \nu^* | \nabla \chi^*| dt \Big | \\
&\leq& 
\limsup_{\dt}
\Big | \int_0^T \int_{\Lambda}  H_{\dt} \zeta \cdot \nu_{\epsilon} |\nabla \chi_{\dt} | dt -  \int_0^T \int_{\Lambda} H^* \zeta \cdot \nu^* | \nabla \chi^* | dt \Big |\\
& & 
+ \limsup_\dt \Big | \int_0^T \int_{\Lambda}  H_{\dt} \zeta \cdot ( \nu_{\dt} - \nu_{\epsilon} ) |\nabla \chi_{\dt} | dt \Big|.
\end{eqnarray*}
The first integral of the last line becomes
\[
\Big | \int_0^T \int_{\Lambda} H^* \zeta \cdot( \nu_{\epsilon} - \nu^* ) | \nabla \chi^*| dt \Big |
\]
which by \eqref{normal.approx}
tends to zero upon taking $\epsilon\Converge 0$.
The second integral also tends to zero due to \EqnRef{normal.error.conv} 
and the ($L^2$-)boundedness of the $H_\dt$.

For the second term of the left hand side of 
(\ref{NumericalWeakSolution2}), we similarly have
\begin{eqnarray*}
& & 
\limsup_\dt
\Big | \int_0^T \int_{\Lambda} \nu_{\dt} \nabla \zeta \nu_{\dt} | \nabla \chi_{\dt} | dt  -  \int_0^T \int_{\Lambda}\nu^* \nabla \zeta \nu^* | \nabla \chi^*|dt \Big |\\
&\leq & 
\limsup_\dt
\Big | \int_0^T \int_{\Lambda} \nu_{\epsilon} \nabla \zeta \nu_{\dt} | \nabla \chi_{\dt} |dt   -  \int_0^T \int_{\Lambda} \nu^* \nabla \zeta \nu^* | \nabla \chi^*|dt \Big | \\ 
&&+ \limsup_\dt
\Big | \int_0^T \int_{\Lambda} (\nu_{\dt} - \nu_{\epsilon} )  \nabla \zeta \nu_{\dt} | \nabla \chi_{\dt} | dt \Big | \\
&\lesssim&
\Big | \int_0^T \int_{\Lambda} (\nu_{\epsilon} - \nu^*) \nabla \zeta \nu^* | \nabla \chi^*| dt \Big |
+ 
\limsup_\dt\Big | \int_0^T \int_{\Lambda} 
\Abs{\nu_{\dt} - \nu_{\epsilon}}^2 
| \nabla \chi_{\dt} | dt \Big |
\end{eqnarray*}
which tends to zero as $\epsilon\Converge 0$, again by 
\eqref{normal.approx} and \eqref{normal.error.conv}.

Finally, the convergence of the first integral of the left hand side of 
\EqnRef{NumericalWeakSolution2} and the whole of 
\EqnRef{NumericalWeakSolution22}
follow from 
the weak-convergence of $\Abs{\nabla\chi_\dt}$ to $\Abs{\nabla\chi^*}$, 
the $L^1$-convergence of $\chi_\dt$ to $\chi^*$ 
and the weak-convergence of $H_\dt\Abs{\nabla \chi_\dt}$ to
$H^*\Abs{\nabla\chi^*}$, respectively.

The above concludes that 
\EqnRef{WeakFormulation-a}-\EqnRef{WeakFormulation-b} 
hold as a limiting equation
of \EqnRef{NumericalWeakSolution2}-\EqnRef{NumericalWeakSolution22}
(with $v$ replaced by $H^*$).

\appendix
\section{Appendix}
\subsection{Properties of the signed distance function}\label{sdist-sec}
Here we give some basic properties of the signed distance function and
prove \Bf{Lemma \ref{HeatOperatorFormLem}}, in particular formula
\eqref{HeatOperatorForm} and \eqref{H-Diff-r}.

Given an open set $\Omega \in \mathbb{R}^{n+1}$ with a smooth 
boundary $M=\Bdry\Omega$, define $\Pi = \Pi_{M}$ to be the projection
operator which maps $x\in\R^{n+1}$ to its closest point on $M$.
By the smoothness assumption of $M$, the map $\Pi$ is well-defined 
in a tubular neighborhood of $M$. 
In this tubular neighborhood, then the signed distance function $r$ to 
$M$ can be expressed as
\begin{equation}
r(x) = -(x - \Pi(x)) \cdot \N(\Pi(x)),
\end{equation}
where $\N(z) = \N_M(z)$ is the outward normal to 
$M$ at $z\in M$.

Many geometrical aspects of $M$ can be recovered
from the signed distance function:
(i) for any $x\in\R^{n+1}$, 
$\nabla r(x) = -\N(\Pi(x))$ and hence $ | \nabla r(x) | = 1 $;
(ii) for any $p\in M$,
$- \nabla_{M}^2 r (p)$ 
is the Weingarten map of $M$ at $p$;
(iii) the mean curvature $H_M$ of $M$ at $p$ is given by
$- \Delta r (p)$.

\Bf{Proof of \eqref{HeatOperatorForm}.}
We consider a general function $f: \R^{n+1}\times\R_+\MapTo\R$ 
in the form
\[
f(x,t) = \bar{f}(r(x,t),t)
\]
where $r$ is the signed-distance function to $M_t = M(\Gamma_0, t)$.
We work in the regime that $r$ and $t$ are small enough so that 
$r$ is a smooth function of its arguments.
Let $(x,t)\in\R^{n+1}\times\R_+$ be an arbitrary point. Then we compute,
\begin{eqnarray*}
& & \partial_t f(x, t) - \Delta_x f(x, t) \\
& = & \partial _t \bar{f} ( r, t) + \partial_r \bar{f}(r,t) 
\partial_t r 
- \partial_{rr}\bar{f}(r,t)\Abs{\nabla_x r}^2 
- \partial_r \bar{f}(r,t) \Delta_x r\\
& = & 
\partial_r \bar{f}(r,t) 
\big(\partial_t r - \Lap_x r \big) 
+ \partial_t \bar{f} (r, t) - \partial_{rr}\bar{f}(r, t)
\end{eqnarray*}
where all the arguments of $r$ are evaluated at $(x,t)$.
Note that:
(i) $\partial_t r = - H_{M_t}(\Pi(x))$
where $H_{M_t}$ is the mean curvature function of $M_t$ and 
$\Pi(x)$ is the projection onto $M_t$; and
(ii) $\Lap_x r = - H_{\curBrac{r=r(x,t)}}(x)$. Then we have
\[
\partial_t f(x, t) - \Delta_x f(x, t)
= \partial_r \bar{f}(r,t) 
\big(H_{\curBrac{r=r(x,t)}}(x)-H_{M_t}(\Pi(x))\big)
+ \partial_t \bar{f} (r, t) - \partial_{rr}\bar{f}(r, t).
\]
Finally, by choosing $\bar{f} = U^0$, the last two terms of the above
vanish altogether as $U^0$ solves the linear heat equation 
\eqref{U0HeatEqn}.

\Bf{Proof of \eqref{H-Diff-r}.}
The statement will follow from a formula for 
$H_{\{ r= r_0\}}$.

Since the time variable $t$ is irrelevant, we simply write $M=M_t$.
We first fix a point $x_0$ and let $r_0 = r(x_0)$.
Then for $\abs{r_0}$ small enough, we can express the manifold 
$\{ r = r_0\}$ as a map over $M$ via 
\[
F: M\MapTo\R^{n+1}: F(y) = y + r_0 \N_M(y).
\]
Note that $x_0\in\CurBrac{r=r_0}$. 
Let $\{ \E_1, ... , \E_n \}$ be an orthnormal basis for the 
tangent plane $\T_M$ of $M$.
Since $M$ and $\{ r  = r_0\}$ share the same normal vector, 
i.e. $\N_{\curBrac{r=r_0}}(x) = \N_M(\Pi(x))$ for all $x\in\CurBrac{r=r_0}$, 
$\E_1,..., \E_n$ is also an orthonormal basis for the 
tangent plane to $\{ r = r_0 \}$, in particular at $x_0$.  
Now let $y_0 = \Pi_M(x_0)$ and we impose that at $y_0$, 
$\{ \E_1, ... , \E_n \}$ corresponds to the 
principal curvatures $\CurBrac{\kappa_1(y_0), ..., \kappa_n(y_0)}$ 
of $M$ at $y_0$.  

Next we compute the second fundamental form $\tilde{A}^i_j$ and 
the Weingarten map $\tilde{A}_{ij}$ of $\CurBrac{r=r_0}$ in the following 
fashion.
\begin{itemize}
\item
The metric $\tilde{g}_{ij}$ of $\CurBrac{r=r_0}$ is given in terms of $F$ as
follows,
\begin{align*}
\tilde{g}_{ij} 
= \Big\langle \partial_i F, \partial_j F \Big\rangle 
= \Big\langle \E_i (1 - r_0 \kappa_i ), \E_j (1 - r_0 \kappa_j) 
\Big\rangle = \delta_{ij} (1 - r_0 \kappa_i )^2.
\end{align*}
Its inverse $\tilde{g}^{ij}$ is then,
\[
( 1 - r_0 \kappa_i )^{-2} \delta_{ij}.
\]
\item
By definition, we have,
\begin{align*}
\tilde{A}_{ij} & = \langle \partial_i N, \partial_j F \rangle = \langle \kappa_i e_i, e_j ( 1 - r_0\kappa_j) \rangle = \delta_{ij} \kappa_i( 1 - r_0 \kappa_i  ).
\end{align*}
Hence the Weingarten map takes the form
\begin{align}\label{WM-r}
\tilde{A}^i_j & = \tilde{g}^{ij} \tilde{A}_{ij} = \delta_{ij} \dfrac{ \kappa_i }{ 1 - r_0 \kappa_i }.
\end{align}
\end{itemize}

Upon taking the trace of $\tilde{A}^i_j$, the mean curvature of 
$\CurBrac{r=r_0}$ at $x_0$ is then given by
\Beqn\label{H.dist.form}
H_{\{ r = r_0 \}}(x_0) 
= 
\left.
\sum_{i = 1}^n \dfrac{ \kappa_i }{ 1 - r_0 \kappa_i }
\right|_{y_0} 
= H_M(y_0) + r_0 
\left.\sum_{i = 1}^n \dfrac{\kappa_i^2}{1 - r_0 \kappa_i}
\right|_{y_0}
\Eeqn
which is \eqref{H-Diff-r}.

\subsection{Some Gaussian Integrations: 
Derivations of \eqref{Tricky Integral 1} and \eqref{Tricky Integral 2}
}\label{Integrations}

We first perform two computations regarding Gaussian integrals.

\begin{thm4}
For $d_0 > 0$, $0 < t < \dt$, we have
\Beqn
\Lover{(\dt - t)^{\frac{1}{2}}}
\int_{-\infty}^{\infty} \exp\left(-\frac{|d_0 - r'|^2}{4(\dt - t)}\right)
\frac{r'}{\sqrt{t}}\exp\left(-\frac{r'^2}{4t}\right) dr'
= 
2\sqrt{\pi}
\exp\left( - \frac{d_0^2}{4\dt} \right) 
\frac{d_0 t}{\dt^{\frac{3}{2}}} 
\label{first integral}
\Eeqn
and
\begin{multline}
\Lover{(\dt - t)^{\frac{1}{2}}}
\int_{-\infty}^{\infty} \frac{d_0 - r'}{(\dt - t) } 
\exp\left(-\frac{|d_0 - r' |^2}{4(\dt - t)}\right)
\frac{r'}{\sqrt{t}}\exp\left(-\frac{r'^2}{4t}\right) dr' \\
= 
2\sqrt{\pi}
\exp \left( - \frac{d_0^2}{4\dt} \right) 
\left( \frac{d_0^2 t}{\dt^{\frac{5}{2}}} 
-  \frac{2t}{\dt^{\frac{3}{2}}} \right)
\label{second integral}
\end{multline}
\end{thm4}
\begin{proof}
We will note the following identities:
\[
\int_{-\infty}^\infty e^{-x^2}\,dx = \sqrt{\pi},\,\,\,\,\,\,
\int_{-\infty}^\infty xe^{-x^2}\,dx = 0,\,\,\,\,\,\,
\int_{-\infty}^\infty x^2e^{-x^2}\,dx = \frac{\sqrt{\pi}}{2},
\]
and 
\[
\frac{(d_0 - r')^2}{4(\dt - t)} + \frac{r'^2}{4t}
= 
\frac{\dt}{4t(\dt - t)}\Brac{r - \frac{d_0 t}{\dt}}^2 + \frac{d_0^2}{4\dt}
\]
so that
\Beqn\label{exp.simplified}
\exp\left(-\frac{ ( d_0 - r' )^2 }{4(\dt - t)} \right) 
\exp\left(\frac{-r'^2}{4t}\right) 
= 
\exp\Brac{-\frac{d_0^2}{4\dt}}
\exp\Brac{-\frac{\dt}{4t(\dt - t)}\Brac{r - \frac{d_0 t}{\dt}}^2}.
\Eeqn

\noindent
\Bf{Proof of (\ref{first integral})}. Using the \eqref{exp.simplified},
we have
\begin{eqnarray*}
& & 
\frac{e^{-\frac{ d_0^2 }{4\dt}}}
{\sqrt{(\dt - t)}}  
\int_{-\infty}^{\infty} \frac{r'}{\sqrt{t}}
\exp\left( - \frac{\dt}{4t (\dt - t) } 
\left(r'-\frac{d_0 t}{\dt}\right)^2 \right) dr'\\
& = & 
\frac{e^{-\frac{ d_0^2 }{4\dt}}}
{\sqrt{ t (\dt - t)}}  
\int_{-\infty}^{\infty} \frac{d_0 t}{\dt}\exp\left( - \frac{\dt}{4t (\dt - t) } 
\left(r'-\frac{d_0 t}{\dt}\right)^2 \right) dr'\\
& = & 
2\sqrt{\pi}
\frac{e^{-\frac{ d_0^2 }{4\dt}}}{\sqrt{ t (\dt - t)}}  
\frac{d_0 t}{\dt}
\frac{\sqrt{t(\dt - t)}}{\sqrt{\dt}}
=
2\sqrt{\pi}
\exp\left( - \frac{d_0^2}{4\dt} \right)
\frac{d_0 t}{\dt^{\frac{3}{2}}}
\end{eqnarray*}
which is (\ref{first integral}).

\noindent
\Bf{Proof of (\ref{second integral})}.
As before, using \eqref{exp.simplified}, we have
\begin{eqnarray*}
& & \frac{e^{- \frac{d_0^2}{4\dt}}}
{(\dt - t)^{\frac{1}{2}}}
\int_{-\infty}^{\infty} \frac{(d_0 - r')r'}{(\dt - t)\sqrt{t}} 
\exp\left( - \frac{\dt}{4t (\dt - t) } 
\left( r' - d_0 \frac{t}{\dt} \right)^2 \right)dr' \\
& = & -\frac{e^{- \frac{d_0^2}{4\dt}}}
{(\dt - t)^{\frac{3}{2}}\sqrt{t}}
\int_{-\infty}^{\infty} 
\left(r' - \frac{d_0 t}{\dt} + \frac{d_0 t}{\dt} - d_0\right)
\left(r' - \frac{d_0 t}{\dt} + \frac{d_0 t}{\dt}\right)
\exp\left( - \frac{\dt}{4t (\dt - t) } 
\left( r' - d_0 \frac{t}{\dt} \right)^2 \right)dr' \\
& = & -\frac{e^{- \frac{d_0^2}{4\dt}}}
{(\dt - t)^{\frac{3}{2}}\sqrt{t}}
\int_{-\infty}^{\infty} 
\left[
\left(r' - \frac{d_0 t}{\dt}\right)^2
+ \frac{d_0^2}{\dt^2}(t-\dt)t
\right]
\exp\left( - \frac{\dt}{4t (\dt - t) } 
\left( r' - d_0 \frac{t}{\dt} \right)^2 \right)dr' \\
& = & -\frac{e^{- \frac{d_0^2}{4\dt}}}
{(\dt - t)^{\frac{3}{2}}\sqrt{t}}
\int_{-\infty}^{\infty} 
\left[r'^2 + \frac{d_0^2}{\dt^2}(t-\dt)t \right]
\exp\left( - \frac{\dt}{4t (\dt - t) } r'^2\right) dr' \\
& = & -\frac{e^{- \frac{d_0^2}{\dt}}}
{(\dt - t)^{\frac{3}{2}}\sqrt{t}}
\SqrBrac{
\Brac{\frac{(\dt - t)t}{\dt}}^{\frac{3}{2}}4\sqrt{\pi}
- 
\frac{d_0^2}{\dt^2}(t-\dt)t 
\Brac{\frac{(\dt - t)t}{\dt}}^{\frac{1}{2}}
2\sqrt{\pi}}\\
& = & 
2\sqrt{\pi}
\exp \left( - \frac{d_0^2}{4\dt} \right) 
\left( \frac{d_0^2t}{(\dt)^{\frac{5}{2}}} 
- \frac{2t}{(\dt)^{\frac{3}{2}}} \right)
\end{eqnarray*}
which is (\ref{second integral}).
\end{proof}

Now we proceed to derive \eqref{Tricky Integral 1} and \eqref{Tricky Integral 2}, which follow almost immediately from the above computations.

\noindent
\Bf{Proof of (\ref{Tricky Integral 1})} -- the inner integral of $J_{111}$. 
Since 
\Beqn\label{d0-replace}
\Sigma_{s,\tau}(p) = s + r\left( \Sigma_{s,\tau}(p), \tau \right)
\N(s,\tau),
\,\,\,\,\,\,\Text{and}\,\,\,\,\,\,
\tilde{s} = s + r'\N(s,\tau), 
\Eeqn
we have that 
$$
| \Sigma_{s,\tau}(p) - \tilde{s}|^2 
= | r( \Sigma_{s,\tau}(p) , \tau) - r' |^2.
$$
Thus \eqref{Tricky Integral 1} can be re-written as,
\begin{equation*}
\Lover{(\dt - \tau)^{\frac{1}{2}}}
\int_{-\infty}^{\infty} 
\exp\left(
-\frac{| r( \Sigma_{s,\tau}(p), \tau) - r' |^2}{4(\dt - \tau)}
\right)
\frac{r'}{\sqrt{\tau}}\exp\left(-\frac{r'^2}{4\tau}\right) dr'.
\end{equation*}
Replacing $r(\Sigma_{s,\tau}(p), \tau)$ by $d_0$, 
this integral is identical to (\ref{first integral}).

\noindent
\Bf{Proof of (\ref{Tricky Integral 2})} -- 
the inner integral appearing in $J_{112}$.  
By \eqref{d0-replace} again, we have
\begin{eqnarray*}
\langle \Sigma_{s,\tau}(p) - \tilde{s}, \E_j \rangle 
&=& \langle \N(s,\tau) , \E_j \rangle 
\big( r(\Sigma_{s,\tau}(p) , \tau) - r' \big).
\end{eqnarray*}
Replacing once again $r(\Sigma_{s,\tau} (p) , \tau)$ by $d_0$,
and substituting this back into \eqref{Tricky Integral 2}, we get,
\begin{equation*}
\frac{\langle \N(s, \tau) , \E_j \rangle}
{(\dt - \tau)^{\frac{1}{2}}}
\int_{-\infty}^{\infty} \frac{d_0 - r'}{(\dt - \tau) } 
\exp\left(-\frac{|d_0 - r' |^2}{4(\dt - \tau)}\right)  
\frac{r'}{\sqrt{\tau}}\exp\left(-\frac{r'^2}{4\tau}\right) dr'
\end{equation*}
which is precisely $\langle \N(s,\tau) , \E_j \rangle$ times the integral
(\ref{second integral}).

\subsection{Proof of Lemma \ref{curvature bound}: 
Bound on the Curvature growth of MMC}\label{ProofCurvatureBound}
What follow in this and next sections are by now classical computations 
about MMC. The readers can refer to \cite{ecker2004regularity} for general 
exposition. For simplicity, we will drop the subscript $M_t$ and write
$A = A_{M_t},\,\,
H = H_{M_t},\,\,
\nabla = \nabla_{M_t}$, and
$\Lap_{M_t} = \Lap$.

By \cite[Corollary 3.5]{huisken1984flow}, 
the quantity $|A|^2$ (defined in \eqref{A2}) satisfies the following equation,
\begin{align*}
\partial_t |A|^2 &= \Delta|A|^2 - 2| \nabla A |^2 + 2 |A|^4.
\end{align*}
Let $x_t$ be the point at which $|A|^2$ attains its maximum.
Then at $(x_t, t)$, we have
\begin{align*}
\partial_t  |A|^2 &\leq 2 |A|^4.
\end{align*}
Hence setting $f(t) = \max_{M_t} |A|^2$, we have the following 
differential inequality:
\begin{align*}
\dfrac{\partial_t f }{f^2} &\leq 2.
\end{align*}
Integrating both sides in $t$ gives
\[
f(t) \leq \dfrac{f(0)}{1-2f(0)t} \leq f(0)\big(1 + C f(0) t \big).
\]
Note that $C$ may be chosen independent of $\|A_0\|$ and $t$ so long as 
$\| A_0 \| \big(\dt |\log \dt| \big)^{\frac{1}{4}} \leq 1$.
Taking square roots we then arrive at the desired result,
\begin{align*}
\|A_t\| \leq \|A_0\| \big( 1 + C \|A_0\|^2 t ).
\end{align*}

\subsection{Proof of Lemma \ref{regularity of curvature}: 
Higher Order Regularity of MMC}\label{ProofCurvatureRegularity}
We follow the arguments given in \cite{ecker1991interior}, 
where a similar bound was proven for the curvature.  
First note that the two estimates are equivalent, 
since by \cite[Lemma 3.3]{huisken1984flow}, we have
$\partial_t \N = \nabla H$.  

First, we quote the following equations given in 
\cite{ecker1991interior}
\begin{align*}
(\partial_t - \Delta) |A|^2 &= -2|\nabla A|^2 + 2|A|^4, \\
(\partial_t - \Delta) | \nabla A |^2 &\leq -2| \nabla^2 A |^2 + C |A^2| |\nabla A|^2.
\end{align*}

Introduce $\psi(t) = \frac{R^2t}{R^2 + t}$ for some $R > 0$. 
Then we compute:
\begin{eqnarray}
(\partial_t - \Delta)\big(\psi | \nabla A |^2\big) 
&\leq& -2 \psi | \nabla^2 A |^2 + C \psi |A^2| |\nabla A|^2 + |\nabla A |^2 \frac{d}{dt}\psi \nonumber\\
&\leq& 
-2 \psi | \nabla^2 A |^2 +  |\nabla A|^2 \big( 1 +  C\psi |A|^2).
\end{eqnarray}

Next define $f = \psi |\nabla A|^2\big( \Lambda_0 + |A|^2)$, 
where $\Lambda_0$ is a constant to be chosen later.  
Then we compute:
\begin{eqnarray}
(\partial_t - \Delta)f 
&\leq& 
\big( |A|^2 + \Lambda_0 \big ) \big \{ |\nabla A |^2 ( 1 + C |A|^2 \psi) - 2 \psi |\nabla^2 A|^2 \big \} \nonumber\\
& & 
+ \psi | \nabla A|^2 \big\{ - 2 | \nabla A|^2 + 2 |A|^4 \big\} - 2 \psi \nabla | \nabla A|^2 \cdot \nabla |A|^2.\label{grad.A.eqn}
\end{eqnarray}
We estimate the last term as follows:
\begin{eqnarray*}
- 2 \psi \nabla | \nabla A|^2 \cdot \nabla |A|^2 
&\leq& 8 \psi |\nabla A| |A| |\nabla |\nabla A|| | \nabla |A| |\\
&\leq & 8 \psi |\nabla A|^2 |A| |\nabla^2 A|
\,\,\,\,\,\,\,\,\, \text{(by Kato's inequality)}\\
& \leq & 2 \psi | \nabla^2 A |^2 \big( |A|^2 + \Lambda_0 \big) + \dfrac{8 \psi |A|^2 |\nabla A|^4}{|A|^2 + \Lambda_0} 
\,\,\,\,\,\,( ab \leq a^2 \epsilon + \frac{b^2}{4 \epsilon}).
\end{eqnarray*}
Substituting the above into \eqref{grad.A.eqn}, we get:
\begin{align*}
(\partial_t - \Delta)f &\leq \big( |A|^2 + \Lambda_0 \big) \big( 1 + C |A|^4 \psi \big) | \nabla A |^2 + 2\psi | \nabla A |^2 |A|^4 - 2 \psi |\nabla A |^2 + \dfrac{8 \psi |A|^2 				|\nabla A|^4}{|A|^2 + \Lambda_0} \\
				& = - \Big\{ 2 - \dfrac{8 |A|^2 }{|A|^2 + \Lambda_0} \Big\} \psi | \nabla A |^4 + \Big\{ \big( |A|^2 + \Lambda_0\big) \big( 1 + C |					A|^2 \psi \big) + 2 \psi |A|^4 \Big\} | \nabla A |^2 \\
				&= -\Big\{ 2 - \dfrac{8  |A|^2 }{|A|^2 + \Lambda_0} \Big\} \dfrac{f^2}{\psi \big( |A|^2 + \Lambda_0 \big)^2} + \psi^{-1}f \Big\{ \big( 1 + C |A|^2 \psi \big) + \dfrac{2 \psi |A|^4}{\lambda_0 + |A|^2} \Big\} \\
				&= - \psi^{-1} \big( \delta f^2 - \bar{K} f \big),
\end{align*}
where $\delta = \Big\{ 2 - \dfrac{8  |A|^2}{|A|^2 + \Lambda_0} \Big\} \big( |A|^2 + \Lambda_0 \big)^{-2}$ and $\bar{K} = \big( 1 + C |A|^2 \psi \big) + \dfrac{2 \psi |A|^4}{\Lambda_0 + |A|^2}$.

Next compute estimates for $f \eta$, where $\eta = (R^2 - ( |x - x_0|^2 + 2nt ) \big)^2 \equiv (R^2 - r(x,t) )^2  $ acts as a ``localization'' function:
\begin{align*}
(\partial_t - \Delta)f \eta &\leq - \eta \psi^{-1} \big( \delta f^2 - \bar{K} f \big) + f (\partial_t - \Delta) \eta - 2 \nabla f \cdot \nabla \eta \\
				    & \leq - \eta \psi^{-1} \big( \delta f^2 - \bar{K} f \big) + 4(R^2 - r(x,t))|x-x_0|^2 f - 2 \nabla (f \eta) \cdot \frac{\nabla \eta}{\eta} \\
				    &\leq - \eta \psi^{-1} \big( \delta f^2 - \bar{K} f \big) + 4R^4 f - 2 \nabla (f \eta) \cdot \frac{\nabla \eta}{\eta},
\end{align*}
the last inequality holds over the set $M_t \cap \{r< R^2\}$.

Now consider 
\MathSty{m(\dt) = \sup_{0 \leq t \leq \dt}
\CurBrac{\sup_{x \in \{M_t | r(x,t) \leq R^2\}} f \eta}}.  
Notice that $\psi \equiv 0$ 
at $t = 0$ and hence $f \eta \equiv 0$ at $t = 0$ also.  
Suppose $m(\dt)$ is attained at some point $(\tilde{x}, \tilde{t} )$ with $\tilde{t} \geq 0$.  At this point, $ (\partial_t - \Delta)f \eta \geq 0$.  Thus we can have the following sequence of implications:
\begin{eqnarray*}
0 &\leq& 
- \eta \psi^{-1} \big( \delta f^2 - \bar{K} f \big) + 4R^4f\\
\psi^{-1} \delta f^2 \eta &\leq& \psi^{-1} \bar{K}  f \eta + 4R^4 f\\
f^2 \eta 
&\leq& \frac{1}{\delta} \Big( \bar{K} f \eta + 4R^4 \psi f \Big)\\
f^2 \eta^2 
& \leq & \frac{1}{\delta} \Big( \bar{K} \eta + 4R^4 \psi \Big) f \eta.
 \end{eqnarray*}
Applying Young's inequality to the last line of the above 
leads to that at $(\tilde{x}, \tilde{t})$,
\begin{eqnarray*}
f^2 \eta^2 \leq \frac{1}{\delta^2} \Big( \bar{K} \eta + 4R^4 \psi \Big)^2.
\end{eqnarray*}
Note that $\eta \leq R^4$.  Substituting this estimate in above to obtain:
\[
m(\dt)^2 \leq \frac{R^8}{\delta^2}\Big( \bar{K} + 4 t \Big)^2,
\,\,\,\,\,\,\text{i.e.}\,\,\,\,\,\,
m(\dt) \leq \frac{R^4}{\delta}\Big( \bar{K} + 4 t \Big).
\]
This means that for all $(x,t) \in \{ M_t | r(x,t) \leq \theta R^2 \}$ (for some $0 < \theta < 1$) we have,
\begin{align*}
\psi | \nabla A|^2 &\leq (1 - \theta)^{-2} \delta^{-1} \big( \Lambda_0 + |A|^2 \big)^{-1} \big( \bar{K} + 4t \big).
\end{align*}

Now set $\Lambda_0 = 8 c_0$.  We then have the following 
estimates/equalities:
\begin{enumerate}
\item \MathSty{
\delta^{-1} \big( \Lambda_0 + |A|^2 \big)^{-1} 
= \big(|A|^2 + \Lambda_0) 
\Big( 2 - \frac{8 |A|^2}{|A|^2 + \Lambda_0} \Big)^{-1}
\leq 6 c_0;
}
\item
\MathSty{
\bar{K}  
= \big( 1 + C |A|^2 \psi \big) + \dfrac{2 \psi |A|^4}{\Lambda_0 + |A|^2} 
\leq C + \frac{2}{9}\psi |A|^2
\leq C 
}
which is true by our assumption on the magnitude of
the curvature in relation to the time step;
\item
\MathSty{
\psi^{-1} = \frac{R^2 + t}{R^2t} = \frac{1}{t} + \frac{1}{R^2}.
}
\end{enumerate}

Utilizing the above estimates we obtain for all $(x,t) \in \{ M_t | r(x,t) \leq \theta R^2 \}$,
\begin{align*}
| \nabla A |^2 &\leq C \cdot c_0 \big( \frac{1}{t} + \frac{1}{R^2} \big).
\end{align*}
Finally let $R \rightarrow \infty$ to obtain the estimate for all of $M_t$.  Taking square roots we get,
\begin{align*}
|\nabla A| &\leq	C \dfrac{c_0}{\sqrt{t}}.
\end{align*}

\subsection{Proof of Lemma \ref{BallLemma}: the Ball Lemma}
\label{ProofBallLemma}
The Ball Lemma \ref{BallLemma} states that a ball of uniform radius 
may be placed in a tangential manner in the interior and exterior of
the numerical manifolds
$\Gamma_{k\dt}$ for 
$0\leq k \leq \left\lfloor\frac{T}{\dt}\right\rfloor$
-- see Definition \ref{BallProperty} to recall the ball property.
This result is used to show that the numerical manifold 
converges to an \Em{embedded} manifold. 

The intuitive reason behind the Ball Lemma is quite simple: 
once a ball of radius $r_k$ 
can be put inside and outside of $\Gamma_{k\dt}$, then by regularity 
of the motion, the same ball but with a smaller radius 
$r_{k+1} = r_k - \CO(\dt)$ 
can still be put inside and outside $\Gamma_{(k+1)\dt}$. 
Such a statement can then be iterated over multiple time steps.
Though this also follows from comparison principle for MMC, we choose
the following route of proof for the sake of its more applicability
in other geometric flows.

The rigorous proof makes use of the
intrinsic distance $d_0$ between two points on the numerical 
manifold $\Gamma_0\subset\R^{n+1}$. This is defined as the length 
of the shortest curve on $\Gamma_0$ joining two points 
$p,\,q\in\Gamma_0$: 
\Beqn\label{DistOnManifold}
d_0(p, q) = 
\min\CurBrac{
\int_0^1 \left|\dot{c}(s)\right| ds:
c:[0,1]\MapTo\Gamma_0,\,\,\,
c(0) = p,\,\,c(1)=q
}.
\Eeqn
Analogously, we will use $d_k$ to denote the intrinsic distance on the 
numerical manifold $\Gamma_{k\dt}$. 

We first present a preliminary result that bounds the amount by which 
the intrinsic distance may change over a short time on a manifold moving 
by its mean curvature. Let $M_0\subset\R^{n+1}$ be a smooth, compact 
$n$-dimensional manifold. Recall the function 
$F: M_0 \times \R^+\MapTo\R^{n+1}$ which solve \eqref{MMC.eqn} so that
$M_t = F(M_0, t)$ evolves according to MMC.
\begin{thm4}\label{IntrinsicDistance}
Let $A_0$ be the Weingarten Map of $M_0$ satisfying 
$\|A_0\| \left (\dt | \log \dt | \right)^{\frac{1}{4}} \leq 1.$
Let $c: [0,1]\MapTo M_0$ be curve on $M_0$.  
Then the following estimate holds.
\begin{align*}
\big | \ell \left( F(c, t) \right) - \ell (c) \big | 
\lesssim \ell(c) \| A_0 \|^2 t,
\end{align*}
where 
\MathSty{\ell(c) = \int \left|\dot{c}(s)\right|\,ds} and
\MathSty{\ell(F(c,t)) = \int \left|DF(c(s),t)\dot{c}(s)\right|\,ds} 
are the lengths of $c$ and $F(c,t)$ repsectively.
\end{thm4}
\begin{proof}
Without loss of generality, suppose $c$ is parametrized by arc length over the interval $[0, \ell(c)] \subset \mathbb{R}$. 
We quickly describe some notation.  
We use $\N$ to denote the unit normal to $M_t$, 
$H$ its mean curvature, 
$\nabla$ denotes the tangential gradient over $M_t$, 
and $\dot{}$ denotes differentiation with respect to the arc length variable.  
Furthermore $\big\langle \cdot, \cdot \big\rangle$ denotes 
the inner product in $\mathbb{R}^{n+1}$.

We first claim that,
\begin{align}\label{shortKinks_claim1}
\left | \frac{d}{dt} \left | \frac{d}{ds} F(c, t) \right |^2 \right | & \lesssim \| A_0 \|^2 \left | \frac{d}{ds} F(c, t) \right |.
\end{align}
Toward this end, consider:
\begin{align*}
\frac{d}{dt} \left | \frac{d}{ds} F(c, t) \right |^2
& = 
{2\frac{d}{dt}\Big\langle 
\nabla F(c(s), t)\dot{c}(s), \nabla F(c(s), t)\dot{c}(s)
\Big\rangle}\\
& = 
{4\Big\langle 
\nabla (\partial_t F)(c(s), t))\dot{c}(s), \nabla F(c(s), t)\dot{c}(s)
\Big\rangle}\\
& = 
{-4\Big\langle 
\nabla (H\N(c(s), t))\dot{c}(s), \nabla F(c(s), t)\dot{c}(s)
\Big\rangle}\\
& = {-4 \left \{ 
\Big\langle \N \otimes \nabla H \dot{c}, \nabla F\dot{c} \Big\rangle 
+ H \Big\langle \nabla \N \dot{c}, \nabla F\dot{c} \Big\rangle 
\right\}.}
\end{align*}
For the first term on the right hand side, we compute:
further compute:
\begin{align*}
\langle \N \otimes \nabla H \dot{c}, \nabla F\dot{c} \rangle
	= \sum_{j=1}^{n+1} \sum_{k=1}^{n+1} \partial_k F^{(j)} \dot{c}^{(k)} 
\N^{(j)} \partial_k H \dot{c}^{(k)}
	= \sum_{k=1}^{n+1} \partial_k H (\dot{c}^{(k)})^2 \langle \partial_k F, \N \rangle
	= 0.
\end{align*}
while for the second term, we have:
\begin{align*}
\left | H \langle \nabla \N \dot{c}, \nabla F\dot{c} \rangle \right |
	\lesssim \| A_0 \|  |\nabla F \dot{c}  | | \nabla \N \dot{c}  |
	 = \|A_0 \| \left | \frac{d}{ds} F(c, t) \right | \left | \frac{d}{ds} \N(c, t) \right |
	 \lesssim \|A_0\|^2 \left | \frac{d}{ds} F(c, t) \right |.
\end{align*}
Combining the above computations leads to (\ref{shortKinks_claim1})

Finally,
\begin{align*}
\left | \ell(F(c,t) - \ell(c) )\right |
	&= \left |  \int_0^t \frac{d}{d\tau} \int_0^{\ell(c)} 
	\left|\frac{d}{ds} F(c,\tau)\right| ds d\tau \right | \\
	& = \left |  \int_0^t \int_0^{\ell(c)}  \frac{d}{d\tau}\left( \left|\frac{d}{ds} F(c,\tau) \right|^2 \right)  
	\Lover{2}\left|\frac{d}{ds} F(c,\tau) \right| ^{-1}  ds d\tau \right | \\
	& \lesssim \ell(c) \|A_0\|^2 t
\end{align*}
which concludes the proof of the Lemma.
\end{proof}

With the above, the Ball Lemma \ref{BallLemma} then
follows from the next two claims. We first let
$r_0:= \frac{1}{C_0 + 1}$,
$B^\text{int}_{r_0, p}$ and $B^\text{ext}_{r_0, p}$ be
the interior and exterior balls with radius $r_0$
which are tangent to $\Gamma_{k \dt}$ at $p$ 
(recall Definition \ref{BallProperty}),
and $r_* := \frac{C_*}{C_0}$ where $C_*$ is from \eqref{RadiusGraph}.

\begin{description}
\item[Claim 1.]
Fix a $p \in \Gamma_{k \dt}$. Then
\Beqn\label{Claim1}
\Big\{q\in\Gamma_{k\dt}:  d_k(p, q) \leq r_* \Big\} \cap 
\Big( B^\text{int}_{r_0,p} \cup B^\text{ext}_{r_0, p} \Big) = \{ p \}.
\Eeqn
\end{description}

\begin{proof}
It follows from the remark after \eqref{RadiusGraph} that 
the connected component of 
$\Gamma_{kh}\Intersect B_{r_*}(p)$ containing $p$ can be written 
as the graph of a function $f$ over the tangent plane $T_p(\Gamma_{kh})$.
By the curvature bound of $\Gamma_{kh}$, we infer that 
\[
\Big\{(x,f(x)): x\in T_p(\Gamma_{kh}),\,\,\abs{x} \leq r_*\Big\}
\Intersect
\Big( B^\text{int}_{r_0, p} \cup B^\text{ext}_{r_0,p} \Big)
= \CurBrac{p}.
\]
Hence any curve in $\Gamma_{kh}$ that joins $p$ and any other
$q \in B^\text{int}_{r_0,p} \cup B^\text{ext}_{r_0,p}$ must have length
at least $r_*$. Hence the claim follows.
\end{proof}

\begin{description}
\item[Claim 2.] 
Define 
$m_k := \min\big\{| p - p|:\,\,
p, q \in \Gamma_{k\dt},\,\,
d_k(p, q) \geq r_*\big\}.
$
Then
\begin{align}\label{step2_claim}
m_{k} \geq m_{k-1} - C(n,C_0)\dt,
\end{align}
for some fixed constant $C(n,C_0)$ depending only on the spatial dimension $n$ and the 
curvature bound $C_0$.
\end{description}

\begin{proof}
We will compare the intrinsic distance between 
$\Gamma_{k \dt}$ and $\Gamma_{(k-1) \dt}$. Recall that 
$\Gamma_{k\dt}$ is obtained as a graph over 
$M_{k\dt}:=M(\Gamma_{(k-1)\dt}, \dt)$, the solution at time $\dt$ of MMC 
with initial data $\Gamma_{(k-1)\dt}$. More precisely, 
\begin{enumerate}
\item
using the function $F$ in \eqref{MMC.eqn},
we have $M_{k\dt}=F(\Gamma_{(k-1)\dt}, \dt)$ and,
\item
there is a function $h: M_{k\dt} \MapTo\R$ such that
$\Gamma_{k\dt} = M_{k\dt} + h\N_{M_{k\dt}}$.
\end{enumerate}
Note that both the transformations from $\Gamma_{(k-1)\dt}$ to $M_{k\dt}$
and from $M_{k\dt}$ to $\Gamma_{k\dt}$ are diffeomorphisms.
Furthermore, by the consistency and stability Theorems
\ref{consistency} and \ref{Finite Time Stability}, we have 
$\Norm{h}_{L^\infty} \lesssim \Norm{\CB_0}^2\dt^\frac{3}{2}$.

Now let $p, q \in \Gamma_{k\dt}$ be such that $d_k(p,q) \geq r_*$. 
We can find $p_1, q_1\in \Gamma_{(k-1)\dt}$ and 
$p_2, q_2\in M_{k\dt}$ satisfying
\[
p_2 = F(p_1,\dt),\,\,
p=p_2 + h(p_2)\N_{M_{k\dt}}(p_2),\,\,
q_2 = F(q_1,\dt),\,\,
q=q_2 + h(q_2)\N_{M_{k\dt}}(q_2).
\]
Then we have
\Beqn\label{distance.lower.bound}
\left |p - q \right | 
\geq
\left |p_2 - q_2\right| - C(n) C_0^2(\dt )^{\frac{3}{2}}
\geq 
\left |p_1 - q_1 \right | - C_1(n, C_0) \dt
\Eeqn
where the first inequality is due to the consistency of the scheme while
the second is due to the regularity of MMC.

If $d_{k-1}(p_1, q_1) \geq r_*$, then 
we automatically have $|p_1-q_1| \geq m_{k-1}$ and hence
we are done.

Now suppose $d_{k-1}(p_1, q_1) < r_*$.
Making use of Lemma \ref{IntrinsicDistance}, we have
\[
\Big|d_k(p,q) - d_{(k-1)}(p_1, q_1)\Big|
\lesssim
d_{(k-1)}(p_1, q_1)\Norm{\CB_0}^2\dt.
\]
from which it follows that
\[
d_{k-1}(p_1,q_1)
\geq 
\frac{d_k(p,q)}{1+\Norm{\CB_0}^2\dt}
\geq
\frac{r_*}{1+\Norm{\CB_0}^2\dt}
\geq
r_* - C_2(n,C_0) \dt
\]
where $C_2(n,C_0)$ is another constant depending on the 
spatial dimension $n$ and the curvature bound $C_0$.
Next, by extending the geodesic curve which joins $p_1$ and $q_1$, 
we can find $\hat{p}_1, \hat{q}_1 \in \Gamma_{(k-1)\dt}$ such that
$|p_1 - \hat{p}_1|,\,|q_1 - \hat{q}_1| \leq C_2(n, C_0)\dt$ 
and $d_{k-1}(\hat{p}_1, \hat{q}_1) \geq r_*$.
Hence
\Beqn
|p_1 - q_1| 
\geq 
|\hat{p}_1 - \hat{q}_1| - C_2(n, C_0)\dt
\geq 
m_{k-1} - C_2(n, C_0)\dt.
\Eeqn
Now, going back to \eqref{distance.lower.bound}, we have
\[
|p-q| 
\geq |p_1 - q_1| - C_1(n, C_0)\dt
\geq m_{k-1} - (C_1(n, C_0)+C_2(n, C_0))\dt
\]
from which \eqref{step2_claim} follows.
\end{proof}

\Bf{Acknowledgment.}
The authors would like to thank Bob Jerrard for pointing us to this direction,
Selim Esedoglu for useful discussion, Tim Laux for reading and suggestions
concerning the manuscript, and the Purdue Research Foundation for 
partial support.

\bibliographystyle{plain}
\bibliography{reference}

\begin{thebibliography}{10}

\bibitem{BG}
Guy Barles and Christine Georgelin.
\newblock A simple proof of convergence for an approximation scheme for
  computing motions by mean curvature.
\newblock {\em SIAM Journal on Numerical Analysis}, 32(2):484--500, 1995.

\bibitem{SouganidisBarles}
Guy Barles and Panagiotis~E Souganidis.
\newblock A new approach to front propagation problems: theory and
  applications.
\newblock {\em Archive for rational mechanics and analysis}, 141(3):237--296,
  1998.

\bibitem{Chambolle}
Eric Bonnetier, Elie Bretin, and Antonin Chambolle.
\newblock Consistency result for a non monotone scheme for anisotropic mean
  curvature flow.
\newblock {\em Interfaces and Free Boundaries}, 14(1):1--35, 2012.

\bibitem{Brezis}
Haim Brezis.
\newblock {\em Functional analysis, Sobolev spaces and partial differential
  equations}.
\newblock Springer Science \& Business Media, 2010.

\bibitem{bronsard1991motion}
Lia Bronsard and Robert~V Kohn.
\newblock Motion by mean curvature as the singular limit of ginzburg-landau
  dynamics.
\newblock {\em Journal of differential equations}, 90(2):211--237, 1991.

\bibitem{chen1992generation}
Xinfu Chen.
\newblock Generation and propagation of interfaces for reaction-diffusion
  equations.
\newblock {\em Journal of Differential equations}, 96(1):116--141, 1992.

\bibitem{chen1991uniqueness}
Yun~Gang Chen, Yoshikazu Giga, ShunÕichi Goto, et~al.
\newblock Uniqueness and existence of viscosity solutions of generalized mean
  curvature flow equations.
\newblock {\em J. Differential Geom}, 33(3):749--786, 1991.

\bibitem{cheung2004some}
Wing-Sum Cheung.
\newblock Some discrete nonlinear inequalities and applications to boundary
  value problems for difference equations.
\newblock {\em Journal of Difference Equations and Applications},
  10(2):213--223, 2004.

\bibitem{de1995geometrical}
Piero De~Mottoni and Michelle Schatzman.
\newblock Geometrical evolution of developed interfaces.
\newblock {\em Transactions of the American Mathematical Society},
  347(5):1533--1589, 1995.

\bibitem{ecker2004regularity}
Klaus Ecker.
\newblock {\em Regularity theory for mean curvature flow}.
\newblock Springer, 2004.

\bibitem{ecker1991interior}
Klaus Ecker and Gerhard Huisken.
\newblock Interior estimates for hypersurfaces moving by mean curvature.
\newblock {\em Inventiones mathematicae}, 105(1):547--569, 1991.

\bibitem{elsey2016threshold}
Matt Elsey and Selim Esedoglu.
\newblock Threshold dynamics for anisotropic surface energies.
\newblock Technical report, preprint, 2016, 2016.

\bibitem{EseOtto}
Selim Esedoglu and Felix Otto.
\newblock Threshold dynamics for networks with arbitrary surface tensions.
\newblock {\em Communications on Pure and Applied Mathematics}, 68(5):808--864,
  2015.

\bibitem{EseRuuthTsai}
Selim Esedoglu, Steven Ruuth, and Richard Tsai.
\newblock Threshold dynamics for high order geometric motions.
\newblock {\em Interfaces and Free Boundaries}, 10(3):263--282, 2008.

\bibitem{Evans}
Lawrence~C Evans.
\newblock Convergence of an algorithm for mean curvature motion.
\newblock {\em Indiana University Mathematics Journal}, 42(2):533--557, 1993.

\bibitem{evans1992phase}
Lawrence~C Evans, H~Mete Soner, and Panagiotis~E Souganidis.
\newblock Phase transitions and generalized motion by mean curvature.
\newblock {\em Communications on Pure and Applied Mathematics},
  45(9):1097--1123, 1992.

\bibitem{evans1991motion}
Lawrence~C Evans, Joel Spruck, et~al.
\newblock Motion of level sets by mean curvature i.
\newblock {\em J. Diff. Geom}, 33(3):635--681, 1991.

\bibitem{evans-gariepy}
Lawrence~Craig Evans and Ronald~F Gariepy.
\newblock {\em Measure theory and fine properties of functions}.
\newblock CRC press, 2015.

\bibitem{Giusti}
E~Giusti.
\newblock {\em Minimal surfaces and functions of bounded variation}.
\newblock Springer Science \& Business Media, 1984.

\bibitem{huisken1984flow}
Gerhard Huisken et~al.
\newblock Flow by mean curvature of convex surfaces into spheres.
\newblock {\em Journal of Differential Geometry}, 20(1):237--266, 1984.

\bibitem{ilmanen1993convergence}
Tom Ilmanen et~al.
\newblock Convergence of the allen-cahn equation to brakkeÕs motion by mean
  curvature.
\newblock {\em J. Differential Geom}, 38(2):417--461, 1993.

\bibitem{ishii1995generalization}
Hitoshi Ishii.
\newblock A generalization of the bence, merriman and osher algorithm for
  motion by mean curvature.
\newblock {\em Curvature flows and related topics (Levico, 1994)}, 5:111--127,
  1995.

\bibitem{ishii1999threshold}
Hitoshi Ishii, Gabriel~E Pires, and Panagiotis~E Souganidis.
\newblock Threshold dynamics type approximation schemes for propagating fronts.
\newblock {\em Journal of the Mathematical Society of Japan}, 51(2):267--308,
  1999.

\bibitem{jost2002partial}
Jurgen Jost.
\newblock {\em Partial Differential Equations}.
\newblock Graduate texts in mathematics. Springer, 2002.

\bibitem{laux2016convergence}
Tim Laux and Felix Otto.
\newblock Convergence of the thresholding scheme for multi-phase mean-curvature
  flow.
\newblock {\em Calculus of Variations and Partial Differential Equations},
  55(5):129, 2016.

\bibitem{leoni2001convergence}
Fabiana Leoni.
\newblock Convergence of an approximation scheme for curvature-dependent
  motions of sets.
\newblock {\em SIAM journal on numerical analysis}, 39(4):1115--1131, 2001.

\bibitem{luckhaus}
Stephan Luckhaus and Thomas Sturzenhecker.
\newblock Implicit time discretization for the mean curvature flow equation.
\newblock {\em Calculus of variations and partial differential equations},
  3(2):253--271, 1995.

\bibitem{BMO}
Barry Merriman, James~Kenyard Bence, and Stanley Osher.
\newblock {\em Diffusion generated motion by mean curvature}.
\newblock CAM Report, 92-18, Department of Mathematics, University of
  California, Los Angeles, 1992.

\bibitem{OsherRuuthXin}
Steven~J Ruuth, Barry Merriman, Jack Xin, and Stanley Osher.
\newblock Diffusion-generated motion by mean curvature for filaments.
\newblock {\em Journal of Nonlinear Science}, 11(6):473--493, 2001.

\bibitem{soner1}
Halil~Mete Soner.
\newblock Ginzburg-landau equation and motion by mean curvature, i:
  convergence.
\newblock {\em The Journal of Geometric Analysis}, 7(3):437--475, 1997.

\bibitem{soner2}
Halil~Mete Soner.
\newblock Ginzburg-landau equation and motion by mean curvature, ii:
  Development of the initial interface.
\newblock {\em The Journal of Geometric Analysis}, 7(3):477--491, 1997.

\end{thebibliography}
\end{document}